\documentclass[11pt]{article}
\usepackage{amssymb}
\usepackage{latexsym}
\usepackage{amsfonts}
\usepackage{amsmath}
\newtheorem{thrm}{Theorem}[section]
\newtheorem{lemma}[thrm]{Lemma}
\newtheorem{prop}[thrm]{Proposition}

\newtheorem{cor}[thrm]{Corollary}
\newtheorem{remark}[thrm]{Remark}

\numberwithin{equation}{section}
\usepackage[dvips]{color}
\usepackage{mathtools}
\mathtoolsset{showonlyrefs}
\usepackage{mathrsfs}
\usepackage{comment}
\usepackage{hyperref}
\usepackage{xcolor}

\def\P{\mathbb{P} }

\def\N{\mathbb{N} }

\makeatletter
\begin{document}
	\allowdisplaybreaks

	\title{\Large \bf Law of iterated logarithm for
	supercritical non-symmetric branching Markov process}
	\author{ \bf   Haojie Hou\footnote{The research of this author is supported by the China Postdoctoral Science Foundation (No.
   2024M764112)\hspace{1mm} } \hspace{1mm}\hspace{1mm}
   Yan-Xia Ren\footnote{The research of this author is supported by NNSFC (Grant No. 12231002) and the Fundamental Research Funds for the Central Universities, Peking University LMEQF.\hspace{1mm} } \hspace{1mm}\hspace{1mm}
	Renming Song\thanks{Research supported in part by a grant from the Simons
	Foundation
	(\#960480, Renming Song).}
	\hspace{1mm}\hspace{1mm}
	\hspace{1mm} }
	\date{}
	\maketitle
	
\begin{abstract}
Let $\{(X_t)_{t\geq 0}, \mathbb{P}_{\delta_x}, x\in E\}$ be a
supercritical branching Markov process (which is not necessary symmetric) on a locally compact metric measure space
$(E,\mu)$ with spatially dependent local branching mechanism. Under some assumptions on the semigroup of the spatial motion,
we first prove law of  iterated logarithm type results for $\langle f, X_t\rangle$ under the second
moment condition on the branching mechanism,
where $f$ is a linear combination of eigenfunctions of the mean semigroup
$\{T_t, t\geq0\}$ of $X$. Then we prove
law of  iterated logarithm type results for $\langle f, X_t\rangle$ under the fourth moment condition, where $f$ belongs to a
larger class of functions.
	\end{abstract}
	
	\medskip
	
	
	\medskip

	
\section{Introduction}
Let $\{Z_n: n\in \mathbb{N}\}$ be a supercritical Galton-Watson process with $Z_0=1$ and
$\mathbb{E}(Z_1)=m\in (1,\infty)$.
It is well-known that, under the assumption $\mathbb{E}(Z_1^2)<\infty$, the process $W_n:= Z_n/m^n$ is a non-negative $L^2$ bounded martingale and thus converges almost surely and in $L^2(\mathbb{P})$ to a non-negative limit $W_\infty$.  Heyde \cite{Heyde1970, Heyde1971}
found the rate at which $W_n -W_\infty$ converges to 0: $m^{n/2}(W_n-W_\infty)$ converges in distribution to $\sqrt{W_\infty} \mathcal{N}(0, \sigma^2)$, where $\mathcal{N}(0, \sigma^2)$ is a normal random variable, independent of $W_\infty$, with variance $\sigma^2:= \frac{1}{m^2-m}(\mathbb{E}(Z_1^2)-m^2)$.  The fluctuation in the almost sure sense of $W_n- W_\infty$
 was established  by Heyde \cite{Heyde1971-2}. Under the assumption $\mathbb{E}(Z_1^3)<\infty$, Heyde \cite{Heyde1971-2} proved that,  on the event $\{W_\infty >0\}$, it holds almost surely that
\begin{align}\label{LIL-Heyde}
\limsup_{n\to\infty}/ \liminf_{n\to\infty}  \frac{m^{n/2}(W_n-W_\infty)}{\sqrt{2\log n}} =+/- \sqrt{\sigma^2W_\infty}.
\end{align}
Later, Heyde and Leslie \cite{HL1971} removed the assumption $\mathbb{E}(Z_1^3)<\infty$ and  proved \eqref{LIL-Heyde} under the second moment condition only.
Since $\lim_{n\to\infty} \frac{\log \log Z_n}{\log n}= 1$ almost surely on $\{W_\infty >0\}$,
it follows from \eqref{LIL-Heyde} that almost surely on $\{W_\infty >0\}$,
\begin{align}
	\limsup_{n\to\infty}/	\liminf_{n\to\infty}  \frac{m^{n/2}(W_n-W_\infty)}{\sqrt{2\log \log Z_n}} = +/- \sqrt{\sigma^2W_\infty}.
\end{align}
Therefore,
results like  \eqref{LIL-Heyde} are called ``laws of iterated logarithm"
(LIL) in the literature. See \cite[Remark 1.3]{IK} and \cite[Remark 2.4]{IKM}.

For supercritical (finite) multitype Galton-Watson processes $\{Z_n: n\in \mathbb{N}\}$,
 Kesten and Stigum  \cite{KS1, KS2} established central limit theorems  by using the Jordan canonical form of the expectation matrix $M$.
 Asmussen \cite{Asmussen1977} extended \eqref{LIL-Heyde}  to $Z_n\cdot a$, where $a$ is a vector satisfying certain conditions.
In the continuous time setting,
central limit type theorems were proved by Athreya \cite{Ath1, Ath2, Ath3} and an analog of \eqref{LIL-Heyde} was given in \cite[Theorem 2]{Asmussen1977}.
	
 There are also some LIL type theorems for more general branching processes.
 Gao and Hu \cite{GH2012} proved \eqref{LIL-Heyde} for  branching processes in random environment.
 For  branching random walks, Iksanov and Kabluchko \cite{IK} proved an  LIL type theorem for Biggins' martingale.
For general Crump-Mode-Jagers branching processes, Iksanov et al \cite{IKM} proved an  LIL type theorem for Nerman's martingale.
All known LIL type results for branching processes, including
 branching random walks and Crump-Mode-Jagers branching processes, are LIL for $L^2$ bounded martingales.
For some related results for $L^2$ bounded martingale in the general case, see
  \cite{Heyde1977, S2017}.

In this paper, we are  interested in
supercritical branching Markov
processes with spatially dependent (local)  branching mechanism.
We always assume that $E$ is a locally compact separable metric space and that
$\mu$ is a $\sigma$-finite Borel measure on $E$ with full support.
We assume that $\partial$ is a point not in $E$ and put
   $E_\partial := E \cup \{\partial \}$.
 Any function $f$ on $E$ is automatically  extended to $E_\partial$ by defining $f(\partial)=0$.
We assume  that $\xi =\{\xi_t , \mathbb{P}_x, x\in E\}$ is a Hunt process on $E$ and that $\zeta: = \inf\{ t >0: \xi_t = \partial\}$ is the lifetime of $\xi$. The semigroup of $\xi$ is denoted by $\{ P_t : t\geq 0\}$. Our  standing
assumption on  $\xi$ is that
there exists a family of continuous strictly positive functions $\{p_t(x,y): t > 0\}$ on $E\times E$ such that
$$P_tf(x) = \int_E p_t(x,y)f(y)\mu(\mathrm{d}y).$$
Let
\[
\widehat{P}_tf(x):= \int_E p_t(y,x) f(y) \mu(\mathrm{d}y)
\]
be the dual operator
of $P_t$.
We use  $\mathbb{C}$ to denote  the set of complex numbers.
Let $L^p(E,\mu; \mathbb{C}):=\big\{f: E\to \mathbb{C}:
\Vert f\Vert_p:= \big(\int_E|f(x)|^p\mu(\mathrm{d}x)\big)^{1/p}
<\infty \big\}$ and $L^p(E,\mu):= \left\{f\in L^p(E,\mu; \mathbb{C}):\ f \mbox{ is real} \right\}.$
 For any complex number $z$, we use $\mathfrak{R}(z)$ and $\mathfrak{I}(z)$ to denote the real and imaginary parts of $z$ respectively.
Our first assumption is as follows:
    \begin{itemize}
        \item [{\bf(H1)}]  {\bf(a)} For all $t>0$ and $x\in E$, $\int_E p_t(y,x)\mu(\mathrm{d}y)\leq 1$.

        {\bf(b)} For any $t>0$, both of the functions
        \[
        x\mapsto a_t(x):= \int_E p_t(x,y)^2\mu(\mathrm{d}y) \quad \mbox{and}\quad x\mapsto \widehat{a}_t(x):=\int_E p_t(y,x)^2\mu(\mathrm{d}y)
        \]
        are continuous in $E$  and belong to  $L^1(E,\mu)$.

        {\bf(c)} There exists $t_0 > 0$ such that  $a_{t_0}, \widehat{a}_{t_0}\in
         L^2(E,\mu).$
    \end{itemize}
    Note that, see \cite[Section 1.1]{RSZ2017}, {\bf(H1)}(c) is equivalent to: There exists $t_0>0$ such that $a_{t}, \widehat{a}_{t}\in L^2(E,\mu)$ for all $t\ge t_0$.

A branching Markov process can be described as follows:
 initially there is a particle located at $x\in  E$ and  it moves according to $\{\xi, \mathbb{P}_x\}$.  When the particle is at site $y$, the branching rate is given by $\beta(y)$, where $\beta$ is a non-negative Borel function, that is, each individual dies in $[t, t + \textup{d}t)$ with probability $\beta(\xi_t)\textup{d}t + o(\textup{d}t)$.  When an individual dies at $y\in E$, it splits into $k$ particles with probability $p_k(y)$.  Once an individual reaches $\partial$, it disappears from the system.  All the individuals, once born, evolve independently.

Our assumption on the branching particle system is as follow:
\begin{itemize}
    	\item [{\bf(H2)}] {\bf(a)} $\beta(x)$  is a non-negative bounded Borel function on $E$.
    	
    	{\bf(b)} $\{p_k(x):k=0,1,...\}$ satisfies
    	$$\sup_{x\in E} \sum_{k=0}^\infty k^2 p_k (x) < \infty.$$
\end{itemize}
 Let $\mathcal{M}_a(E)$ be the space of finite atomic measures on
$E$.
 For $t\ge 0$ and $B \in \mathcal{B}(E),$ let $X_t(B)$ denote the number of particles alive at time $t$ and located in $B$. Then $X=\{X_t: t\geq 0 \}$ is an $\mathcal{M}_a (E)$-valued Markov process. For any $x \in E$, we denote by $\mathbb{P}_{\delta_x}$ the law of $X$ with initial value $X_0=\delta_x$.
 For any Borel function $f$ on $E$ and $\nu\in \mathcal{M}_a(E)$, define $\langle f, \nu\rangle := \int_{E} f(y)\nu(\mathrm{d}y)$. For $f, g\in L^2(E, \mu; \mathbb{C})$, define
 $\langle f, g\rangle_\mu:= \int_E f(x) \overline{g(x)} \mu(\mathrm{d}x)$.
 For any non-negative
 bounded Borel function $f$ on $E$, let
 $$
 \omega(t,x) := \mathbb{E}_{\delta_x} \left(e^{-\langle f, X_t \rangle}\right),
 $$
 then it is well-known (for example, see \cite[Section 1.1]{GHK2022}) that $\omega(t,x)$ is the unique positive solution to the equation
 \begin{equation}
    	\omega(t,x) = \mathbb{E}_x\left(\int_0^t \psi\left(\xi_s, \omega(t-s, \xi_s) \right) \mathrm{d}s \right)+ \mathbb{E}_x\left(e^{-f(\xi_t)}\right),
 \end{equation}
 here $\psi(x,z) = \beta(x) \left(\sum_{k=0}^\infty p_k(x) z^k -z\right)$ if $x \in E, z\in [0,1]$, and $\psi(\partial, z) = 0, z\in [0,1]$.
 For $k=1, 2, \dots$, define
 \begin{equation}\label{def-A}
 A^{(k)}(x):=\frac{\partial^k}{\partial z^k}\psi(x,z)|_{z=1}.
\end{equation}
By {\bf(H2)}, $A^{(1)}$ and $A^{(2)}$ are finite, and
 \begin{equation}
  A^{(1)}(x) = \beta(x) \left(\sum_{k=1}^\infty kp_k(x) -1 \right),\quad  A^{(2)}(x) = \beta(x) \sum_{k=2}^\infty k(k-1)p_k(x).
 \end{equation}
 For any complex-valued Borel function $f$ on $E$
 and $(t,x) \in (0,\infty) \times E$, define
 \begin{equation}
  T_t f(x) : =  \mathbb{E}_x \left[e^{\int_0^t A^{(1)}(\xi_s)\mathrm{d}s} f(\xi_t) \right].
  \end{equation}
Then it is well-known that for any $t\ge 0$ and $x\in E$, $T_t f(x) = \mathbb{E}_{\delta_x} \left(\langle f, X_t\rangle\right)$,
see \cite[Lemma 1]{GHK2022} for example.
	
Under the assumptions {\bf(H1)} and {\bf(H2)}, there exists a family of continuous strictly positive functions $\{q_t(x,y): t \geq 0 \}$ on $E\times E$ such that
$$
T_t f(x) = \int_E q_t(x,y) f(y)\mu(\mathrm{d}y).
$$
Let
\[
\widehat{T}_tf(x):= \int_E q_t(y,x)  f(y)\mu(\mathrm{d}y)
\]
be the dual of $T_t$.  As summarized in \cite[Section 1]{RSZ2017},  both $(T_t)_{t\ge 0}$ and $(\widehat{T}_t)_{t\ge 0}$ are strongly continuous semigroups on $L^2(E,\mu;\mathbb{C})$ and, for any $t>0$, $T_t$ and $\widehat{T}_t$ are compact operators
with $L^2$ norm $\Vert T_t\Vert_2 = \Vert \widehat{T}_t\Vert_2 \leq e^{\Vert A^{(1)}\Vert_\infty t}$.
Let $\mathcal{L}$ and  $\widehat{\mathcal{L}}$ denote the infinitesimal generator of $(T_t)_{t\ge 0}$ and $(\widehat{T}_t)_{t\ge 0}$   in $L^2(E,\mu;\mathbb{C})$, then
the spectra $\sigma(\mathcal{L}) $ and $ \sigma(\widehat{\mathcal{L}})$ of $\mathcal{L}$ and $ \widehat{\mathcal{L}}$ both consist of
 eigenvalues of finite multiplicity only. $\mathcal{L}$ and $\widehat{\mathcal{L}}$ have the same number, say $N$, of eigenvalues. Of course $N$ might be finite or infinite.  We write $\mathbb{I}:=\{1,2,...,N\}$ when $N<
 \infty$ and $\mathbb{I}:=\{1,2,...,\}$ otherwise.
 The common value $-\lambda_1=\sup\mathfrak{R}(\sigma(\mathcal{L}))=\sup\mathfrak{R}(\sigma(\widehat{\mathcal{L}}))$ is
an eigenvalue of multiplicity one for both $\mathcal{L}$ and  $\widehat{\mathcal{L}}$. An eigenfunction $\phi_1$ of $\mathcal{L}$ associated with $-\lambda_1$ can be chosen to be strictly positive and continuous, and an eigenfunction $\widehat{\phi}_1$ of $\widehat{\mathcal{L}}$ associated with $-\lambda_1$ can also be chosen to be strictly positive and continuous.  Without loss of generality, we assume $\Vert \phi_1 \Vert_2 = 1$ and $\langle \phi_1, \widehat{\phi}_1 \rangle_\mu=1$.
We list the eigenvalues $\{-\lambda_k,k \in \mathbb{I}\}$ of $\mathcal{L}$ in an order so that $ -\lambda_1 > -\mathfrak{R}(\lambda_2) \geq
-\mathfrak{R}(\lambda_3)
 \geq  \dots$, then $\{-\overline{\lambda}_k, k\in \mathbb{I}\}$ are the eigenvalues of $\widehat{\mathcal{L}}$.
For simplicity we set $\mathfrak{R}_k:= \mathfrak{R}(\lambda_k)$ and $\mathfrak{I}_k:= \mathfrak{I}(\lambda_k)$ for $k\geq 2$.
In any finite vertical strip of the complex plane, there are at most finitely many $\lambda_k$'s.
Thus, in the case when $\mathbb{I}$ is infinite, $\mathfrak{R}_k\to \infty$ as $k\to\infty$.
 Define
 \begin{align}\label{Def-of-b-t}
 b_t(x):=\int_E q_t(x,y)^2\mu(\mathrm{d}y),\quad \widehat{b}_t(x):= \int_E q_t(y,x)^2\mu(\mathrm{d}y).
 \end{align}
 Using {\bf(H1)} and {\bf(H2)}, one can check, see \cite[Section 1.1]{RSZ2017}, that,
 for any $t>0$, $b_t$ and $\widehat{b}_t$ are continuous in $E$,  belong to $L^1(E,\mu)$ and that
 $b_t,\widehat{b}_t\in L^2(E,\mu)$ for all $t\geq t_0$.

Now we recall some spectral theoretic results
for $(T_t)_{t\geq 0}$ and $(\widehat{T}_t)_{t\geq 0}$ from \cite[Section 1]{RSZ2017}.
For each $k\in \mathbb{I}$, by \cite[(1.19),(1.27) and Lemma 1.11]{RSZ2017}, there exist  integers $n_k, r_k$,  $\{d_{k,j}, 1\leq j\leq r_k\}$,
families of continuous functions $\left\{\phi_j^{(k)}, 1\leq j\leq n_k\right\}\subset {\cal D}(\mathcal{L})\subset L^2(E,\mu; \mathbb{C})$ and $\left\{\widehat{\phi}_j^{(k)}, 1\leq j\leq n_k\right\} \subset {\cal D}(\widehat{\mathcal{L}})\subset L^2(E,\mu; \mathbb{C})$ such that
\[
\langle \phi_j^{(k)}, \phi_\ell^{(k)}\rangle_\mu =\delta_{j,l}:=
 1_{\left\{j=\ell\right\}}
= \langle \widehat{\phi}_j^{(k)}, \widehat{\phi}_\ell^{(k)}\rangle_\mu = \langle \phi_j^{(k)}, \widehat{\phi}_\ell^{(k)}\rangle_\mu,
\]
 $\sum_{j=1}^{r_k} d_{k,j}= n_k$ and that the $\mathbb{C}^{n_k}$-valued functions
\begin{align}
	\Phi_k(x):= (\phi_j^{(k)}(x), 1\leq j\leq n_k)^T\quad \mbox{and}\quad \widehat{ \Phi}_k(x):= (\widehat{\phi}_j^{(k)}(x), 1\leq j\leq n_k)^T
\end{align}
satisfy for all $x\in E$,
\begin{align}\label{T-t-Phik}
	T_t(\Phi_k)^T(x):= (T_t\phi_j^{(k)}(x), 1\leq j\leq n_k)= e^{-\lambda_k t} \left(\Phi_k(x)\right)^T D_k(t)
\end{align}
and
\begin{align}
	\widehat{T}_t(\widehat{\Phi}_k)(x):= (\widehat{T}_t\widehat{\phi}_j^{(k)}(x), 1\leq j\leq n_k)^T=  e^{-\overline{\lambda}_k t} D_k(t)\widehat{\Phi}_k(x),
\end{align}
where $D_k(t):=\mbox{diag}\left\{J_{k,j}(t), 1\leq j\leq r_k\right\} $ is an invertible matrix
with $D_k(t)D_k(s)=D_k(t+s)$ for all $s,t\in \mathbb{R}$ and,
 for $1\le j\le r_k$,   $J_{k,j}(t)$ is a $d_{k,j}\times d_{k,j}$
matrix given by $(J_{k,j}(t))_{p,q}:= 1_{\{q\geq p\}} t^{q-p}/(q-p)!$.
Moreover, $\phi_j^{(k)}, \widehat{\phi}_n^{(\ell)}\in L^2(E,\mu;\mathbb{C})\cap L^4(E,\mu;\mathbb{C})$ are continuous functions with $\langle \phi_j^{(k)}, \widehat{\phi}_n^{(\ell)} \rangle_\mu  =\delta_{k,\ell}\delta_{j,n}$.
By \cite[Remark 1.10]{RSZ2017}, for each $k\in \mathbb{I}$, there exists a unique $k'\in \mathbb{I}$ such that $\lambda_{k'}=\overline{\lambda}_k$. Since $D_k(t)=D_{k'}(t)$, we can choose $\Phi_{k'}(x) = \overline{\Phi_k(x)},$ which implies that $\widehat{\Phi}_{k'}(x) =\overline{ \widehat{\Phi}_{k}(x) }.$
The functions $\left\{\phi_j^{(k)}, 1\leq j\leq n_k\right\}$ are sometimes referred as the generalized eigenfunctions associated with $-\lambda_k$.

    We assume that the branching Markov process is supercritical, that is
    \begin{itemize}
    	\item [{\bf(H3)}] $\lambda_1 <0.$
    \end{itemize}

For a list of symmetric and non-symmetric spatial processes satisfying  {\bf(H1)},  see \cite{KS2008a, KS2008, KS2009} and \cite[Section 1.4]{RSZ20152}.
For $k\in \mathbb{I}$, we define
\begin{align}
	H_t^{(k)}:= e^{\lambda_k t}\left( \langle \phi_j^{(k)}, X_t\rangle, 1\leq j\leq n_k\right)(D_k(t))^{-1}.
\end{align}
According to  \cite[Lemma 3.1]{RSZ2017}, when $\lambda_1 > 2\mathfrak{R}_k$, for any $\nu\in \mathcal{M}_a(E)$ and $v\in\mathbb{C}^{n_k}$, $H_t^{(k)}v$ is an $L^2(\mathbb{P}_\nu)$-bounded martingale, which implies that the limit $H_\infty^{(k)}:=\lim_{t\to\infty} H_t^{(k)}$ exists $\P_\nu$-a.s. and in $L^2(\P_\nu)$.
For simplicity, we set $W_t:= H_t^{(1)}$ and $W_\infty:= H_\infty^{(1)}.$
Define $\mathcal{E}:= \{W_\infty =0\}.$

Spatial central limit theorems for linear functionals of $X$ were established in \cite{RSZ2017} when the spatial motion is not necessarily symmetric, generalizing the results of \cite{AM2015, RSZ} in the symmetric case.
To state the main results of \cite{RSZ2017}, we first introduce some notations.

For any $f \in L^2(E, \mu; \mathbb{C})$ and $k\in \mathbb{I}$, we define
\[
\langle f ,\widehat{\Phi}_k\rangle_\mu:= \left(\langle f, \widehat{\phi}_j^{(k)}\rangle_\mu, 1\leq j\leq n_k\right)^T\quad \mbox{and}\quad \gamma(f):= \inf \{k \in \mathbb{I}: \langle f ,\widehat{\Phi}_k\rangle_\mu \neq 0 \},
\]
here we use the usual convention $\inf \emptyset = \infty$. If $\gamma(f)<\infty$, define
\[
\zeta(f):= \sup\{k\in \mathbb{I}: \mathfrak{R}_k= \mathfrak{R}_{\gamma(f)}\}.
\]
Since for each $k\in \mathbb{I}$, every component of the function $t\mapsto D_k(t) \langle f, \widehat{\Phi}_k\rangle_\mu$ is a polynomial of $t$, we denote the degree of the $\ell$-th component of $ D_k(t) \langle f, \widehat{\Phi}_k\rangle_\mu$ by $\tau_{k,\ell} (f)$ and define
\begin{align}\label{Def-of-tau}
	\tau(f):= \sup\{\tau_{k,\ell}(f): \gamma(f)\leq k\leq \zeta(f), 1\leq \ell\leq n_k \}.
\end{align}
Then for any $k$ with $\mathfrak{R}_k=\mathfrak{R}_{\gamma(f)}$, the limit
\begin{align}\label{def-F(f,k)}
	F_{f,k}:= \lim_{t\to\infty} t^{-\tau(f)} D_k(t)\langle f, \widehat{\Phi}_k\rangle_\mu
\end{align}
exists and there exists $k$ such that $F_{f,k}\neq 0$.
Define
\begin{align}
	\mathcal{C}_{la}&:=\left\{g(x)= \sum_{k\in \mathbb{I}:\lambda_1>2\mathfrak{R}_k} \left(\Phi_k(x)\right)^T v_k: v_k\in \mathbb{C}^{n_k} \mbox{ with } \overline{v}_k=v_{k'}\right\}, \\
		\mathcal{C}_{cr}&:=\left\{g(x)= \sum_{k\in \mathbb{I}:\lambda_1=2\mathfrak{R}_k} \left(\Phi_k(x)\right)^T v_k: v_k\in \mathbb{C}^{n_k} \mbox{ with } \overline{v}_k=v_{k'}\right\},
		\\
			\mathcal{C}_{sm}&:=\left\{g \in
 L^2(E,\mu)\cap L^4(E,\mu):
 \lambda_1<2\mathfrak{R}_{\gamma(g)}\right\}.
\end{align}
 Note that $\mathcal{C}_{la}$, $\mathcal{C}_{cr}$ and $\mathcal{C}_{sm}$ consist of  real-valued functions, and that $\mathcal{C}_{la}$ and $\mathcal{C}_{cr}$ are of finite dimension and
 $\mathcal{C}_{cr}$ may be empty. $\mathcal{C}_{la}$ only involves $\Phi_k$'s associated with
 ``large'' eigenvalues $-\lambda_k$ satisfying $\lambda_1>2\mathfrak{R}_k$.  $\mathcal{C}_{cr}$ only involves $\Phi_k$'s associated with
 ``critical'' eigenvalues $-\lambda_k$ satisfying  $\lambda_1=2\mathfrak{R}_k$, if any.
 Any $f\in L^2(E,\mu)\cap L^4(E,\mu)$ can be decomposed as $f= f_{sm}+f_{cr}+f_{la}$ with
    \begin{align}
  	&	f_{la}(x):= \sum_{2\mathfrak{R}_k < \lambda_1} \left(\Phi_k(x)\right)^T v_k\in  \mathcal{C}_{la}, \quad
  	f_{cr}(x)
  	:= \sum_{2\mathfrak{R}_k = \lambda_1} \left(\Phi_k(x)\right)^T v_k\in  \mathcal{C}_{cr},\label{e:la-cr}\\
  	&\mbox{and}\quad  f_{sm}(x):= f(x)-f_{la} (x)-f_{cr}(x) \in \mathcal{C}_{sm}. \label{e:sm}
  \end{align}
  $f_{la}$, $f_{cr}$ and $f_{sm}$ are called the large, critical and  small  components of $f$ respectively.
  We define $\sigma_{sm}^2(f),\sigma_{cr}^2(f)$ and $\sigma_{la}^2(f)$ by
\begin{align}
		\sigma_{sm}^2(f)&:= \int_0^\infty e^{\lambda_1s}\langle A^{(2)}\cdot \left|T_s f_{sm}\right|^2, \widehat{\phi}_1 \rangle_\mu \mathrm{d}s + \langle |f_{sm}|^2, \widehat{\phi}_1\rangle_\mu,	\label{e:varsm}\\
		\sigma_{cr}^2(f) & := (1+2\tau(f_{cr}))^{-1} \sum_{k:\lambda_1=2\mathfrak{R}_k}  \langle A^{(2)} \cdot \left|(\Phi_k)^T F_{f_{cr}, k}\right|^2, \widehat{\phi}_1\rangle_\mu,
		\label{e:varcr}\\
		\sigma_{la}^2(f) &:= \int_0^\infty e^{-\lambda_1 s}\Big\langle A^{(2)}\cdot |I_s f_{la}|^2, \widehat{\phi}_1\Big\rangle_\mu \mathrm{d}s -\langle \left|f_{la}\right|^2, \widehat{\phi}_1 \rangle_\mu,\label{e:varla}
	\end{align}
	where
	\[
	I_sf_{la}(x):= \sum_{k:\lambda_1>2\mathfrak{R}_k} e^{\lambda_k s}\left(\Phi_k(x)\right)^TD_k(s)^{-1} v_k.
	\]
For any $f\in  L^2(E,\mu)\cap L^4(E,\mu)$, it was shown in \cite[Theorem 1.16]{RSZ2017} that $\sigma_{sm}^2(f)\in (0,\infty)$ if $f_{sm}\neq 0$ and similar results hold for $f_{cr}$ and $f_{la}$.
	Define
	\begin{align}\label{def-of-e-t}
	E_t(f_{la}):=  \sum_{k:\lambda_1>2\mathfrak{R}_k} \left(e^{-\lambda_k t}H_\infty^{(k)}D_k(t)v_k\right).
	\end{align}

{\it In this paper, for any $f\in L^2(E,\mu)\cap L^4(E,\mu)$, we will always use the	notations	$\mathcal{C}_{la}$, $\mathcal{C}_{cr}$, $\mathcal{C}_{sm}$, $\sigma_{sm}^2(f)$, $\sigma_{cr}^2(f)$ and	 $\sigma_{la}^2(f)$ defined above. }

	Recall that $\mathcal{E}=\left\{W_\infty =0\right\}$.
    The spatial central limit theorem of \cite{RSZ2017} is follows.
\begin{thrm} (\cite[Theorem 1.16]{RSZ2017})\label{CLT}
	If $f\in L^2(E,\mu)\cap L^4(E,\mu)$,
	then $\sigma_{sm}^2(f),\sigma_{cr}^2(f), \sigma_{la}^2(f) \in
	[0,\infty)$.
	 Furthermore, under $\mathbb{P}_{\delta_x}(\cdot| \mathcal{E}^c)$, as $t\to\infty$,
	\begin{align}
		& \left(e^{\lambda_1 t}\langle \phi_1, X_t\rangle, \frac{\langle f_{la},X_t\rangle-E_t(f_{la})}{\sqrt{\langle \phi_1, X_t\rangle}}, \frac{\langle f_{cr}, X_t\rangle}{\sqrt{t^{1+2\tau(f_{cr})}\langle \phi_1, X_t\rangle}}, \frac{\langle f_{sm}, X_t\rangle}{\sqrt{\langle \phi_1, X_t\rangle}} \right)\\
		&\quad\quad \stackrel{\mathrm{d}}{\Rightarrow} (W^*, G_{la}, G_{cr}, G_{sm}),
	\end{align}
	where  $W^*$ has the same law as $W_\infty$ conditioned on $\mathcal{E}^c, G_{la} \sim \mathcal{N}\left(0, \sigma_{la}^2(f) \right), G_{cr} \sim \mathcal{N}\left(0, \sigma_{cr}^2(f)\right), G_{sm} \sim \mathcal{N}\left(0, \sigma_{sm}^2(f) \right)$ and that $W_\infty^*, G_{la}, G_{cr}$ and $G_{sm}$ are independent.
\end{thrm}

The main purpose of this paper is to complement the CLT type results above for $ \langle f, X_t \rangle$ with law of iterated logarithm type results for $\langle f, X_t \rangle$.

\section{Main results}
Our first two results are LIL type results in the special case when
$f$ is of the form
\begin{align}\label{e:rsspecial}
f(x)= \sum_{k=1}^m \left(\Phi_k(x)\right)^T v_k, \quad \mbox{ for some } m\in \N \mbox{ and } v_k\in \mathbb{C}^{n_k} \mbox{ with } \overline{v}_k=v_{k'}.
\end{align}
In the symmetric case, functions of the form \eqref{e:rsspecial} are dense in  $L^2(E, \mu)$.

  \begin{thrm}\label{thm2}
	Suppose {\bf(H1)}--{\bf(H3)} hold and $f$ is of the form \eqref{e:rsspecial}.
	If  $f_{cr}=0$, then  $\mathbb{P}_{\delta_x}\left(\cdot \vert \mathcal{E}^c\right)$-almost surely,
  	 \begin{equation}\label{e:t2.3lims}
  		\limsup_{t\to\infty}/\liminf_{t\to\infty} \frac{e^{\lambda_1  t/2}\left(\langle f,X_t\rangle-E_t(f_{la})\right)}{\sqrt{2\log t}} = +/- \sqrt{\left(\sigma_{sm}^2(f)+ \sigma_{la}^2(f)\right) W_\infty}.
  	\end{equation}
  \end{thrm}

  \begin{remark}\label{remark-LIL}
  	 Note that 	Theorem \ref{thm2}  is equivalent to that, for $f$  of the form \eqref{e:rsspecial} with $f_{cr}=0$, $\mathbb{P}_{\delta_x}\left(\cdot \vert \mathcal{E}^c\right)$-almost surely,
  	 \begin{align*}
			\limsup_{t\to\infty}/\liminf_{t\to\infty} \frac{e^{\lambda_1  t/2}\left(\langle f,X_t\rangle-E_t(f_{la})\right)}{\sqrt{2\log \log \langle \phi_1,  X_t\rangle}}
			&=  +/-	 \sqrt{\left(\sigma_{sm}^2(f)+ \sigma_{la}^2(f)\right) W_\infty}.
  	 \end{align*}
	Thus, the result above is a law of  iterated logarithm in some sense. In this paper, we will call results like Theorem \ref{thm2} ``law of iterated logarithm" following the convention of \cite{IK, IKM}.
  \end{remark}

Our next  theorem gives the law of iterated logarithm
for $\langle f , X_{t}\rangle$ for the case when  $f_{cr}\neq 0$.

\begin{thrm}\label{thm4}
	Suppose {\bf(H1)}--{\bf(H3)} hold and $f$ is of the form \eqref{e:rsspecial}.
	If $f_{cr}\neq 0$, then  $\mathbb{P}_{\delta_x}\left(\cdot \vert \mathcal{E}^c\right)$-almost surely,
	\begin{equation*}
		\limsup_{t\to\infty}/ \liminf_{t\to\infty} \frac{e^{\lambda_1  t/2}\left(\langle f,X_t\rangle - E_t(f_{la})\right)}{\sqrt{2t^{1+2\tau(f_{cr})}\log\log t}} = +/ - \sqrt{\sigma_{cr}^2(f)W_\infty},
	\end{equation*}
	where $\tau(f)$ is given as in \eqref{Def-of-tau}.
\end{thrm}

  \begin{remark}\label{Remark1}
 In the special case
 where $X$  is  a (finite) multitype branching process, our results are consistent with \cite[Theorem 2]{Asmussen1977}.
For test functions (vectors) with non-trivial ``large component'',  \cite[Theorem 4]{Asmussen1977}
is only for eigenvectors corresponding to large eigenvalues.
We need some new idea to handle general test functions $f$, especially when the critical component $f_{cr}$ is non-trivial.
\end{remark}

Theorems \ref{thm2} and \ref{thm4} are for functions of the form \eqref{e:rsspecial} only and the proofs crucially use this assumption.
To extend Theorems \ref{thm2}  and \ref{thm4} to
more general functions, we need the following stronger assumption and a different argument.

  \begin{itemize}
  	\item [{\bf(H4)}]
	{\bf(a)} $\widehat{\phi}_1$ is bounded; {\bf(b)} $\sup_{x\in E} \sum_{k=0}^\infty k^4 p_k(x)<\infty$.
  \end{itemize}

First, we give an example showing that LIL is not true for all test functions. Consider the 1-dimensional branching OU-process with branching rate $\beta=1, p_2=1$ and suppose that $f(x)= 1_{x\neq 0}+ \infty 1_{x=0}$.  Since the  1-dimensional OU-process is Harris recurrent, $\mathbb{P}_{\delta_x}$-almost surely the set $\mathcal{J}:= \{t<\infty: X_t(\{0\})\neq 0\}$ contains
a sequence of times $t_k$ increasing $\infty$.  Thus $\mathbb{P}_{\delta_x}$-almost surely, $\langle f, X_{t_k}\rangle = \infty$ for each $k$ and so there is no LIL-type result for this function $f$. Thus, for LIL, we do need some regularity assumption on the test function $f$. The following condition will play an important role in our argument below:
\begin{align}\label{Semigroup}
T_{s}f(x)- f(x)= \int_0^{s}T_r (\mathcal{L}f)(x)\mathrm{d}r, \qquad \text{ for all } s>0 \text{ and } \ x\in E.
\end{align}
Recall that $\mathcal{L}$ is the generator of $(T_t)_{t\ge 0}$ in $L^2(E, \mu, \mathbb{C})$ and the fact that the equality above is valid for all $s>0$ and almost every $x\in E$ is well known.

Now we introduce our  space of test functions.
Let $\mathcal{M}$ be the space of real valued functions in the closure  of the linear span of $\{\phi_j^{(k)}: k\in \mathbb{I}, 1\leq j\leq n_k\}$ in $L^2(E, \mu; \mathbb{C})$.
In the symmetric case,
$\mathcal{M}=L^2(E, \mu)$. Define
   \begin{align}
		\mathcal{T}	:=\bigg\{ &	f	\in
		\mathcal{M}\cap
		\mathcal{D}(\mathcal{L}):
		 \frac{f}{b_{4t_0}^{1/2}}\in L^\infty (E,\mu),
		 \  \mathcal{L}f\in  L^4(E,\mu),
		 \  f \mbox{ satisfies } \eqref{Semigroup}
	\bigg\},
   \end{align}
 where  $b_t$ is defined in \eqref{Def-of-b-t}.

Note that any function $f\in \mathcal{M}$ is the $L^2$ limit of a
sequence $\{f_k, k\in \mathbb{N}\}$ of functions
of form \eqref{e:rsspecial} and that $\gamma(f)<\infty$.
Using Lemma \ref{Useful-known-result} (2) below, it is easy to see
that any function of the form \eqref{e:rsspecial} is in $\mathcal{T}$.
Let $g\in \mathcal{M}$, then there exists a  sequence of functions $g_k$ of form \eqref{e:rsspecial}
converging to $g$ in $L^2(E, \mu)$. It is easy to see that, for any $r>0$ and $\lambda>-\lambda_1$,  $f_k:= T_r R_\lambda g_k$ is also of form \eqref{e:rsspecial} and that $f_k$ converges in $L^2$ to $f:=T_rR_\lambda g$.  Using Lemma \ref{lemma3}, one can easily check that, if $r>8t_0$, then
$f:=T_rR_\lambda g$ belongs to $\mathcal{T}$.  Thus, for any $r>8t_0$ and $\lambda>-\lambda_1$, $T_rR_\lambda(\mathcal{M})\subset \mathcal{T}$.
 In the case when $\mathbb{I}$ is finite, all the functions in $\mathcal{T}$ are of the form
 \eqref{e:rsspecial}.

We mention here that if $\mathbb{I}$ is finite, Theorems \ref{thm2} and \ref{thm4} give the full law of iterated logarithm theorem.
The set ${\cal T}$ is only used to treat the case when $\mathbb{I}$ is infinite.
Here is our  law of iterated logarithm theorem for general $f\in {\cal T}$.

\begin{thrm}\label{thm5}
If {\bf(H1)}--{\bf(H4)} hold, then the conclusions of
Theorems \ref{thm2}   and \ref{thm4} hold for any $f\in {\cal T}$.
  \end{thrm}

The proof of Theorem \ref{thm5} is different from
that of Theorems \ref{thm2}  and \ref{thm4}.
One of the key differences is that we choose a different discretization scheme.

We mention here that  {\bf(H4)}{\bf(a)} is used once only to show  $\sigma_{sm}^2(f)\lesssim \Vert f\Vert_2^2$ in the proof of Theorem \ref{thm5}, while {\bf(H4)}{\bf(b)} is used only in the proof of Lemma \ref{Useful-lemma} to bound $\mathbb{E}_{\delta_x}(|\langle f, X_t\rangle|^3)$ and $\mathbb{E}_{\delta_x}(|\langle f, X_t\rangle|^4)$ from above.

Now we compare our results with existing results. The most closely related paper is
Asmussen \cite{Asmussen1977}
 on multi-type branching processes. \cite[Theorem 1]{Asmussen1977} contains LILs for test functions (vectors) with trivial large components. For test functions (vectors) with nontrivial large components, \cite[Theorem 4]{Asmussen1977} only considered the eigenfunction functions (eigenvectors) associated with real-eigenvalues and proved an LIL for the martingales associated with these eigenfunctions. Our model is more general in that  our spatial motion is a general non-symmetric Markov process and our branching mechanism is spatially dependent. For test functions with non-trivial large components, we allow them to be linear
combinations of (generalized) eigenfunctions associated all (real or complex) eigenvalues.  The papers  \cite{IK, IKM} contain LIL-type results for non-negative martingales of
general branching processes.
To the best of our knowledge, the main results of this paper are the first almost sure fluctuation results
for signed linear functional of branching Markov processes.

We end this section with a brief description of the strategy and organization of this paper.  In Section \ref{ss:3}, we gather some useful results and give a general law of iterated logarithm for sequence of random variables. In Subsection \ref{ss: 4.1} we give some general results and we
 prove Theorems \ref{thm2}   and \ref{thm4}
 in Subsections \ref{ss: 4.2} and \ref{ss:4.3} respectively. The proof of Theorem \ref{thm5} is given in Section \ref{ss: 5}.

We believe that the general idea of this paper
can be adapted to
non-local branching Markov process \cite{DH2025} and superprocesses
\cite{GHK2022,PY2020, RSZ2015, YT2025}.
Our approach can also be adapted to prove an LIL for the non-negative martingale associated to the
principal eigenfunction (or ground state) for the branching symmetric Markov processes treated in \cite{CS}.
We will not pursue these in this paper.

    \section{Preliminary}\label{ss:3}

 Throughout this paper, we always assume that {\bf(H1)}--{\bf(H3)} hold.
  We use $F(x)\lesssim_{r,f,\kappa,...} G(x), x\in E,$ to denote that there exists some constant $C=C(r,f,\kappa,...)$ such that $F(x)\leq CG(x)$ for all $x\in E$.

  We first give some preliminary results on the moments of $\langle f, X_t \rangle$ for $f\in L^2(E, \mu)$.
    For any $f\in L^2(E,\mu;\mathbb{C})$, define
\begin{align}
		\widetilde{f}(x):= f(x)-\sum_{j=\gamma(f)}^{\zeta(f)} (\Phi_j(x))^T\langle f, \widehat{\Phi}_j\rangle_\mu.\label{step_4}
\end{align}
Note that, if $f \in L^2(E,\mu)$, then  $\widetilde{f}$ is real-valued.
  For any real-valued random variable $Y$, we define
    $$ \textup{Var}_x(Y|\mathcal{F}):= \mathbb{E}_{\delta_x} \left[Y^2 \big\vert \mathcal{F}\right] - \left(\mathbb{E}_{\delta_x}\left[Y \big\vert \mathcal{F}\right]\right)^2.$$
Here and throughout the paper we use the notation
 $\textup{Var}_x(Y) = \mathbb{E}_{\delta_x}\left(Y^2\right)- \left(\mathbb{E}_{\delta_x}(Y)\right)^2$.

\begin{lemma}\label{lemma3}
Assume $f\in L^2(E,\mu)$.
\newline
$(1)$ For any $a\in(\lambda_1, \mathfrak{R}_2)$ and any $t_1>0$, we have
\[
\left|e^{\lambda_1t} T_t f(x) -\langle f, \widehat{\phi}_1\rangle_\mu \phi_1(x)\right|
 \lesssim_{a, t_1}
e^{-(a-\lambda_1)t} \Vert f\Vert_2 b_{t_1}^{1/2}(x),\quad  t>2t_1, x\in E.
\]
$(2)$ If $\gamma(f)<\infty$, then for any $t_1 > 0,$ we have
\begin{equation}
t^{-\tau(f)} e^{\mathfrak{R}_{\gamma(f)}t} | T_t f(x) |\lesssim_{\gamma(f), t_1} \Vert f \Vert_2 b_{t_1}^{1/2} (x),\  \quad  t>2t_1, x\in E.
\end{equation}
 Consequently, for any $f\in L^2(E,\mu;\mathbb{C})$, $T_tf\in L^2(E,\mu;\mathbb{C})$ for any $t>0$ and $T_tf\in L^4(E,\mu;\mathbb{C})$ for any $t>2t_0$.
   \end{lemma}
   \textbf{Proof:}
   (1) follows from \cite[(2.16)]{RSZ2017}, so we prove (2) here.
   By \cite[Lemma 2.2]{RSZ2017}, for
   any fixed $a\in (\mathfrak{R}_{\gamma(f)}, \mathfrak{R}_{\zeta(f)+1})$,
   we have
   \begin{align}\label{Ttf-main}
   	  \left|T_t f(x)- \sum_{j=\gamma(f)}^{\zeta(f)}e^{-\lambda_j t}\left(\Phi_j(x)\right)^T D_j(t) \langle f, \widehat{\Phi}_j\rangle_\mu\right|
   	  & \lesssim_{\gamma(f), a,t_1} e^{-at}b_{t_1}^{1/2}(x)\int_E |f(y)| \widehat{b}_{t_1}^{1/2}(y)\mu(\mathrm{d} y) \nonumber\\
   	  &\lesssim_{ \gamma(f), a, t_1} \Vert f\Vert_2 e^{-at}b_{t_1}^{1/2}(x).
   \end{align}
Using  \cite[(1.20)]{RSZ} with $t=t_1$, we get $\|\Phi_j(x)\|_\infty\lesssim_{j} b_{t_1}^{1/2}(x)$, and then, using $|\langle f, \widehat{\phi}_j^{(k)}\rangle_\mu| \lesssim_{j,k} \Vert f \Vert_2$, we get
$
   |\left(\Phi_j(x)\right)^T D_j(t) \langle f, \widehat{\Phi}_j\rangle_\mu|\lesssim_{j, t_1}  t^{\tau(f)}b_{t_1}^{1/2}(x) \Vert f\Vert_2.
 $
Therefore, the assertions of (2) hold by \eqref{Ttf-main}.

   \hfill$\Box$

As a consequence of Lemmas  \ref{lemma3} (2),
we have the following  inequality:	for any $R>3t_0$ and $s\in (3t_0, R]$,
\begin{align}\label{T-s-b-t-0}
	T_{s}(b_{t_0})(x)
	\lesssim_{R, t_0} b_{t_0}^{1/2}(x)\land T_{s-3t_0}(b_{t_0}^{1/2})(x).
\end{align}
We collect some useful estimates  obtained in \cite{RSZ2017}.

\begin{lemma}\label{Useful-known-result}
$(1)$
For any $R>0$ and $f\in L^2(E,\mu, \mathbb{C})\cap L^4(E,\mu, \mathbb{C})$, we have
	$\mathbb{E}_{\delta_x}(|\langle f, X_r\rangle^2|)\lesssim_R T_r(|f|^2)(x)$ for all $r\in (0, R]$.
	\\
$(2)$
 For each $k\in \mathbb{I}$, $\sup_{1\leq j\leq n_k} |\phi_j^{(k)}|\lesssim_{k, t_0} b_{t_0}^{1/2}$.
\\
$(3)$
 For any $t>0, x\in E$, $b_{4t}(x) \lesssim_{t}
 T_{2t}(a_{2t})(x)$
and $b_{4t}(x) \lesssim_{t} T_{3t}(a_{t})(x)$.
\end{lemma}
\textbf{Proof: }
For (1) and (2), see  \cite[(2.19)]{RSZ2017} and \cite[(1.20)]{RSZ2017} respectively.
For the first inequality of  (3), see the display below \cite[(2.23)]{RSZ2017}, and the second equality of (3) follows similarly.

\hfill$\Box$

    \begin{lemma}\label{lemma2}
	Assume $f\in  L^2(E,\mu) \cap  L^4(E,\mu)$.
    \newline
    $(1)$ For any $x\in E,$
    $$
    	\lim_{t\to\infty} e^{\lambda_1t/2}\left|\mathbb{E}_{\delta_x} \left(\langle f_{sm}, X_t\rangle \right)\right|= 0,\quad
    	\lim_{t \to \infty} e^{\lambda_1 t}\mathbb{E}_{\delta_x}\left(\langle f_{sm}, X_t\rangle^2\right) = \sigma_{sm}^2(f) \phi_1(x).
    $$
   Moreover, for any $t>10 t_0$ and  $x\in E$,
   \begin{equation}
  	  	e^{\lambda_1 t} \mathbb{E}_{\delta_x}
  	  \left(\langle f_{sm}, X_t \rangle^2\right) \lesssim_{f, t_0} b_{t_0}^{1/2}(x)+ b_{t_0}(x),
   \end{equation}
   and when $\mathfrak{R}_{\gamma(f_{sm})}>0$, it holds that
   \begin{equation}\label{step_10}
   	 \sigma_{sm}^2(f)\lesssim \Vert f_{sm}\Vert_2^2 + \langle |f_{sm}|^2, \widehat{\phi}_1\rangle_\mu .
   \end{equation}
    $(2)$ For any $t>10t_0$ and $x\in E$, it holds that
    \begin{align*}
    	\left|t^{-(1+2\tau(f_{cr}))}e^{\lambda_1 t}\textup{Var}_x\left( \langle f_{cr}, X_t \rangle\right) - \sigma_{cr}^2 (f)\phi_1(x) \right|\lesssim_{t_0, f} t^{-1}\left(b_{t_0}^{1/2}(x)+ b_{t_0}(x)\right).
    \end{align*}
    $(3)$  For any $t>10 t_0$  and $x\in E$,
  \begin{equation*}
  	t^{-2\tau(f)}e^{2\mathfrak{R}_{\gamma(f)} t}\mathbb{E}_{\delta_x}\left(\langle f_{la}, X_t\rangle^2 \right)\lesssim_{t_0, f} b_{t_0}^{1/2}(x).
  \end{equation*}
    \end{lemma}
\textbf{Proof:}
All the assertions, except \eqref{step_10}, follow from \cite[Lemmas 2.5, 2.6 and 2.7]{RSZ2017}.
Now we prove \eqref{step_10}.
Combining the inequality $|A^{(2)}| |T_s f_{sm}|^2 \leq e^{\Vert A^{(1)}\Vert_\infty s}\Vert A^{(2)}\Vert_\infty T_s(|f_{sm}|^2) \lesssim_{t_0} T_s(|f_{sm}|^2)$ for all $s\leq 2t_0$ and that $\langle T_s (|f_{sm}|^2), \widehat{\phi}_1\rangle_\mu =e^{-\lambda_1 s}\langle |f_{sm}|^2, \widehat{\phi}_1\rangle_\mu$,  we conclude that
\begin{align}\label{step_22}
	\sigma_{sm}^2(f)& = \int_0^\infty e^{\lambda_1s}\langle A^{(2)}\left|T_s f_{sm}\right|^2, \widehat{\phi}_1 \rangle_\mu \mathrm{d}s + \langle |f_{sm}|^2, \widehat{\phi}_1\rangle_\mu \nonumber \\
	& \lesssim_{t_0}  \int_{2t_0}^\infty e^{\lambda_1s}\langle \left|T_s f_{sm}\right|^2, \widehat{\phi}_1 \rangle_\mu \mathrm{d}s +  \langle |f_{sm}|^2, \widehat{\phi}_1\rangle_\mu.
\end{align}
Let $k_0:= \sup\{k:\mathfrak{R}_k\leq 0\}$.  Taking $a=0, k = k_0$ and $t_1=t_0$ in
\cite[Lemma 2.2]{RSZ2017}, we have for all $t>2 t_0$ and $x\in E$,
\begin{align}
	|T_t f_{sm}(x)| & = \left|\int_E \left(q_t(x,y) -\sum_{j=1}^{k_0} e^{-\lambda_j t}\left(\Phi_j(x)\right)^T D_j(t)\overline{\widehat{\Phi}_j(y)} \right) f_{sm}(y)\mu(\mathrm{d}y)\right|\\
	& \lesssim b_{t_0}^{1/2}(x)\int_E \widehat{b}_{t_0}^{1/2}(y) \left| f_{sm}(y)\right|\mu(\mathrm{d}y)\lesssim \Vert f_{sm}\Vert_2b_{t_0}^{1/2}(x).
\end{align}
Plugging this back to \eqref{step_22} yields that
\begin{align}
		\sigma_{sm}^2(f)
	& \lesssim_{t_0}  \Vert f_{sm}\Vert_2^2 \int_{2t_0}^\infty e^{\lambda_1s}\langle b_{t_0}^{1/2}, \widehat{\phi}_1 \rangle_\mu \mathrm{d}s +  \langle |f_{sm}|^2, \widehat{\phi}_1\rangle_\mu  \lesssim \Vert f_{sm}\Vert_2^2 +\langle |f_{sm}|^2, \widehat{\phi}_1\rangle_\mu,
\end{align}
which implies \eqref{step_10}.

\hfill$\Box$

 \begin{lemma}\label{lemma4}
Suppose that $f \in L^2(E,\mu) \cap  L^4(E,\mu)$
with $\lambda_1 > 2\mathfrak{R}_{\gamma(f)}$
and recall $\widetilde{f}$ is defined in \eqref{step_4}.
Then there
 exists $ c(f) > 0$ such that  for any $t >10 t_0$ and $x\in E$,
	$$ e^{2\mathfrak{R}_{\gamma(f)} t}\mathbb{E}_{\delta_x} \left[\langle \widetilde{f}, X_t \rangle^2 \right] \lesssim_{f, t_0} e^{-c(f) t} \left(b_{t_0}^{1/2}(x)+ b_{t_0}(x)\right).$$
\end{lemma}
\textbf{Proof:}
See the proof of \cite[Theorem 1.14, (3.11)]{RSZ2017}.
Moreover, one can choose $c(f) <2(\mathfrak{R}_{\gamma(\widetilde{f})}-\mathfrak{R}_{\gamma(f)})$ if $\lambda_1>2\mathfrak{R}_{\gamma(\widetilde{f})}$ and $c(f) < \mathfrak{R}_{\gamma(\widetilde{f})}- \mathfrak{R}_{\gamma(f)}$ if $\lambda_1=2\mathfrak{R}_{\gamma(\widetilde{f})}$  and $c(f) <\lambda_{1}- 2\mathfrak{R}_{\gamma(f)}$ if $2\mathfrak{R}_{\gamma(\widetilde{f})}>\lambda_1> 2\mathfrak{R}_{\gamma(f)}$.

\hfill$\Box$

As an application of Lemma \ref{lemma4}, we have the following strong law of large numbers type result.
\begin{lemma}\label{lemma10}
For any
$f \in L^2(E,\mu) \cap  L^4(E,\mu)$ and $\delta >0$, we have
$$
\lim_{n\to\infty} e^{\lambda_1 n\delta}\langle f, X_{n\delta}\rangle = \langle f, \widehat{\phi}_1\rangle_\mu W_\infty,\quad \P_{\delta_x}\mbox{-a.s.}
$$
\end{lemma}
\textbf{Proof:}
We only treat the case $f\geq 0$ since for general $f$,
we can treat the positive and negative parts of $f$ separately.
Note that $\widetilde{f}(x)= f(x) - \langle f,\widehat{\phi}_1\rangle_\mu \phi_1(x)$, by Lemma \ref{lemma4},   for any $n\in \N$ with $n>10 t_0/\delta$,
\begin{equation}
	e^{2\lambda_1 n\delta} \mathbb{E}_{\delta_x} \left[\langle \widetilde{f}, X_{n\delta}\rangle^2 \right] \lesssim_{f,t_0} e^{-c(f) n\delta}\left(b_{t_0}^{1/2}(x)+ b_{t_0}(x)\right) .
\end{equation}
 Thus, for any $\varepsilon >0,$
 by Markov's inequality,
$$
\sum_{n\geq0} \mathbb{P}_{\delta_x} \left(\left\vert  e^{\lambda_1 n\delta} \langle \widetilde{f}, X_{n\delta}\rangle \right\vert > \varepsilon \right) \lesssim_{f, t_0} 1+\frac{10 t_0}{\delta}+  \frac{1}{\varepsilon^2} \sum_{n\geq 0} e^{-c(f) n\delta} \left(b_{t_0}^{1/2}(x)+ b_{t_0}(x)\right)< \infty,
$$
which implies that  $ e^{\lambda_1 n\delta} \langle \widetilde{f}, X_{n\delta}\rangle$ converges to $0\  \mathbb{P}_{\delta_x}$-a.s.
Since $e^{\lambda_1 n\delta}\langle \widetilde{f}, X_{n\delta} \rangle = e^{\lambda_1 n\delta} \langle f, X_{n\delta}\rangle -  \langle f, \widehat{\phi}_1 \rangle_\mu W_{n\delta} $ and $ \langle f, \widehat{\phi}_1 \rangle_\mu  W_{n\delta}$ converges to $\langle f, \widehat{\phi}_1\rangle_\mu  W_\infty$ almost surely,
the assertion of the lemma follows immediately.

\hfill$\Box$

Now we give some useful limit results for sequence of real-valued random variables.

    \begin{lemma}\label{lemma6}
    (\cite[Lemma A.2.]{IK})
	Let $X_1,X_2,...$ be independent
	real-valued random variables
	with $\mathbb{E}X_i = 0$ and $ \mathbb{E}|X_i|^3 < \infty , i = 1,2,...$ . If $\sum_{i\geq 1} \mathbb{E} X_i^2 < \infty$, then there exists an absolute constant
	$C$ such that
	\begin{equation*}
		\sup_{y\in\mathbb{R}}  \left\vert \mathbb{P} \left[\frac{\sum_{i\geq 1} X_i }{\sqrt{\sum_{i\geq 1}\mathbb{E} X_i^2 }} \leq y \right] - \Phi(y)\right\vert \leq
		C
		 \frac{\sum_{i\geq 1} \mathbb{E}|X_i|^3}{ \sqrt{\left( \sum_{i\geq 1}\mathbb{E} X_i^2\right)^3 }},
	\end{equation*}
	where $\Phi(y) = (1/\sqrt{2\pi})\int_{-\infty}^y e^{-x^2/2} \mathrm{d}x, y \in \mathbb{R}.$
\end{lemma}

\begin{lemma}\label{lemma11}
	For any $\delta>0$ and any
	real-valued random variable $Y$ with
	$\mathbb{E}[Y^2]< \infty$, it holds that
	\begin{equation}\label{step_24}
		\sum_{n\geq 0} e^{\lambda_1 n\delta /2}\mathbb{E}\left[ |Y|^3 1_{\{|Y| \leq e^{-\lambda_1 n\delta/2}\}}\right] + \sum_{n\geq 0} e^{-\lambda_1 n\delta /2}\mathbb{E}\left[ |Y| 1_{\{|Y|>e^{-\lambda_1 n\delta/2}\}}\right] \lesssim_{\delta} \mathbb{E}[Y^2].
	\end{equation}
\end{lemma}
\textbf{Proof:} Define $n_y: = \inf\{
n\in \mathbb{N}:
 e^{-\lambda_1 n\delta /2}\geq y\}$. Combining the inequalities $\mathbb{E}[|Y|^3 1_{\{|Y| \leq K \}}] \leq 3\int_0^K y^2 \mathbb{P}(|Y| > y)\mathrm{d}y$ and
    $\sum_{n\geq 0} e^{-\lambda_1 n\delta /2}\mathbb{E}\left[ |Y| 1_{\{|Y|
    >  e^{-\lambda_1 n\delta/2}\}}\right] \leq  \mathbb{E}\left[ |Y| \sum_{n=0}^{n_Y} e^{-\lambda_1 n\delta/2}\right]$, we have
\begin{align*}
	& \sum_{n\geq 0} e^{\lambda_1 n\delta /2}\mathbb{E}\left[ |Y|^3 1_{\{|Y| \leq e^{-\lambda_1 n\delta/2}\}}\right] + \sum_{n\geq 0} e^{-\lambda_1 n\delta /2}\mathbb{E}\left[ |Y| 1_{\{|Y|>	e^{-\lambda_1 n\delta/2}\}}\right] \\
	& \leq  3 \sum_{n\geq 0} e^{\lambda_1 n\delta /2} \int_0^{ e^{-\lambda_1 n\delta/2}} y^2 \P(|Y| >y) \mathrm{d}y +
	\mathbb{E}\left[ |Y| \sum_{n=0}^{n_Y} e^{-\lambda_1 n\delta/2}\right]
	\\ & =
	3\int_0^\infty  y^2 \P(|Y| >y) \sum_{n\geq n_y} e^{\lambda_1 n\delta /2}\mathrm{d}y +
	\frac{1}{e^{-\lambda_1 \delta/2}-1}\mathbb{E}\left[ |Y| \left(e^{-\lambda_1 (n_Y+1)\delta/2 }-1\right)\right]
	\\ & \leq 3 \sum_{n\geq 0}  e^{\lambda_1 n\delta /2} \int_0^\infty  y^2 \P(|Y| >y) e^{\lambda_1 n_y\delta  /2}\mathrm{d}y +
	\frac{1}{e^{-\lambda_1 \delta/2}-1}\mathbb{E}\left[ |Y| e^{-\lambda_1 (n_Y+1)\delta/2 }\right].
\end{align*}
 Since $ e^{-\lambda_1 \delta /2}y \geq e^{-\lambda_1 n_y\delta  /2} \geq y$ and $\mathbb{E}[Y^2] = \int_0^\infty 2y \mathbb{P}(|Y|>y)\mathrm{d} y,$  we conclude that
 \begin{align}
 		& \sum_{n\geq 0} e^{\lambda_1 n\delta /2}\mathbb{E}\left[ |Y|^3 1_{\{|Y| \leq e^{-\lambda_1 n\delta/2}\}}\right] + \sum_{n\geq 0} e^{-\lambda_1 n\delta /2}\mathbb{E}\left[ |Y| 1_{\{|Y| \geq e^{-\lambda_1 n\delta/2}\}}\right] \\
 	& \leq 3 \sum_{n\geq 0}  e^{\lambda_1 n\delta /2} \int_0^\infty  y \P(|Y| >y) \mathrm{d}y +
 	 \frac{e^{-\lambda_1 \delta}}{e^{-\lambda_1 \delta/2}-1}\mathbb{E}\left[ Y^2 \right]
 	\nonumber\\
 	& = \frac{1}{2}\left(3 \sum_{n\geq 0}  e^{\lambda_1 n\delta /2}+\frac{2e^{-\lambda_1 \delta}}{e^{-\lambda_1 \delta/2}-1}
 	\right) \mathbb{E}(Y^2),
 \end{align}
 which implies the desired result.
\hfill$\Box$

\begin{lemma}\label{lemma1}
	Let $\{\mathcal{G}_n: n =0,1,...\}$ be an increasing sequence of $\sigma$-fields and $B$ be an event with positive probability. Let $\{T_n: n =0,1,...\}$ be a sequence of
	real-valued random variables such that
	\begin{equation*}
		1_{B}\sum_{n\geq 0} \sup_{y\in\mathbb{R}} \left\vert \mathbb{P}[T_n \leq y \vert \mathcal{G}_n] - \Phi(y)\right\vert < \infty \quad \mathbb{P}\mbox{-a.s.}
		\end{equation*}
	 Then
	$$\limsup_{n\to\infty} \frac{T_n}{\sqrt{2\log n}} \leq 1 \quad \P(\cdot| B)\mbox{-a.s.}$$
	If, furthermore, there exists
 a constant $k\geq 1$
such that $T_n$ is $\mathcal{G}_{n+k}$-measurable for each $n =0,1,...$, then
	$$\limsup_{n\to\infty} \frac{T_n}{\sqrt{2\log n}} = 1 \quad \P(\cdot| B)\mbox{-a.s.}$$
\end{lemma}
\textbf{Proof:} From \cite[p. 430, 1.5]{AH}, for any
sequence $B_n$ of events
and any  filtration $\mathcal{G}_n$,
\begin{align}\label{B-C-1}
	\left\{B_n,\  \mbox{i.o.} \right\} \subset\left\{\sum_{n=1}^\infty \mathbb{P}\left(B_n \mid \mathcal{G}_n\right)=\infty \right\}
\end{align}
and the two events above are  $\P$-a.s. equal if there exists
 a constant $k\geq 1$
 such that $B_{n}\in \mathcal{G}_{n+k}$ for all $n$.
Thus,
\begin{align}\label{B-C-2}
	&   \left\{B_n\cap B,\  \mbox{i.o.} \right\}=B \cap \left\{B_n,\  \mbox{i.o.} \right\} \nonumber\\
	& \subset B\cap \left\{\sum_{n=1}^\infty \mathbb{P}\left(B_n \mid \mathcal{G}_n\right)=\infty \right\}= \left\{1_B\sum_{n=1}^\infty \mathbb{P}\left(B_n \mid \mathcal{G}_n\right)=\infty \right\}.
\end{align}
Applying this fact to $B_n= \{T_n > (1+\eta)\sqrt{2\log n}\}$ and noting that for any $\eta>0$,
\[\sum_{n=1}^\infty \left(1-\Phi((1+\eta)\sqrt{2\log n})\right)<\infty,
\]
we conclude that $\P\left(B \cap \left\{B_n,\  \mbox{i.o.} \right\} \right)=0$, which implies
the first result. For the second result, let $B_n= \{T_n > (1-\eta)\sqrt{2\log n}\}$,
then according to the fact that $B_n\in \mathcal{G}_{n+k}$, we have
\[
\left\{B\cap B_n,\  \mbox{i.o.} \right\} =\left\{1_{B}\sum_{n=1}^\infty \mathbb{P}\left(B_n \mid \mathcal{G}_n\right)=\infty \right\}.
\]
Noticing that $\sum_{n=1}^\infty\left(1- \Phi((1-\eta)\sqrt{2\log n})\right)=\infty$ for any $\eta >0$, we conclude that
\[
B= \left\{1_B\sum_{n=1}^\infty \mathbb{P}\left(B_n \mid \mathcal{G}_n\right)=\infty \right\}=\left\{B\cap B_n,\  \mbox{i.o.} \right\},
\]
which implies the desired result.
\hfill$\Box$

\section{Proof of Theorems \ref{thm2}  and \ref{thm4}}\label{ss: 4}

In this section, we always assume that {\bf(H1)}---{\bf(H3)} hold.

\subsection{General theory}\label{ss: 4.1}

Combining the branching property and
the property that $D_k(t)D_k(s)=D_k(t+s)$, we get that
for any $f\in L^2(E,\mu)\cap L^4(E,\mu)$ and $0<r< s\leq \infty$,
\begin{align}\label{Decomposition}
		& \langle f, X_{t+r} \rangle - \langle T_r (f_{sm}+f_{cr}), X_{t}\rangle -\sum_{2\mathfrak{R}_k < \lambda_1} e^{-\lambda_k (t+r)} H_{t+s}^{(k)}D_k(t+r)\langle f, \widehat{\Phi}_k\rangle_\mu \nonumber \\
	& = \sum_{i=1}^{M_{t}} \left[\langle f, X_r^i\rangle  - T_r(f_{sm}+f_{cr})(X_t(i))  -\sum_{2\mathfrak{R}_k < \lambda_1} e^{-\lambda_k r}H_s^{(k),i}D_k(r)\langle f, \widehat{\Phi}_k\rangle_\mu  \right].
\end{align}
Here $M_{t}$ is the number of particles alive at time $t$. For $i=1, \dots, M_t$, $X_{t}(i)$ is the position of the $i$-th particle,
and $\left(X_r^i,  H_s^{(k),i}\right)$ has the same distribution as  $\left(X_r, H_s^{(k)}\right)$ under $\mathbb{P}_{\delta_{X_{t}(i)} }$.
Furthermore, by the branching property, the random variables $H_s^{(k),i}$ are independent conditioned on
$\mathcal{F}_{t}:=\sigma\left( X_s: s\leq t \right)$.

For $0<r<s\leq \infty$,
 we define for $i=1, \dots, M_t$,
\begin{align}
	 &
	 Y_{t}^{f,i}(s,r) :=\langle f, X_r^i\rangle  - T_r(f_{sm}+f_{cr})(X_t(i))  -\sum_{2\mathfrak{R}_k < \lambda_1} e^{-\lambda_k r}H_s^{(k),i}D_k(r)\langle f, \widehat{\Phi}_k\rangle_\mu,
\label{def-Y(f,i)}
\\
	 &
	 Z_{t}^{f,i}(s,r) := Y_t^{f,i}(s, r)1_{\{|Y_t^{f,i}(s, r)|\leq e^{-\lambda_1 t/2}\}} ,\\
&
	 U_{t}^f(s,r) :=\sum_{i=1}^{M_{t}} \left(Z_t^{f,i}(s, r) - \mathbb{E}_{\delta_x}\left[Z_t^{f,i}(s, r)\big\vert \mathcal{F}_{t}\right]\right).\label{def-U^f_t}
\end{align}
Note that, for $i=1,\cdots, M_t$, $Y_{t}^{f, i}(s, r)$, $Z_{t}^{f, i}(s, r)$ and $U_{t}^f(s,r)$ contain information about the branching Markov process after time $t$ and therefore are not in $\mathcal{F}_{t}$. Note also that
$\mathbb{E}_{\delta_x}\left[Y_{t}^{f, i}(s, r) \big\vert \mathcal{F}_{t}\right] = 0$ and hence $\mathbb{E}_{\delta_x}\left[Z_t^{f,i}(s,r)\big\vert  \mathcal{F}_{t}\right] =-\mathbb{E}_{\delta_x}\left[Y_{t}^{f, i}(s, r) 1_{\{ |Y_t^{f,i}(s, r)|> e^{-\lambda_1 t/2}\}}\big\vert \mathcal{F}_{t}\right] .$

For $0<r< s\leq \infty$, we define
\begin{align}\label{V-s-r}
	&
	Y^f(s,r):= \langle f, X_r\rangle  - T_r(f_{sm}+f_{cr})(x)  -\sum_{2\mathfrak{R}_k < \lambda_1} e^{-\lambda_k r}H_s^{(k)}D_k(r)\langle f, \widehat{\Phi}_k\rangle_\mu,\nonumber\\
	&
	V_{s,r}^f(x):=\textup{Var}_x \left(Y^f(s, r)\right)= \mathbb{E}_{\delta_x} \left((Y^f(s, r))^2\right).
\end{align}
From the definitions of $\mathcal{C}_{la}, \mathcal{C}_{cr} $ and $\mathcal{C}_{sm}$,
$Y_{t}^{f,i}(s,r)$, $Z_{t}^{f,i}(s,r)$, $U_{t}^f(s,r)$, $Y^f(s,r)$ and $V_{s,r}^f$ are all real-valued random variables.
It follows from Lemma \ref{lemma2} that, for any $f \in L^2(E,\mu)\cap L^4(E,\mu)$
and $0<r< s\leq \infty$,  $V_{s,r}^f \in L^2(E,\mu)$.
We claim that if $f$ is of the form \eqref{e:rsspecial}, then
for any $0<r<s\leq\infty,$
\begin{align}\label{Claim1}
	V^f_{s, r}\in L^2(E,\mu ) \cap L^4(E,\mu).
\end{align}
In fact, for any $k\in \mathbb{I}$ and $v\in \mathbb{C}^{n_k}$, define
\[
g_k(x):= (\Phi_k(x))^Tv\quad \mbox{and}\quad h_k(x):= e^{ 14 \lambda_k t_0} (\Phi_k(x))^TD_k(-14t_0) v.
\]
By Jensen's inequality, for any $f\in L^2(E,\mu)\cap L^4(E,\mu), p\geq 1$ and $t>0$, we have
\begin{align}\label{T-t-upp}
|T_t f|^p
\leq   e^{\Vert A^{(1)}\Vert_\infty (p-1) t} T_t(|f|^p).
\end{align}
  Combining Lemma \ref{Useful-known-result} (2) and  \eqref{T-t-upp},
   we get $|h_k|\lesssim_{t_0, v,k} b_{t_0}^{1/2} $  and
 \begin{align}\label{phi-is-l-8}
 	|g_k|^2 = |T_{14 t_0} h_k|^2 \lesssim_{t_0} T_{14 t_0} (|h_k|^2) \lesssim_{t_0,v,k} T_{14 t_0} (b_{t_0}).
 \end{align}
 Therefore,  it follows from
 \eqref{T-s-b-t-0} and Lemma  \ref{Useful-known-result} (1)
 that for any $R>0$ and any $r\in (0, R]$,
\begin{align}\label{upp-second-phi}
	  \mathbb{E}_{\delta_x}\large(\left|\langle g_k, X_r\rangle\right|^2\large) & \lesssim_{R} T_{r}(|g_k|^2) (x)  \lesssim_{v,k,t_0,R}
	  	 T_{r +14t_0}(b_{t_0}) (x) \lesssim_{v,k,t_0, R}  b_{t_0}^{1/2}(x) \land T_{r+11t_0}(b_{t_0}^{1/2})(x),\quad
\end{align}
which implies \eqref{Claim1} for $s<\infty$. The case  $s=\infty$ follows from Lemma \ref{lemma2} (3) and \eqref{upp-second-phi}.

 Note that for two
 real-valued random variables $Y_1$ and $Y_2$,
\begin{align}\label{step_3}
	|\textup{Var}(Y_1 + Y_2) - \textup{Var}(Y_1)| \leq \textup{Var}(Y_2) + 2\sqrt{\textup{Var}(Y_1) \textup{Var}(Y_2)}.
\end{align}
Therefore, by the definition of $V^f_{s, r}$, we have
\begin{align}\label{Limit-of-V}
	\lim_{s\to\infty} V_{s,r}^f= V_{\infty ,r}^f,\quad \forall r\in (0,\infty), x\in E.
\end{align}

	\begin{lemma}\label{lemma5}
	If $f$ is of form \eqref{e:rsspecial}, then  for any $0<r< s\leq \infty$ and $\delta>0$,
	\begin{align}\label{step_34}
		\lim_{n\to\infty} e^{\lambda_1 n\delta} \textup{Var}_x \left[U^f_{n\delta}(s, r) \big| \mathcal{F}_{n\delta} \right]= \langle V^f_{s, r},\widehat{\phi}_1 \rangle_\mu W_\infty,\quad \P_{\delta_x}\text{-a.s.}
	\end{align}
	\end{lemma}
	\textbf{Proof:}
	We first prove that
	\begin{align}\label{step_36}
			\lim_{n \to \infty} e^{\lambda_1 n\delta}\textup{Var}_x\left[\sum_{i=1}^{M_{n\delta}}Y_{n\delta}^{f,i}(s, r) \bigg| \mathcal{F}_{n\delta}\right] = \langle V^f_{s, r},\widehat{\phi}_1 \rangle_\mu W_\infty,\quad \mathbb{P}_{\delta_x}\mbox{-a.s.}
	\end{align}
	Note that,  conditioned on $\mathcal{F}_{n\delta}$,
	$\big\{Y_{n\delta}^{f,i}(s, r), i=1, \dots, M_{n\delta}\big\}$
	are independent. Thus,
	\begin{align*}
		e^{\lambda_1 n\delta}\textup{Var}_x\left[\sum_{i=1}^{M_{n\delta}}Y_{n\delta}^{f,i} (s, r)\bigg| \mathcal{F}_{n\delta}\right]= e^{\lambda_1 n\delta} \sum_{i=1}^{M_{n\delta}} V^f_{s, r}\left(X_{n\delta}(i)\right) = e^{\lambda_1 n\delta} \langle V^f_{s, r}, X_{n\delta}\rangle.
	\end{align*}
	Combining Lemma \ref{lemma10} and \eqref{Claim1}, we get \eqref{step_36}.

	Define
	$Y_1:= \sum_{i=1}^{M_{n\delta}}Y_{n\delta}^{f,i}(s,r)$ and $Y_2 :=   U_{n\delta}^f(s,r)	- \sum_{i=1}^{M_{n\delta}}Y_{n\delta}^{f,i}(s,r)$.
	By \eqref{step_3}, to prove \eqref{step_34}, it suffices to show that
	\begin{equation}\label{step_6}
		\lim_{n \to \infty} e^{\lambda_1 n\delta}\textup{Var}_x\left[\sum_{i=1}^{M_{n\delta}}Y_{n\delta}^{f,i}(s, r) -
		U_{n\delta}^f(s, r)
		\bigg| \mathcal{F}_{n\delta}\right] =0 \quad \mathbb{P}_{\delta_x}\textup{-a.s.}
	\end{equation}
	 Using the inequality $\mbox{Var}_x(X)\leq \mathbb{E}_{\delta_x}(X^2)$, we get
	\begin{align}\label{step_7}
			e^{\lambda_1 n\delta}\textup{Var}_x\left[\sum_{i=1}^{M_{n\delta}}Y_{n\delta}^{f,i}(s, r) -	U_{n\delta}^f(s, r)
			\bigg| \mathcal{F}_{n\delta}\right]  &  = e^{\lambda_1 n\delta} \sum_{i=1}^{M_{n\delta}} \textup{Var}_x\left[Y_{n\delta}^{f,i}(s, r) 1_{\{|Y_{n\delta}^{f,i}(s, r)| > e^{-\lambda_1 n\delta/2} \}}\Big| \mathcal{F}_{n\delta}\right]\nonumber \\ &
		\leq e^{\lambda_1 n\delta}\langle V^{f,n\delta}_{s, r}, X_{n\delta} \rangle ,
	\end{align}
	where for $A>0$,
	$$
	V^{f, A}_{s, r}(x) : = \mathbb{E}_{\delta_x} \left[\left(Y^f(s, r)\right)^2 1_{\{|Y^f(s, r)| > e^{-\lambda_1 A/2} \}}\right]\leq V^f_{s, r}(x).
	$$
	For any fixed $A >0$, if $n>A/\delta$, then we have $V^{ f,n\delta}_{s, r}  \leq V^{f,A}_{s, r}$.
	Applying Lemma \ref{lemma10} to $ V^{f,A}_{s, r}$, we get that $e^{\lambda_1n\delta }\langle V^{f,A}_{s, r}, X_{n\delta} \rangle$ converges to $\langle  V^{f,A}_{s,r}, \widehat{\phi}_1\rangle_\mu W_\infty \ \mathbb{P}_{\delta_x}$-a.s. Hence,
	$$
	\limsup_{n\to\infty} e^{\lambda_1 n\delta }\langle V^{f,n\delta}_{s, r}, X_{n\delta} \rangle  \leq \limsup_{n\to\infty} e^{\lambda_1 n\delta }\langle V^{f,A}_{s, r}, X_{n\delta} \rangle = \langle V^{f,A}_{s, r}, \widehat{\phi}_1\rangle_\mu W_\infty, \quad \mathbb{P}_{\delta_x}\mbox{-a.s.}$$
	Letting $A \to \infty $, together with \eqref{step_7}, we get \eqref{step_6} and this completes the proof of the lemma.
   \hfill$\Box$

	\begin{lemma}\label{lemma7}
		Let $f$ is of form \eqref{e:rsspecial} and $0<r< s\leq \infty$.
		For any $\delta>0$, define
		\begin{equation*}
			\Delta_{n\delta}^f(s,r) : = \sup_{y\in\mathbb{R}} \left\vert
			\mathbb{P}_{\delta_x}\left[\frac{U_{n\delta}^f(s, r)}{\sqrt{\textup{Var}_x\left[U_{n\delta}^f(s, r)\big| \mathcal{F}_{n\delta}\right]}} \leq y \bigg\vert \mathcal{F}_{n\delta}\right]   -\Phi(y) \right\vert.
		\end{equation*}
		Then $\mathbb{P}_{\delta_x}$-almost surely,
		\begin{equation}\label{step_23}
			1_{\mathcal{E}^c}\sum_{n\geq 1} \Delta^f_{n\delta}(s, r) < \infty.
		\end{equation}
	\end{lemma}
	\textbf{Proof:}
\emph{Step 1}: The goal of this step is to prove that
\begin{equation}\label{step_8}
	\sum_{n > 2t_0/\delta} e^{3\lambda_1 n\delta/2} \sum_{i=1}^{M_{n\delta}}\mathbb{E}_{\delta_x} \left[\left|Z_{n\delta}^{f,i}(s, r)\right|^3 \Big| \mathcal{F}_{n\delta}\right]  < \infty ,\quad \mathbb{P}_{\delta_x}\textup{-a.s.}
\end{equation}
It suffices to show that
\begin{equation}\label{step_8a}
	\mathbb{E}_{\delta_x} \left(\sum_{n > 2t_0/\delta}e^{3\lambda_1 n\delta/2} \sum_{i=1}^{M_{n\delta}}\mathbb{E}_{\delta_x} \left[\left|Z_{n\delta}^{f,i}(s, r)\right|^3 \Big| \mathcal{F}_{n\delta}\right] \right)  < \infty.
\end{equation}
Define
\begin{equation}\label{Def-of-g}
	g_{s,r}^{f,n\delta}(x):= \mathbb{E}_{\delta_x} \left(\left|Y^f(s, r)\right|^3 1_{\{|Y^f(s, r)| \leq e^{-\lambda_1 n\delta/2} \}}\right).
\end{equation}
Then
\begin{align}
	&	\mathbb{E}_{\delta_x} \left( \sum_{n > 2t_0/\delta} e^{3\lambda_1 n\delta/2} \sum_{i=1}^{M_{n\delta}}\mathbb{E}_{\delta_x} \left[\left|Z_{n\delta}^{f,i}(s, r)\right|^3 \Big| \mathcal{F}_{n\delta}\right]  \right) = \sum_{n > 2t_0/\delta} e^{3\lambda_1 n\delta/2} \mathbb{E}_{\delta_x} \left(\sum_{i=1}^{M_{n\delta}}\mathbb{E}_{\delta_x} \left[\left|Z_{n\delta}^{f,i}(s, r)\right|^3 \Big| \mathcal{F}_{n\delta}\right]  \right) \nonumber\\ &
	= \sum_{n > 2t_0/\delta} e^{3\lambda_1 n\delta/2}  \mathbb{E}_{\delta_x} \left(\sum_{i=1}^{M_{n\delta}} g^{f,n\delta}_{s, r} \left(X_{n\delta}(i) \right) \right) = \sum_{n > 2t_0/\delta} e^{3\lambda_1 n\delta/2} T_{n\delta} g^{f,n\delta}_{s, r}(x).
\label{e:step_8b}\end{align}
Recall that, for any $f\in L^2(E, \mu; \mathbb{C})$, $\widetilde{f}$ is defined in \eqref{step_4}.
Fix  $a\in (\lambda_1, \mathfrak{R}_2)$, by Lemma \ref{lemma3} (1),
\begin{equation}\label{step_9}
        \left|T_{n\delta} \widetilde{g^{f,n\delta}_{s, r}}(x)\right|
	\lesssim_{a, t_0} e^{-a {n\delta}}\big\Vert g^{f,n\delta}_{s, r} \big\Vert_2  b_{t_0}^{1/2}(x),\quad n> \frac{2t_0}{\delta}, x\in E.
\end{equation}
Using the definition of $g^{f,n\delta}_{s,r}$, it is easy to see that
\begin{align} \label{upp-of-g-s-r}
g^{f,n\delta}_{s, r}(x)
\leq e^{-\lambda_1 n\delta/2} \mathbb{E}_{\delta_x} \left(\left|Y^f(s, r)\right|^2\right) =  e^{-\lambda_1 n\delta/2} V^f_{s, r}(x).
\end{align}
Plugging the inequality above into  \eqref{step_9} and applying \eqref{Claim1}, we get that
\begin{equation*}
        \left|T_{n\delta} \widetilde{g^{f,n\delta}_{s, r}}(x)\right|
\lesssim_{a, t_0} e^{-a n\delta }e^{-\lambda_1 n\delta /2} \Vert V^f_{s, r} \Vert_2 b_{t_0}^{1/2}(x), \quad n> \frac{2t_0}{\delta}, x\in E.
\end{equation*}
Therefore,
\begin{align}\label{step_11}
	\sum_{n > 2t_0/\delta} e^{3\lambda_1 n\delta /2}
\left|T_{n\delta} \widetilde{g^{f,n\delta}_{s, r}}(x)\right|
& \lesssim_{a, t_0} b_{t_0}^{1/2}(x)\sum_{n > 2t_0/\delta} e^{3\lambda_1 n\delta/2}e^{-(a+\lambda_1/2) n\delta}  \nonumber \\ & =b_{t_0}^{1/2}(x)\sum_{n > 2t_0/\delta} e^{(\lambda_1 - a) n\delta} < \infty.
\end{align}
We claim that
\begin{equation}\label{Goal}
	\sum_{n \geq 1} e^{3\lambda_1 n\delta /2}
\left|T_{n\delta } \left({g}^{f,n\delta}_{s, r} -\widetilde{g^{f,n\delta}_{s, r}}\right)(x) \right|
 = \sum_{n \geq 1} e^{\lambda_1 n\delta/2}\langle {g}^{f,n\delta}_{s, r}, \widehat{\phi}_1 \rangle_\mu \phi_1(x) < \infty.
\end{equation}
In fact, combining Lemma \ref{lemma11} (with $Y= Y^f(s, r)$) and the definition of $g^{f,n\delta}$ in \eqref{Def-of-g}, we get that
\begin{align*}
	\sum_{n \geq 1} e^{\lambda_1 n\delta/2}{g}^{f,n\delta}_{s, r}(x) \lesssim_{\delta} V^f_{s, r}(x) .
\end{align*}
 Since $V^f_{s,r}(x)$ and $\widehat{\phi}_1(x)$ both belong to $L^2(E ,\mu)$, we have  $\langle V^f_{s,r} ,\widehat{\phi}_1 \rangle_\mu < \infty$.
Now \eqref{Goal} follows from Fubini's theorem.
Combining \eqref{e:step_8b},  \eqref{step_11} and \eqref{Goal}, we get \eqref{step_8a}.

\emph{Step  2}: In this step, we  prove the conclusion of the lemma.
It is trivial that $	\Delta^f_{n\delta}(s, r) \leq 2$. Since $\{M_{n\delta} > 0 \}\in \mathcal{F}_{n\delta}$,  by Lemma \ref{lemma6},  under $\P_{\delta_x},$ on the event $\{M_{n\delta}>0 \},$
\begin{align}\label{step_13}
	\Delta^f_{n\delta}(s, r) &\lesssim \frac{\sum_{i=1}^{M_{n\delta}}\mathbb{E}_{\delta_x} \left[\left|Z_{n\delta}^{f,i}(s, r) - \mathbb{E}_{\delta_x}\left[Z_{n\delta}^{f,i}(s, r) \big\vert \mathcal{F}_{n\delta}\right] \right|^3 \Big| \mathcal{F}_{n\delta}\right] }{\sqrt{\left(\textup{Var}_x\left[U^f_{n\delta}(s, r)\big| \mathcal{F}_{n\delta} \right]\right)^3}}\nonumber\\
	& \lesssim \frac{\sum_{i=1}^{M_{n\delta}}\mathbb{E}_{\delta_x} \left[\left|Z_{n\delta}^{f,i}(s, r)\right|^3 \Big| \mathcal{F}_{n\delta}\right] }{\sqrt{\left(\textup{Var}_x\left[U^f_{n\delta}(s, r)\big| \mathcal{F}_{n\delta}\right]\right)^3}},
\end{align}
where in the second inequality, we used the inequality $\mathbb{E}|Y-\mathbb{E}Y|^3 \leq 8\mathbb{E}|Y|^3$ for any
real-valued $Y$ with $\mathbb{E}|Y|^3 < \infty.$
Since $\mathcal{E}^c \subset \{M_{n\delta}>0 \}$,
\eqref{step_13} holds on the event $\mathcal{E}^c$ under $\P_{\delta_x}$.
Now suppose that $\Omega_0$ is an event with $\P_{\delta_x}\left(\Omega_0\right) = 1$ such that, for any $\omega \in \Omega_0$, the assertion of Lemma \ref{lemma5}, \eqref{step_8} and \eqref{step_13} hold.
Then for each $\omega \in \Omega_0\cap \mathcal{E}^c,$ there exists a large integer $N = N(\omega)>2t_0/\delta$ such that for $n \geq N,$
$$ \textup{Var}_x\left[
U_{n\delta}^f(s, r)
\big| \mathcal{F}_{n\delta}\right](\omega) \geq \frac{e^{-\lambda_1 n\delta}}{2}\langle V^f_{s, r},\widehat{\phi}_1 \rangle_\mu  W_\infty(\omega)>0.$$
Therefore, on $\Omega_0\cap\mathcal{E}^c$, by \eqref{step_13},
\begin{equation*}
	\sum_{n\geq 1} \Delta^f_{n\delta}(s, r) \lesssim \left(1+N\right) +
\frac{	\sqrt{8}}{\sqrt{\left[\langle V_{s,r}^f,\widehat{\phi}_1 \rangle_\mu  W_\infty\right]^3} }\sum_{n \geq N} e^{3\lambda_1 n\delta/2} \sum_{i=1}^{M_{n\delta}}\mathbb{E}_{\delta_x} \left[\left|Z_{n\delta}^{f,i}(s, r)\right|^3 \Big| \mathcal{F}_{n\delta}\right] .
\end{equation*}
Applying \eqref{step_8}, we get that \eqref{step_23} holds $\mathbb{P}_{\delta_x}$-almost surely.
\hfill$\Box$

Now we are going to prove an LIL for  $\sum_{i=1}^{M_{t}} Y_{t}^{f,i}(s, r)$
for functions of  the form \eqref{e:rsspecial}. We first deal with discrete times $\{n\delta, n\in \N\}$ for any given $\delta>0$,
then we prove the continuous-time LIL.
The argument for discrete-time is inspired by \cite{IK} and the argument for continuous time is inspired by \cite[Section 12]{AH}
(for example, see   the proof of \cite[Theorem 12.4, p.340]{AH}) and  \cite[p.20--p.22]{IKM}.

\begin{lemma}\label{lemma12}
	 If  $f$ is of form \eqref{e:rsspecial}, then for any
	 $0<r<s\leq \infty$ and $\delta>0$,
	\begin{align}\label{LIMSUP3}
		 \limsup_{n\to\infty}/\liminf_{n \to \infty} \frac{e^{\lambda_1 n\delta/2}
		 			 	\sum_{i=1}^{M_{n\delta}} Y_{n\delta}^{f,i}(s, r)
		 			 	}{\sqrt{2\log (n\delta)}} =+/- \sqrt{\langle V^f_{s, r} ,\widehat{\phi}_1 \rangle_\mu W_\infty},
		\quad \mathbb{P}_{\delta_x}\left(\cdot| \mathcal{E}^c\right)
			\mbox{-a.s.}
	\end{align}
\end{lemma}
\textbf{Proof: }
We only prove the $\limsup$ assertion. The proof of the $\liminf$ assertion is  similar.

\emph{Step 1}. In this step, we prove that for any $0<r< s\leq \infty$,
\begin{align}\label{Goal-1}
	\lim_{n\to\infty} \left|e^{\lambda_1 n\delta/2}
	\sum_{i=1}^{M_{n\delta}} Y_{n\delta}^{f,i}(s, r)
		 - e^{\lambda_1 n\delta/2}
	U_{n\delta}^f(s, r)\right|
	 = 0,\quad\mathbb{P}_{\delta_x}\left(\cdot |\mathcal{E}^c\right) \mbox{ a.s.}
\end{align}
Note that
\begin{align*}
	 &	e^{\lambda_1n\delta/2}
		 \sum_{i=1}^{M_{n\delta}} Y_{n\delta}^{f,i}(s, r)
		  - e^{\lambda_1 n\delta/2}
	U_{n\delta}^f(s, r)
	\nonumber\\ & = e^{\lambda_1 n\delta /2} \left(\sum_{i=1}^{M_{n\delta}} Y_{n\delta}^{f,i}(s, r)1_{\{|Y_{n\delta}^{f,i}(s, r) | >e^{-\lambda_1 n\delta/2}\}} -\mathbb{E}_{\delta_x}\left[Y_{n\delta}^{f,i}(s, r) 1_{\{| Y_{n\delta}^{f,i}(s, r) | >e^{-\lambda_1 n\delta /2}\}} \Big| \mathcal{F}_{n\delta} \right] \right).
\end{align*}
Using the inequality $\left|\mathbb{E}[Y|\mathcal{F}]\right| \leq \mathbb{E}[|Y| \big| \mathcal{F}],$ we get that
\begin{align*}
	 \mathbb{E}_{\delta_x}\left|e^{\lambda_1 n\delta/2}
	 	 \sum_{i=1}^{M_{n\delta}} Y_{n\delta}^{f,i}(s, r)
		 - e^{\lambda_1 n\delta/2}	U_{n\delta}^f(s, r)\right|
	\leq 2 e^{\lambda_1 n\delta/2} \mathbb{E}_{\delta_x} \left[\sum_{i=1}^{M_{n\delta}} \left|Y_{n\delta}^{f,i}(s, r)\right|1_{\{| Y_{n\delta}^{f,i}|(s, r) >e^{-\lambda_1 n\delta/2}\}} \right].
\end{align*}
Therefore, to prove \eqref{Goal-1}, we only need to show that
\begin{equation}\label{step_17}
	\sum_{n > 2t_0/\delta} e^{\lambda_1 n\delta/2} \mathbb{E}_{\delta_x} \left[\sum_{i=1}^{M_{n\delta}} \left|Y_{n\delta}^{f,i}(s, r)\right|1_{\{| Y_{n\delta}^{f,i}(s, r)| >e^{-\lambda_1 n\delta/2}\}} \right] < \infty.
\end{equation}
Define
\begin{equation}\label{step_27}
	l_{s,r}^{f,n\delta}(x) := \mathbb{E}_{\delta_x} \left(\left|Y^f(s, r)\right| 1_{\{|Y^f(s, r)| > e^{-\lambda_1 n\delta/2}\}} \right),
\end{equation}
then $l^{f,n\delta}_{s, r}(x)\leq e^{\lambda_1 n\delta/2} V^f_{s, r}(x)$ for any $ n \in \N$ and $x\in E.$
Fix any $a\in (\lambda_1, \mathfrak{R}_2)$, then by Lemma \ref{lemma3} (1), we have
\begin{align*}
\left|T_{n\delta} \widetilde{l^{f,n\delta}_{s, r}}(x) \right|
\lesssim_{a,\delta, t_0} e^{-a n\delta} e^{\lambda_1 n\delta/2} \Vert V^f_{s, r} \Vert_2 b_{t_0}^{1/2}(x),\quad n> \frac{2t_0}{\delta}, x\in E.
\end{align*}
Thus,
\begin{equation}\label{step_15}
	\sum_{n > 2t_0/\delta} e^{\lambda_1 n\delta /2}	
\left|T_{n\delta} \widetilde{l^{f,n\delta}_{s, r}}(x)\right|
\lesssim_{a, \delta, t_0} \Vert V^f_{s, r} \Vert_2 b_{t_0}^{1/2}(x) \sum_{n > 2t_0/\delta}
 e^{(\lambda_1-a) n\delta}  < \infty.
\end{equation}
Since $T_{n\delta}\left(l^{f,n\delta}_{s,r} -\widetilde{l^{f,n\delta}_{s, r}}\right)(x)
= e^{-\lambda_1 n\delta}\langle l^{f,n\delta}_{s, r} , \widehat{\phi}_1\rangle_\mu\phi_1(x)$, by Lemma \ref{lemma11} (with  $Y= Y^f(s,r)$), we get that
\begin{align}\label{step_16}
	\sum_{n > 2t_0/\delta} e^{\lambda_1 n\delta/2} e^{-\lambda_1 n\delta}l^{f, n\delta}_{s, r}	\lesssim_{t_0, \delta} V^f_{s, r} (x) < \infty.	
\end{align}
Combining \eqref{step_15} and \eqref{step_16}, we obtain that
\begin{align}\label{step_28}
	&\sum_{n >2t_0/\delta} e^{\lambda_1 n\delta/2} \mathbb{E}_{\delta_x} \left[\sum_{i=1}^{M_{n\delta}} \left|Y_{n\delta}^{f,i}(s, r)\right|1_{\{| Y_{n\delta}^{f,i}(s, r)| >e^{-\lambda_1 n\delta/2}\}} \right]
 =  \sum_{n >2t_0/\delta} e^{\lambda_1 n\delta/2}
T_{n\delta} {l}^{f,n\delta}_{s, r}(x)
	\nonumber \\ &
	\lesssim_{\delta,a, t_0}
	 \Vert V^f_{s, r} \Vert_2 b_{t_0}^{1/2}(x) \sum_{n >2t_0/\delta} e^{(\lambda_1-
	a) n\delta} + \langle V^f_{s, r}, \widehat{\phi}_1\rangle_\mu \phi_1(x)<\infty,
\end{align}
which implies \eqref{Goal-1}.

\emph{Step 2}: In this step, we prove the  assertion of lemma for  $s\in (r,\infty)$.
 Combining Lemma \ref{lemma1} (with $B = \mathcal{E}^c$) and Lemma \ref{lemma7}, we get that, for $s\in (r,\infty)$,
\begin{equation*}
	\limsup_{n\to\infty} \frac{U^f_{n\delta}(s, r)}{\sqrt{2\log n \textup{Var}_x\left[U^f_{n\delta}(s, r)\big| \mathcal{F}_{n\delta}\right]}} = 1,\quad \mathbb{P}_{\delta_x}\left(\cdot \vert \mathcal{E}^c\right)\textup{-a.s.}
\end{equation*}
Noticing that $\lim_{n\to\infty} \log (n\delta) / \log n =1$, by Lemma \ref{lemma5}, we have
\begin{equation}\label{LIMSUP1}
	\limsup_{n\to\infty} \frac{e^{\lambda_1 n\delta /2}U^f_{n\delta}(s, r)}{\sqrt{2\log (n\delta) }} = \sqrt{\langle V^f_{s, r}, \widehat{\phi}_1\rangle_\mu W_\infty},\quad \mathbb{P}_{\delta_x}\left(\cdot \vert \mathcal{E}^c\right)\textup{-a.s.}
\end{equation}
Now combining \eqref{Goal-1} and \eqref{LIMSUP1}, we get the desired result for $s\in (r,\infty)$.

\emph{Step 3}:  In this step, we prove the assertion of the  lemma for  $s=\infty$.
Combining Lemma \ref{lemma1} and Lemma \ref{lemma7}, we get
\begin{equation}\label{LIMSUP2}
	\limsup_{n\to\infty} \frac{e^{\lambda_1 n\delta /2}U_{n\delta}^f(\infty, r)}{\sqrt{2\log (n\delta) }}\leq \sqrt{\langle V_{\infty, r}^f, \widehat{\phi}_1\rangle_\mu W_\infty},\quad \mathbb{P}_{\delta_x}\left(\cdot \vert \mathcal{E}^c\right)\textup{-a.s.}
\end{equation}
Together with \eqref{Goal-1}, we obtain
\begin{equation}\label{LIMSUP4}
	\limsup_{n\to\infty} \frac{e^{\lambda_1 n\delta/2}
			\sum_{i=1}^{M_{n\delta}} Y_{n\delta}^{f,i}(\infty,r)
			}{\sqrt{2\log (n\delta)}} \leq \sqrt{\langle V_{\infty ,r}^f, \widehat{\phi}_1\rangle_\mu W_\infty},\quad \mathbb{P}_{\delta_x}\left(\cdot| \mathcal{E}^c\right)\textup{-a.s.}
\end{equation}
Using the same argument with  $U^f_{n\delta}(\infty, r)$ replaced to $-U^f_{n\delta}(\infty, r)$, we also see that
\begin{equation}\label{LIMSUP5}
	\liminf_{n\to\infty} \frac{e^{\lambda_1 n\delta/2}
			\sum_{i=1}^{M_{n\delta}} Y_{n\delta}^{f,i}(\infty,r)}{\sqrt{2\log (n\delta)}} \geq - \sqrt{\langle V_{\infty ,r}^f, \widehat{\phi}_1\rangle_\mu W_\infty},\quad \mathbb{P}_{\delta_x}\left(\cdot| \mathcal{E}^c\right)\textup{-a.s.}
\end{equation}
Note that for $s= \ell \delta+r , \ell\in \N$, it holds that
\begin{align*}
	&\frac{e^{\lambda_1 n\delta/2}
			\sum_{i=1}^{M_{n\delta}} Y_{n\delta}^{f,i}(\infty,r)}{\sqrt{2\log (n\delta)}} \\ &= \frac{e^{\lambda_1 n\delta/2}
				\sum_{i=1}^{M_{n\delta}} Y_{n\delta}^{f,i}(\ell \delta+r  ,r)}{\sqrt{2\log (n\delta)}} \\
	&\qquad + \sum_{2\mathfrak{R}_k< \lambda_1} e^{\lambda_k \ell \delta } \frac{e^{\lambda_1 n\delta /2}\sum_{i=1}^{M_{(n+\ell)\delta}}\left( \langle \Phi_k^T, X_r^{i} \rangle- e^{-\lambda_k r}H_{\infty}^{(k),i}D_k(r)\right)D_k(\ell\delta)^{-1} \langle f,\widehat{\Phi}_k\rangle_\mu}{\sqrt{2\log (n\delta)}}\\
	& = \frac{e^{\lambda_1 n\delta/2}
				\sum_{i=1}^{M_{n\delta}} Y_{n\delta}^{f,i}(\ell \delta+r  ,r)}{\sqrt{2\log (n\delta)}} \\
	&\quad + \sum_{2\mathfrak{R}_k< \lambda_1} e^{\lambda_k \ell \delta} \sum_{j=1}^{n_k} \frac{e^{\lambda_1 n\delta /2}\sum_{i=1}^{M_{(n+\ell)\delta}} Y_{(n+\ell)\delta}^{\phi_j^{(k)}, i} (\infty, r) \left( D_k(\ell\delta)^{-1} \langle f,\widehat{\Phi}_k\rangle_\mu\right)_j}{\sqrt{2\log (n\delta)}},
\end{align*}
where we use the notation $(v)_j=v_j$ for any vector $v=(v_1, v_2,..., v_{n_k})^T\in \mathbb{C}^{n_k}$.
Using the inequality
$$ \limsup_{n\to\infty}
\sum_{i=1}^p x_n^i
\geq \limsup_{n\to\infty}
x_n^1
+\sum_{i=2}^p \liminf_{n\to\infty}
x_n^i $$
and applying \eqref{LIMSUP3} to $ Y_{n\delta}^{f,i}(\ell \delta+r,r)$ and \eqref{LIMSUP5} to $Y_{(n+\ell)\delta}^{\phi_j^{(k)}, i} (\infty, r)$, we conclude that $\mathbb{P}_{\delta_x}\left(\cdot| \mathcal{E}^c\right)$-almost surely,
\begin{align*}
	&\limsup_{n\to\infty} \frac{e^{\lambda_1 n\delta/2}
			\sum_{i=1}^{M_{n\delta}} Y_{n\delta}^{f,i}(\infty,r)}{\sqrt{2\log (n\delta)}} \\ &\geq  \sqrt{\langle V_{\ell \delta +r,r}^f ,\widehat{\phi}_1 \rangle W_\infty}-\sum_{2\mathfrak{R}_k< \lambda_1} \sum_{j=1}^{n_k} \left|\left( D_k(\ell\delta)^{-1} \langle f,\widehat{\Phi}_k\rangle_\mu\right)_j\right| e^{-(\lambda_1/2 -\mathfrak{R}_k) \ell \delta }\sqrt{\langle V_{\infty,r}^{\phi_j^{(k)}} ,
		\widehat{\phi}_1\rangle_\mu
		 W_\infty} .
\end{align*}
It follows from \eqref{Limit-of-V} that $V_{\ell\delta+r ,r}^f(x) $ converges to $V_{\infty,r}^f (x) $.  Letting $\ell\to \infty$ in the above inequality and noticing that $|D_k(\ell\delta)^{-1} | $ is of polynomial growth, we get that
$$
\limsup_{n\to\infty} \frac{e^{\lambda_1 n\delta/2}
		\sum_{i=1}^{M_{n\delta}} Y_{n\delta}^{f,i}(\infty,r)}{\sqrt{2\log (n\delta)}}  \geq \sqrt{\langle V_{\infty  ,r}^f ,\widehat{\phi}_1 \rangle_\mu W_\infty} ,\quad \mathbb{P}_{\delta_x}\left(\cdot| \mathcal{E}^c\right)\textup{-a.s.}
$$
Combining the above with \eqref{LIMSUP4}, we get that \eqref{LIMSUP3} holds for $s=\infty$.
The proof is complete.

\hfill$\Box$

Now we are ready to treat the continuous-time case.

\begin{lemma}\label{lem:conti-mart}
Assume $g_k(x):= \left(\Phi_k(x)\right)^T v$ for some $k\in \mathbb{I}$ and $v\in \mathbb{C}^{n_k}$.
	Then
	\begin{align}
		\lim_{\delta \to 0} \limsup_{n\to\infty} \sup_{t\in [n\delta, (n+1)\delta]} \frac{e^{\lambda_1 n\delta/2}\left| \langle g_k, X_t \rangle - \langle T_{t-n\delta}g_k, X_{n\delta} \rangle\right|}{\sqrt{2\log (n\delta)}}=0,\quad  \mathbb{P}_{\delta_x}\left(\cdot| \mathcal{E}^c\right)
			\mbox{-a.s.}
	\end{align}
\end{lemma}
\textbf{Proof: }
\emph{Step 1:} We deal with discrete times in this step.
Note that $\gamma(\mathfrak{R}(g_k)) =\gamma(\mathfrak{I}(g_k))
=\gamma(\mathfrak{R}(T_\delta g_k)) =\gamma(\mathfrak{I}(T_\delta g_k))=\mathfrak{R}_k$.
When $2\mathfrak{R}_k\geq \lambda_1$,
using \eqref{Decomposition} for $g_k$ with $t=n\delta$ and $r=\delta$, and applying
Lemma \ref{lemma12} for $f=\mathfrak{R}(g_k)$ with  $r=\delta$
and $s=2\delta$,
 we get
that $\mathbb{P}_{\delta_x}\left(\cdot| \mathcal{E}^c\right)$-almost surely,
\begin{align}
	&\limsup_{n\to\infty} \frac{e^{\lambda_1 n\delta/2}\left|\mathfrak{R}\left(\langle g_k, X_{(n+1)\delta} \rangle - \langle T_{\delta}g_k, X_{n\delta} \rangle \right)\right|}{\sqrt{2\log (n\delta)}}
	= \sqrt{\langle V^{\mathfrak{R}(g_k)}_{2\delta, \delta} ,\widehat{\phi}_1 \rangle_\mu W_\infty}.
\end{align}
Similarly, applying Lemma \ref{lemma12} with $r=\delta$ and $f=\mathfrak{I}(g_k)$,
we also have $\mathbb{P}_{\delta_x}\left(\cdot| \mathcal{E}^c\right)$-almost surely,
\begin{align}
	&\limsup_{n\to\infty} \frac{e^{\lambda_1 n\delta/2}\left|\mathfrak{I}\left(\langle g_k, X_{(n+1)\delta} \rangle - \langle T_{\delta}g_k, X_{n\delta} \rangle \right)\right|}{\sqrt{2\log (n\delta)}}
	=
	 \sqrt{\langle V^{	\mathfrak{I}(g_k)	}_{2\delta, \delta} ,\widehat{\phi}_1 \rangle_\mu W_\infty}.
\end{align}
Therefore, when $2\mathfrak{R}_k\geq \lambda_1$, we conclude that
\begin{align}\label{e2}
	&\limsup_{n\to\infty} \frac{e^{\lambda_1 n\delta/2}\left|\langle g_k, X_{(n+1)\delta} \rangle - \langle T_{\delta}g_k, X_{n\delta} \rangle \right|}{\sqrt{2\log (n\delta)}} \nonumber\\
	& \leq
	 \sqrt{\langle V^{	\mathfrak{R}(g_k)
	}_{2\delta, \delta} ,\widehat{\phi}_1 \rangle_\mu W_\infty}+
	 \sqrt{\langle V^{\mathfrak{I}(g_k)
	}_{2\delta, \delta} ,\widehat{\phi}_1 \rangle_\mu W_\infty}
  =: \Gamma_\delta(g_k) \sqrt{W_\infty},\quad
   \mathbb{P}_{\delta_x}\left(\cdot| \mathcal{E}^c\right)\textup{-a.s.}
\end{align}
When $2\mathfrak{R}_k<\lambda_1$,
using \eqref{Decomposition} for $T_\delta g_k$ with $t=(n+1)\delta$ and $r=\delta$, applying Lemma \ref{lemma12} for $f=\mathfrak{R} (T_\delta g_k)$ with $ r=\delta$ and $s=2\delta $,
we get
\begin{align*}
	\limsup_{n\to\infty} \frac{e^{\lambda_1 n\delta/2}\left|\mathfrak{R}\left(\langle T_\delta g_k, X_{(n+1)\delta}\rangle - \langle g_k, X_{(n+2)\delta}\rangle \right) \right|}{\sqrt{2\log (n\delta)}} = \sqrt{\langle V^{\mathfrak{R} (T_\delta g_k) }_{2\delta, \delta} ,\widehat{\phi}_1 \rangle_\mu W_\infty},\quad \mathbb{P}_{\delta_x}\left(\cdot| \mathcal{E}^c\right)\textup{-a.s.}
\end{align*}
Similarly, we have
\begin{align*}
	\limsup_{n\to\infty} \frac{e^{\lambda_1 n\delta/2}\left|\mathfrak{I}\left(\langle T_\delta g_k, X_{(n+1)\delta}\rangle - \langle g_k, X_{(n+2)\delta}\rangle \right) \right|}{\sqrt{2\log (n\delta)}} = \sqrt{\langle V^{\mathfrak{I} (T_\delta g_k) }_{2\delta, \delta} ,\widehat{\phi}_1 \rangle_\mu W_\infty},\quad \mathbb{P}_{\delta_x}\left(\cdot| \mathcal{E}^c\right)\textup{-a.s.}
\end{align*}
Combining the two displays above, we get that, in the case  $2\mathfrak{R}_k<\lambda_1$, \eqref{e2} holds with
\[
\Gamma_\delta(g_k):= e^{\lambda_1\delta/2} \sqrt{\langle V^{\mathfrak{R} (T_\delta g_k) }_{2\delta, \delta} ,\widehat{\phi}_1 \rangle_\mu }  +e^{\lambda_1\delta/2} \sqrt{\langle V^{\mathfrak{I} (T_\delta g_k) }_{2\delta, \delta} ,\widehat{\phi}_1 \rangle_\mu }  .
\]
Define $W_t^{(k)}:= \langle T_{(n+1)\delta-t}g_k, X_t\rangle = \mathbb{E}_{\delta_x}\left( \langle g_k, X_{(n+1)\delta}\rangle \Big|\mathcal{F}_t\right)$ for  $t\in [n\delta, (n+1)\delta]$. Then $(W_t^{(k)}: t\in [n\delta, (n+1)\delta])$ is a martingale.
By \eqref{e2}, we have
\begin{align}\label{e1}
	&\limsup_{n\to\infty} \frac{e^{\lambda_1n\delta/2 }\left| W_{(n+1)\delta}^{(k)}- W_{n\delta}^{(k)}\right|}{\sqrt{2\log (n\delta)}}
	\leq \Gamma_\delta(g_k)\sqrt{W_\infty},\quad \mathbb{P}_{\delta_x}\left(\cdot| \mathcal{E}^c\right)\textup{-a.s.}
\end{align}
For $\rho >0$, define
\begin{align} \label{Def-of-epsilon-k}
	\epsilon_n (k,\delta) :=
	e^{-\lambda_1 n\delta/2}\sqrt{2\log (n\delta)}\left(\Gamma_\delta(g_k) \sqrt{W_{n\delta}}+	\rho\right).
\end{align}
By the second Borel-Cantelli lemma (see e.g. \cite[Theorem 5.3.2]{Durrett}),  we have
\begin{align*}
	& \left\{\left|W_{n\delta}^{(k)}-W_{(n+1)\delta}^{(k)}\right|> \epsilon_n(k,\delta), \textup{ i.o.} \right\}\\
	& = \left\{\sum_{n=1}^\infty \mathbb{P}_{\delta_x}\left(\left|W_{n\delta}^{(k)}-W_{(n+1)\delta}^{(k)}\right|> \epsilon_n(k,\delta) \Big|\mathcal{F}_{n\delta} \right) = +\infty\right\}.
\end{align*}
Combining this  with \eqref{e1}, we get that on $\mathcal{E}^c$, $\mathbb{P}_{\delta_x}$- almost surely,
\begin{equation}\label{step_20}
	\sum_{n=1}^\infty \mathbb{P}_{\delta_x}\left(\left|W_{n\delta}^{(k)}-W_{(n+1)\delta}^{(k)}\right|> \epsilon_n(k,\delta) \Big|\mathcal{F}_{n\delta} \right) < +\infty.
\end{equation}

\emph{Step 2:} Now we consider continuous time.
For any $t\in [n\delta, (n+1)\delta)$, define
\begin{align*}
	Z_t^{(k)}&:= \mathbb{E}_{\delta_x}\left[\left|W_{(n+1)\delta}^{(k)}-W_t^{(k)}\right|^2 \Big| \mathcal{F}_t\right],\quad
	B_n^{(k)} := \sup_{t\in[n\delta, (n+1)\delta)} \left[\left|W_{n\delta}^{(k)}-W_{t}^{(k)}\right|- \sqrt{2Z_t^{(k)}}\right],
	\\	\Gamma_n^{(k)} &:= \inf\left\{s\in [n\delta, (n+1)\delta): \left|W_{n\delta}^{(k)}-W_{s}^{(k)}\right|-\sqrt{2 Z_s^{(k)}}> \epsilon_n(k,\delta) \right\}\wedge ((n+1)\delta).
\end{align*}
We have
\begin{align}\label{step_19}
	&\mathbb{P}_{\delta_x}\left(\left|W_{n\delta}^{(k)}-W_{(n+1)\delta}^{(k)}\right|> \epsilon_n(k,\delta) \Big|\mathcal{F}_{n\delta}\right) \geq \mathbb{P}_{\delta_x}\left(\left|W_{n\delta}^{(k)}-W_{(n+1)\delta}^{(k)}\right|> \epsilon_n(k,\delta), \Gamma_n^{(k)} < (n+1)\delta  \Big|\mathcal{F}_{n\delta} \right)\nonumber \\
	&\geq  \mathbb{P}_{\delta_x}\left(\left| W_{\Gamma_n^{(k)}}^{(k)}-W_{(n+1)\delta}^{(k)} \right|< \sqrt{2Z_{\Gamma_n^{(k)}}^{(k)}}, \Gamma_n^{(k)} < (n+1)\delta \big|\mathcal{F}_{n\delta} \right)\nonumber \\
	&=  \mathbb{E}_{\delta_x}\left(\mathbb{P}_{\delta_x}\left(\left|W_{\Gamma_n^{(k)}}^{(k)}-W_{(n+1)\delta}^{(k)} \right|<\sqrt{2Z_{\Gamma_n^{(k)}}^{(k)}}\big| \mathcal{F}_{\Gamma_n^{(k)}}\right)1_{\left\{\Gamma_n^{(k)} < (n+1)\delta  \right\}} \Big|\mathcal{F}_{n\delta} \right).
\end{align}
By Markov's inequality and the strong Markov property, it is easy to see that
\begin{align}\label{Markov-Ineq}
	&\mathbb{P}_{\delta_x}\left(\left|W_{\Gamma_n^{(k)}}^{(k)}-W_{(n+1)\delta}^{(k)} \right|<\sqrt{2Z_{\Gamma_n^{(k)}}^{(k)}}\big| \mathcal{F}_{\Gamma_n^{(k)}}\right)
	=  1 - \mathbb{P}_{\delta_x}\left(\left|W_{\Gamma_n^{(k)}}^{(k)}-W_{(n+1)\delta}^{(k)} \right|\geq \sqrt{2Z_{\Gamma_n^{(k)}}^{(k)}}\big| \mathcal{F}_{\Gamma_n^{(k)}}\right) \nonumber\\
	& \geq 1- \mathbb{E}_{\delta_x}\left[\frac{\left|W_{(n+1)\delta}^{(k)}- W_{\Gamma_n^{(k)}}^{(k)}\right|^2 }{2Z_{\Gamma_n^{(k)}}^{(k)}}\bigg|\mathcal{F}_{\Gamma_n^{(k)}}\right]= \frac{1}{2}.
\end{align}
Therefore,
\begin{align}\label{step_21}
	&\mathbb{P}_{\delta_x}\left(\left|W_{n\delta}^{(k)}-W_{(n+1)\delta}^{(k)}\right|> \epsilon_n(k,\delta) \Big|\mathcal{F}_{n\delta}\right)  \geq \frac{1}{2} \mathbb{P}_{\delta_x}\left( \Gamma_n^{(k)} < (n+1)\delta \big|\mathcal{F}_{n\delta} \right)\nonumber\\
	&= \frac{1}{2} \mathbb{P}_{\delta_x}\left( B_n^{(k)} > \epsilon_n(k,\delta) \big|\mathcal{F}_{n\delta} \right).
\end{align}
Together with \eqref{step_20} and \eqref{step_21} we obtain that on $\mathcal{E}^c$,  $\mathbb{P}_{\delta_x}$-almost surely,
$$
\sum_{n=1}^\infty \mathbb{P}_{\delta_x}\left( B_n^{(k)} > \epsilon_n(k,\delta) \big|\mathcal{F}_{n\delta} \right) < +\infty.
$$
Since $\{B_n^{(k)} > \epsilon_n(k,\delta)\}\in \mathcal{F}_{(n+1)\delta}$, using the second Borel-Cantelli lemma again, we get that for any $\rho>0$ and $\delta>0$, $\mathbb{P}_{\delta_x}(\cdot|\mathcal{E}^c)$-almost surely,
\begin{align}\label{Step13-2}
	& \limsup_{n\to\infty} \sup_{t\in [n\delta, (n+1)\delta)} \frac{e^{\lambda_1 n\delta/2}\left|\langle T_{(n+1)\delta-t}g_k, X_{t} \rangle - \langle T_{\delta}g_k, X_{n\delta} \rangle \right|}{\sqrt{2\log (n\delta)}} \nonumber\\
	&\leq
	\Gamma_\delta (g_k)\sqrt{W_\infty}+\rho
	+ 	\limsup_{n\to\infty} \sup_{t\in [n\delta, (n+1)\delta)} \frac{\sqrt{2e^{\lambda_1n\delta}Z_t^{(k)}}}{\sqrt{2\log (n\delta)}}.
\end{align}
Using the inequality $\mbox{Var}(Y)\leq \mathbb{E}(Y^2)$, the branching property and \eqref{upp-second-phi} (with $R=1$), we have
\begin{align}\label{Z-T}
	e^{\lambda_1 t}Z_t^{(k)}& = e^{\lambda_1 t}\mbox{Var}_x\left[\langle g_k, X_{(n+1)\delta}\rangle \Big| \mathcal{F}_t\right] \leq e^{\lambda_1 t}\langle  \mathbb{E}_{\delta_\cdot}\left(\left|\langle g_k, X_{(n+1)\delta-t}\rangle\right|^2\right), X_t \rangle \nonumber\\
	& \lesssim_{v,k,t_0} e^{\lambda_1 ((n+1)\delta+11t_0)}\langle   T_{ (n+1)\delta -t+11t_0}(b_{t_0}^{1/2}), X_t \rangle \nonumber\\
	& = e^{\lambda_1 ((n+1)\delta+11t_0)} \mathbb{E}_{\delta_x}\left( \langle   \widetilde{b_{t_0}^{1/2}}, X_{(n+1)\delta+11t_0} \rangle \Big| \mathcal{F}_t \right)  + \langle  b_{t_0}^{1/2}, \widehat{\phi}_1 \rangle_\mu W_t\nonumber\\
	& =: \mathcal{M}_t^{b} + \langle  b_{t_0}^{1/2}, \widehat{\phi}_1 \rangle_\mu W_t,\qquad t\in[n\delta, (n+1)\delta).
\end{align}
Since $\mathcal{M}_t^b$ is a martingale for $t\in [0, (n+1)\delta+11t_0]$, it follows from Lemma  \ref{lemma4} and the $L^2$-maximal inequality that for any $n\geq 0$,
\begin{align}
	& \mathbb{E}_{\delta_x}\left(\sup_{t\in [n\delta,(n+1)\delta]}  \left(\mathcal{M}_t^b\right)^2\right)  \leq 4 \mathbb{E}_{\delta_x}\left( \left(\mathcal{M}_{(n+1)\delta +11t_0}^b\right)^2\right) \nonumber\\
	& = 4e^{2\lambda_1 ((n+1)\delta+11t_0)} \mathbb{E}_{\delta_x}\left( \langle   \widetilde{b_{t_0}^{1/2}}, X_{(n+1)\delta +11t_0} \rangle ^2  \right)  \lesssim_{t_0} e^{-c(b_{t_0}^{1/2}) (n+1)\delta } (b_{t_0}^{1/2}(x)+b_{t_0}(x)),
\end{align}
which implies that
$\sup_{t\in [n\delta, (n+1)\delta)} \mathcal{M}_t^{b}\to 0$ almost surely as $n\to \infty$.
Plugging this back to \eqref{Z-T} yields that
\begin{align}\label{Step14}
	\limsup_{t\to\infty} \left(e^{\lambda_1 t}Z_t^{(k)}\right)\leq  \langle  b_{t_0}^{1/2}, \widehat{\phi}_1 \rangle_\mu W_\infty ,\ \ \mathbb{P}_{\delta_x}\left(\cdot| \mathcal{E}^c\right)\textup{-a.s.}
\end{align}
For  $k\in \mathbb{I}$, let $e_k^{(j)}$ be the vector with $(e_k^{(j)})_i=\delta_{i,j}$ for $1\leq i\leq n_k$. Taking $g_k=\phi_j^{(k)}= (\Phi_k(x))^T e_k^{(j)}$ in  \eqref{Step13-2} and combining \eqref{Step14} with $\rho\downarrow 0$ first and then letting $\delta\downarrow 0$, the dominated convergence theorem implies  $\mathbb{P}_{\delta_x}(\cdot| \mathcal{E}^c)$-almost surely,
\begin{align}\label{e4}
		&\lim_{\delta\to0} \sup_{1\leq j\leq n_k}  \limsup_{n\to\infty} \sup_{t\in [n\delta, (n+1)\delta)} \frac{e^{\lambda_1 n\delta/2}\left|\langle T_{(n+1)\delta-t}\phi_j^{(k)}, X_{t} \rangle - \langle T_{\delta}\phi_j^{(k)}, X_{n\delta} \rangle \right|}{\sqrt{2\log (n\delta)}} =0.
\end{align}
Now let $\{s_k^{(j)}((n+1)\delta-t), 1\leq j\leq n_k\}$ be a collection of coefficients
such that
\[
e^{\lambda_k ((n+1)\delta-t)}D_k((n+1)\delta-t)^{-1}v =\sum_{j=1}^{n_k} s_k^{(j)}((n+1)\delta-t) e_k^{(j)}.
\]
 Then it is simple to see that $g_k= \sum_{j=1}^{n_k} s_k^{(j)}((n+1)\delta-t) T_{(n+1)\delta-t}\phi_j^{(k)}$. By \eqref{T-t-Phik},
\begin{align}
& |\langle g_k, X_t\rangle -\langle T_{t-n\delta}g_k, X_{n\delta}\rangle| \\
&= |\langle (\Phi_k)^Tv, X_t\rangle -e^{-\lambda_k(t-n\delta  )}\langle (\Phi_k)^T D_k(t-n\delta)v, X_{n\delta}\rangle| \nonumber\\
& = |\sum_{j=1}^{n_k} s_k^{(j)}((n+1)\delta-t) \langle T_{(n+1)\delta -t}\phi_j^{(k)}, X_t\rangle -\sum_{j=1}^{n_k} s_k^{(j)}((n+1)\delta-t)\langle T_\delta \phi_j^{(k)}, X_{n\delta}\rangle|\nonumber\\
&\leq \sum_{j=1}^{n_k} |s_k^{(j)}((n+1)\delta-t)| \cdot | \langle T_{(n+1)\delta -t}\phi_j^{(k)}, X_t\rangle - \langle T_\delta \phi_j^{(k)}, X_{n\delta}\rangle|.
\end{align}
Since $\sup_{r\in (0,1), 1\leq j\leq n_k}|s_k^{(j)}(r)|<\infty$,
combining the display above with \eqref{e4} yields the desired result.

\hfill$\Box$

\begin{lemma}\label{lemma9}
	If $f$ is of the form \eqref{e:rsspecial},
	then for any $r\in (0,\infty)$,
	\begin{align}\label{LIMSUP7}
		\limsup_{t\to\infty}/\liminf_{t\to\infty} \frac{e^{\lambda_1 t/2}
					\sum_{i=1}^{M_{t}} Y_{t}^{f,i}(\infty,r)}{\sqrt{2\log t}} = +/- \sqrt{\langle V_{\infty,r}^f ,\widehat{\phi}_1 \rangle_\mu W_\infty},
		\quad \mathbb{P}_{\delta_x}\left(\cdot| \mathcal{E}^c\right)
				\mbox{-a.s.}
	\end{align}
\end{lemma}	
\textbf{Proof :}
We only prove the $\limsup$ assertion, the proof of the $\liminf$ assertion is  similar.
Let $\delta= r/\ell$ for some $\ell\in \mathbb{N}$.
It follows from \eqref{Decomposition} and definition \eqref{def-Y(f,i)}
that for $n\delta \le t$,
\begin{align}
	& \sum_{i=1}^{M_{t}} Y_{t}^{f,i}(\infty,r)- \sum_{i=1}^{M_{n\delta}} Y_{n\delta}^{T_{t-n\delta}f,i}(\infty,r)\nonumber\\
	& = \left(\langle f, X_{t+r} \rangle - \langle T_{t-n\delta}f, X_{n\delta+r} \rangle\right)+\left(\langle
	T_r
	(f_{sm}+f_{cr}), X_{t} \rangle - \langle T_{t-n\delta+r}(f_{sm}+f_{cr}), X_{n\delta} \rangle\right).
\end{align}
Note that by Lemma \ref{lem:conti-mart},
\begin{align}
	 \lim_{\ell \to\infty} 	\limsup_{n\to\infty} \sup_{t\in [n\delta, (n+1)\delta]} \frac{e^{\lambda_1 n\delta/2}\left| \langle f, X_{t+r} \rangle - \langle T_{t-n\delta}f, X_{n\delta+r} \rangle\right|}{\sqrt{2\log (n\delta)}}=0,\quad  \mathbb{P}_{\delta_x}\left(\cdot| \mathcal{E}^c\right)\textup{-a.s.}
\end{align}
and that $ \mathbb{P}_{\delta_x}\left(\cdot| \mathcal{E}^c\right)$-almost surely,
\begin{align}
	& \lim_{\ell\to\infty} \limsup_{n\to\infty} \sup_{t\in [n\delta, (n+1)\delta]} e^{\lambda_1 n\delta/2} \nonumber
	\frac{\left| \langle T_{r}(f_{sm}+f_{cr}), X_{t} \rangle - \langle T_{t-n\delta+r}(f_{sm}+f_{cr}), X_{n\delta} \rangle\right|}{\sqrt{2\log (n\delta)}}=0.
\end{align}
Therefore,  $\mathbb{P}_{\delta_x}\left(\cdot| \mathcal{E}^c\right)$-almost surely,
\begin{align}\label{conti-small-br1}
	&\lim_{\ell \to \infty} \limsup_{n\to\infty} \sup_{t\in [n\delta, (n+1)\delta]} \frac{e^{\lambda_1 n\delta/2}\left| \sum_{i=1}^{M_{t}} Y_{t}^{f,i}(\infty,r)- \sum_{i=1}^{M_{n\delta}} Y_{n\delta}^{T_{t-n\delta}f,i}(\infty,r) \right|}{\sqrt{2\log (n\delta)}}=0.
\end{align}
In light of Lemma \ref{lemma12} and \eqref{conti-small-br1}, to prove \eqref{LIMSUP7}, it suffices to show that $ \mathbb{P}_{\delta_x}\left(\cdot| \mathcal{E}^c\right)$-almost surely,
\begin{align}\label{conti-small-br2}
	&\lim_{\ell \to \infty} \limsup_{n\to\infty} \sup_{t\in [n\delta, (n+1)\delta]} \frac{e^{\lambda_1 n\delta/2}\left|  \sum_{i=1}^{M_{n\delta}} Y_{n\delta}^{T_{t-n\delta}f,i}(\infty,r)-\sum_{i=1}^{M_{n\delta}} Y_{n\delta}^{f,i}(\infty,r) \right|}{\sqrt{2\log (n\delta)}}=0.
\end{align}
Recall that $e_k^{(j)}$ is a $\mathbb{C}^{n_k}$-valued vector with $(e_k^{(j)})_i = \delta_{i,j}$ for $1\leq i\leq n_k$.
Define
\begin{align}
	T_{t-n\delta} f- f & =  \sum_{k\in \mathbb{I}: k\leq m} \left(\Phi_k(x)\right)^T (e^{-\lambda_k (t-n\delta)} D_k(t-n\delta)- I) v_k\nonumber\\
	&=: \sum_{k\in \mathbb{I}: k\leq m} \sum_{j=1}^{n_k}\widehat{s}_j^{(k)}(t-n\delta) \left(\Phi_k(x)\right)^T e_k^{(j)}.
\end{align}
Then by the linearity of $Y^f(\infty,r)$ with respect to $f$ (see definition \eqref{def-Y(f,i)}), we have
\begin{align}
	&\left|  \sum_{i=1}^{M_{n\delta}} Y_{n\delta}^{T_{t-n\delta}f,i}(\infty,r)-\sum_{i=1}^{M_{n\delta}} Y_{n\delta}^{f,i}(\infty,r) \right|=\left|  \sum_{i=1}^{M_{n\delta}} Y_{n\delta}^{T_{t-n\delta}f-f,i}(\infty,r)\right|\nonumber\\
	&\leq \sum_{k\in \mathbb{I}: k\leq m} \sum_{j=1}^{n_k}\left| \widehat{s}_j^{(k)}(t-n\delta)\right| \left|  \sum_{i=1}^{M_{n\delta}} Y_{n\delta}^{\phi_j^{(k)},i}(\infty,r)\right|\nonumber\\
	& \leq \sup_{
	\tilde{t}
	\in (0,\delta), k\in \mathbb{I}, k\leq m, 1\leq j\leq n_k}\left| \widehat{s}_j^{(k)}(\tilde{t})\right| \sum_{k\in \mathbb{I}: k\leq m} \sum_{j=1}^{n_k} \left|  \sum_{i=1}^{M_{n\delta}} Y_{n\delta}^{\phi_j^{(k)},i}(\infty,r)\right|.
\end{align}
Applying Lemma \ref{lemma12} to $\sum_{i=1}^{M_{n\delta}} Y_{n\delta}^{\phi_j^{(k)},i}(\infty,r)$ for $k=1, \dots, m$ and $j=1, \dots, n_k$ in the inequality above, we see that, to prove \eqref{conti-small-br2}, it suffices to show that for
$k=1, \dots, m$ and $1\leq j\leq n_k$,
\begin{align}\label{conti-small-br3}
	\lim_{\delta\to 0} \sup_{\tilde{t}\in (0,\delta)}\left| \widehat{s}_j^{(k)}(\tilde{t})\right| =0.
\end{align}
Since $\widehat{s}_j^{(k)}(\tilde{t})$ is a polynomial of $\tilde{t}$ with $\widehat{s}_j^{(k)}(0)=0$, \eqref{conti-small-br3} holds trivally.
Hence \eqref{conti-small-br2} is valid. The proof is now complete.

\hfill$\Box$

As a consequence of Lemma \ref{lemma9}, we have the following useful collory:

\begin{cor}\label{Cor}
	If $k\in \mathbb{I}$ with $2\mathfrak{R}_k>\lambda_1$, then for each $1\leq j\leq n_k$, it holds that
	\begin{align}
		D_+(j,k):= \limsup_{t\to\infty} 	\frac{ e^{\frac{\lambda_1}{2} t}   \left| \langle \phi_j^{(k)}, X_{t} \rangle\right| }{\sqrt{\log t}} < \infty,  \quad \mathbb{P}_{\delta_x}\left(\cdot| \mathcal{E}^c\right)
				\mbox{-a.s.}
	\end{align}
\end{cor}
\textbf{Proof: } Fix $k\in \mathbb{I}$.
In light of Lemma \ref{lem:conti-mart}, \eqref{conti-small-br1} and \eqref{conti-small-br2}, to prove the desired assertion, it suffices to show that for any small $\delta>0$,
\begin{align}\label{step_43}
	\limsup_{n\to\infty} 	\frac{ e^{\frac{\lambda_1}{2} n\delta}   \left| \langle \phi_j^{(k)}, X_{n\delta} \rangle\right| }{\sqrt{\log (n\delta)}} < \infty, \quad \mathbb{P}_{\delta_x}\left(\cdot| \mathcal{E}^c\right)\textup{-a.s.}
\end{align}
 We prove \eqref{step_43} by induction. When $j=1$, then $T_t \phi_1^{(k)}= e^{-\lambda_k t}\phi_1^{(k)}$. Fix an arbitrary $L\in \N$. By
\eqref{e2} (with $\delta$ replaced by $\delta/L$ and $g_k=\phi_1^{(k)}$), we have
\begin{align}\label{Induction-1}
	\limsup_{n\to\infty} \frac{e^{\lambda_1 n\delta/(2L)}\left|\langle \phi_1^{(k)}, X_{(n+1)\delta/L} \rangle - e^{-\lambda_k \delta/L}\langle \phi_1^{(k)}, X_{n\delta/L} \rangle \right|}{\sqrt{2\log (n\delta/L)}} <\infty,\quad \mathbb{P}_{\delta_x}(\cdot|\mathcal{E}^c)-\textup{a.s.}
\end{align}
Therefore,
 there exists a finite random variable
 $\mathcal{U}=\mathcal{U}(k,\delta, L)$
 such that almost
 surely, for $n$ large enough,
\begin{align}\label{step_41}
	& \left| e^{\lambda_k n\delta/L} \langle \phi_1^{(k)}, X_{n\delta/L} \rangle -e^{\lambda_k n\delta} \langle \phi_1^{(k)}, X_{n\delta} \rangle \right| \nonumber\\
	&\leq  \sum_{q= n}^{nL-1} 	\left| e^{\lambda_k q\delta/L} \langle \phi_1^{(k)}, X_{q\delta/L} \rangle -e^{\lambda_k (q+1)\delta/L} \langle \phi_1^{(k)}, X_{(q+1)\delta/L} \rangle \right|
	\\
&  \leq \mathcal{U} \sum_{q= n}^{nL-1}  \sqrt{2\log (q\delta/L)} e^{\mathfrak{R}_k(q+1)\delta/L -\lambda_1 q\delta/ (2L)}\nonumber\\
&\leq
\mathcal{U}
e^{\mathfrak{R}_k\delta/L}  \sqrt{2\log (n\delta)}\sum_{q= n}^{nL-1}  e^{(-\lambda_1/2+ \mathfrak{R}_k)q\delta/L } .
\end{align}
Note that the right hand side of the above inequality is bounded by
$
\mathcal{U}'
\sqrt{\log(n\delta)} e^{(-\lambda_1/2+ \mathfrak{R}_k)n\delta} $ for some random variable
$\mathcal{U}'$.
Therefore, we have almost surely,
\begin{align}\label{Step15}
 &\limsup_{n\to\infty} \frac{ e^{\lambda_1 n\delta/2}  \left| e^{-\lambda_k n\delta(L-1)/L} \langle \phi_1^{(k)}, X_{n\delta/L} \rangle - \langle \phi_1^{(k)}, X_{n\delta} \rangle \right| }{\sqrt{\log (n\delta )}}\nonumber\\
  &=
 \limsup_{n\to\infty} \frac{ e^{-\frac{2\mathfrak{R}_k -\lambda_1}{2} n\delta}  \left| e^{\lambda_k n\delta/L} \langle \phi_1^{(k)}, X_{n\delta/L} \rangle -e^{\lambda_k n\delta} \langle \phi_1^{(k)}, X_{n\delta} \rangle \right| }{\sqrt{\log (n\delta )}}  <\infty .
\end{align}
Combining Lemma \ref{lemma10} (with $f= |\phi_1^{(k)}|$ and $\delta$ replaced by $\delta/L$) and the assumption $2\mathfrak{R}_k>\lambda_1$,  taking $L> \frac{2(\mathfrak{R}_k-\lambda_1)}{2\mathfrak{R}_k-\lambda_1}$,
we get that
\begin{align}\label{Step16}
	& e^{\frac{\lambda_1}{2} n\delta} \left|  e^{-\lambda_k n\delta(L-1)/L}\langle \phi_1^{(k)}, X_{n\delta/L} \rangle \right| \leq e^{\frac{\lambda_1}{2} n\delta}   e^{-\mathfrak{R}_k n\delta(L-1)/L}\langle \left|\phi_1^{(k)}\right| , X_{n\delta/L} \rangle \nonumber\\ &
		\lesssim_{k}
	 e^{\frac{\lambda_1}{2} n\delta}   e^{-\mathfrak{R}_k n\delta(L-1)/L} e^{-\lambda_1 n\delta/L} = e^{-\frac{n\delta}{2L}\left((2\mathfrak{R}_k-\lambda_1)L-2(\mathfrak{R}_k-\lambda_1)\right)  } \stackrel{n\to\infty}{\longrightarrow}0.
\end{align}
Combining Lemma \ref{lem:conti-mart},  \eqref{Step15} and \eqref{Step16}, we get \eqref{step_43} for $j=1$.

Suppose that \eqref{step_43} holds for all
 $\ell =1,...,j-1$. It suffices to show \eqref{Induction-1} holds with $\phi_1^{(k)}$ replaced by $\phi_j^{(k)}$.
We will use \eqref{e2} for $g_k=\phi_j^{(k)}$. Note that
\[
T_\delta \phi_j^{(k)}= e^{-\lambda_k \delta} (\Phi_k(x))^TD_k(\delta)e_k^{(j)} =  e^{-\lambda_k \delta} \sum_{q=1}^{n_k} (D_k(\delta))_{q,j} (\Phi_k(x))^T e_{k}^{(q)} = e^{-\lambda_k \delta} \sum_{q=1}^{j} (D_k(\delta))_{q,j} \phi_q^{(k)}
,\]
 where in the last equality we used the fact that $(D_k(\delta))_{q,j}=0$ when $q>j$.  Therefore, it follows from \eqref{e2} and
the induction hypothesis that
\begin{align*}
		&\limsup_{n\to\infty} \frac{e^{\lambda_1 n\delta/2}\left|\langle \phi_j^{(k)}, X_{(n+1)\delta} \rangle - e^{-\lambda_k \delta}\langle \phi_j^{(k)}, X_{n\delta} \rangle \right|}{\sqrt{2\log (n\delta)}} \nonumber\\
		&\leq \limsup_{n\to\infty} \frac{e^{\lambda_1 n\delta/2}\left|\langle \phi_j^{(k)}, X_{(n+1)\delta} \rangle - \langle T_{\delta}\phi_j^{(k)}, X_{n\delta} \rangle \right|}{\sqrt{2\log (n\delta)}} \nonumber\\
		& \quad + e^{-\mathfrak{R}_k \delta}\sum_{q=1}^{j-1}\left|(D_k(\delta))_{q,j}\right|\limsup_{n\to\infty} 	\frac{ e^{\frac{\lambda_1}{2} n\delta}   \left| \langle \phi_j^{(k)}, X_{n\delta} \rangle\right| }{\sqrt{\log (n\delta)}}<\infty.
\end{align*}
Thus \eqref{Induction-1} is valid. The proof is now complete.

\hfill$\Box$

\subsection{The case of test functions with no critical components}\label{ss: 4.2}
\textbf{Proof of Theorem \ref{thm2}:}
We only prove the $\limsup$ assertion, the proof of the $\liminf$ assertion is  similar.
Recall the definitions of $E_t$ and $Y_t^{f,i}(\infty, r)$ in \eqref{def-of-e-t} and \eqref{def-Y(f,i)}, respectively. By Lemma \ref{lemma9},
\begin{align}\label{e5}
	&	\limsup_{t\to\infty}\! \frac{e^{\lambda_1 t/2}\!\left(\! \langle f, X_{t+r} \rangle - \langle T_r f_{sm}, X_{t}\rangle -E_{t+r}(f_{la})\right)}{\sqrt{2\log t}} = \sqrt{\!\langle V_{\infty,r}^f ,\widehat{\phi}_1 \rangle_\mu W_\infty},\quad \mathbb{P}_{\delta_x}\!\left(\cdot| \mathcal{E}^c\right)\textup{-a.s.}
\end{align}
Recall the definition of $D_+(j,k)$ in Corollary \ref{Cor},
and note that
 $D_+(j,k)<\infty$ almost surely for all $ 2\mathfrak{R}_k>\lambda_1$ and $1\leq j\leq n_k$.
We write
\[
T_rf_{sm}=
\sum_{k\leq m: 2\mathfrak{R}_k>\lambda_1} e^{-\lambda_k r}\left(\Phi_k(x)\right)^T D_k(r)v_k
\]
in the form
\[
\sum_{k\leq m: 2\mathfrak{R}_k>\lambda_1}\sum_{j=1}^{n_k} e^{-\lambda_k r}R_j^{(k)}(r)  \phi_j^{(k)}.
\]
Then each $R_j^{(k)}(r)$ is a polynomial of $r$ of degree at most $n_k$.
Therefore, there exists some
 constant $\Gamma$ depending on $v_1,...,v_m$ such that almost surely,
\[
\limsup_{t\to\infty} \frac{e^{\lambda_1 (t+r)/2}\left|\langle T_r f_{sm}, X_{t}\rangle \right|}{\sqrt{2\log t}} \leq \Gamma \sum_{k\leq m: 2\mathfrak{R}_k>\lambda_1}e^{(\lambda_1/2 -\mathfrak{R}_k) r}(1+r)^{n_k} \sum_{j=1}^{n_k} D_+(j,k).
\]
Multiplying both sides of \eqref{e5} by $e^{\lambda_1 r/2}$ and applying  the inequality
\begin{align}\label{Limsup-ineq}
	\limsup_{t\to\infty} x_t  + \liminf_{t\to\infty} y_t \leq \limsup_{t\to\infty} (x_t +y_t) \leq \limsup_{t\to\infty} x_t + \limsup_{t\to\infty} y_t,
\end{align}
we get that  for any $r>0$,
\begin{align*}
	& \limsup_{t\to\infty} \frac{e^{\lambda_1 t/2}\left( \langle f, X_{t} \rangle  -E_{t}(f_{la})\right)}{\sqrt{2\log t}} - \ \Gamma \sum_{k\leq m: 2\mathfrak{R}_k>\lambda_1}e^{(\lambda_1/2 -\mathfrak{R}_k) r}(1+r)^{n_k} \sum_{j=1}^{n_k} D_+(j,k)\nonumber\\
	&
	\leq \sqrt{\langle e^{\lambda_1 r}V_{\infty,r}^f ,\widehat{\phi}_1 \rangle W_\infty} \nonumber\\
	&\leq \limsup_{t\to\infty} \frac{e^{\lambda_1 t/2}\left( \langle f, X_{t} \rangle  -E_{t}(f_{la})\right)}{\sqrt{2\log t}} + \Gamma \sum_{k\leq m: 2\mathfrak{R}_k>\lambda_1}e^{(\lambda_1/2 -\mathfrak{R}_k) r}(1+r)^{n_k} \sum_{j=1}^{n_k} D_+(j,k).
\end{align*}
Letting $r\to\infty$ in the display above
yields that $\mathbb{P}_{\delta_x}(\cdot|\mathcal{E}^c)$-almost surely,
\begin{align}
	\limsup_{t\to\infty} \frac{e^{\lambda_1 t/2}\left( \langle f, X_{t} \rangle  -E_{t}(f_{la})\right)}{\sqrt{2\log t}} = \lim_{r\to\infty} \sqrt{\langle e^{\lambda_1 r}V_{\infty,r}^f ,\widehat{\phi}_1 \rangle_\mu W_\infty} .
\end{align}
Therefore, to get the desired result,
it suffices to show that
\begin{align}\label{Identity-of-lim}
	\lim_{r\to\infty} \langle e^{\lambda_1 r}V_{\infty,r}^f ,\widehat{\phi}_1 \rangle_\mu=	\sigma_{sm}^2(f)+ \sigma_{la}^2(f) .
\end{align}
Define $Q:=\langle f_{la}, X_r\rangle -\sum_{2\mathfrak{R}_k < \lambda_1} e^{-\lambda_k r}H_\infty^{(k)}D_k(r)v_k$, then $\mathbb{E}_{\delta_x}(Q\big| \mathcal{F}_r)=0$.
Therefore,
\begin{align}
	V_{\infty, r}^f(x) & = \mathbb{E}_{\delta_x}\left(  \mathbb{E}_{\delta_x} \left( (\langle f_{sm} , X_r\rangle	-T_rf_{sm}(x)	+Q )^2 \big|\mathcal{F}_r\right)\right)= \mbox{Var}_x\left(\langle f_{sm} , X_r\rangle\right)+ \mbox{Var}_x\left(Q\right).
\end{align}
Noticing that $e^{\lambda_1 r} \mbox{Var}_x\left(\langle f_{sm} , X_r\rangle\right) \to
\sigma_{sm}^2(f)
\phi_1(x)$ and that $e^{\lambda_1 r} \mbox{Var}_x\left(\langle f_{sm} , X_r\rangle\right) \lesssim_{f} b_{t_0}^{1/2}(x) + b_{t_0}(x)$ for all $t>10t_0$,
 applying the dominated convergence theorem, we get that
\begin{align}\label{Identity-of-lim1}
		\lim_{r\to\infty} \langle e^{\lambda_1 r}\mbox{Var}_{\cdot}\left(\langle f_{sm} , X_r\rangle\right)  ,\widehat{\phi}_1 \rangle_\mu=	\sigma_{sm}^2(f).
\end{align}
For $Q$,
by the branching property, we have
\begin{align}
	\mbox{Var}_x\left(Q\right)& = \mathbb{E}_{\delta_x}\left( \mbox{Var}_x\left(Q\Big|
	\mathcal{F}_r
	\right) \right) =\mathbb{E}_{\delta_x}\bigg(\langle \mbox{Var}_{\cdot} \Big(\sum_{k\leq m: 2\mathfrak{R}_k<\lambda_1} H_\infty^{(k)}v_k\Big), X_r\rangle\bigg)\nonumber\\
	& = T_r\bigg(\mbox{Var}_{\cdot} \Big(\sum_{k\leq m: 2\mathfrak{R}_k<\lambda_1} H_\infty^{(k)}v_k\Big)\bigg)(x).
\end{align}
Therefore,
combining Lemma \ref{lemma3} (1) and the dominated convergence theorem,  we get that
\begin{align}\label{Identity-of-lim2}
		\lim_{r\to\infty} \langle e^{\lambda_1 r}\mbox{Var}_{\cdot}\left(Q\right)  ,\widehat{\phi}_1 \rangle_\mu= \langle \mbox{Var}_{\cdot} \Big(\sum_{k\leq m: 2\mathfrak{R}_k<\lambda_1} H_\infty^{(k)}v_k\Big),\widehat{\phi}_1\rangle_\mu=\sigma_{la}^2(f),
\end{align}
where the last equality follows from \cite[(3.48)]{RSZ2017}. Combining \eqref{Identity-of-lim1} and \eqref{Identity-of-lim2},
 we get \eqref{Identity-of-lim}. The proof is complete.

\hfill$\Box$

\subsection{The case of test functions with non-trivial critical components}\label{ss:4.3}

The main goal in this subsection is to prove the following theorem.
\begin{thrm}\label{thm3}
	 If $f\in L^2(E, \mu)$  and $f_{cr}\neq 0$, then
	      $\mathbb{P}_{\delta_x}\left(\cdot \vert \mathcal{E}^c\right)$-almost surely,
	\begin{align*}
		\limsup_{t\to\infty}/\liminf_{t\to\infty} \frac{e^{\lambda_1 t /2}\langle f_{cr} , X_{t}\rangle}{\sqrt{2 t^{1+2\tau(f_{cr})} \log \log t}} = +/- \sqrt{\sigma_{cr}^2(f) W_\infty}.
	\end{align*}
\end{thrm}

We first give the proof of Theorem \ref{thm4} using Theorem \ref{thm3}.

\textbf{Proof of Theorem \ref{thm4}:}
Applying Theorem \ref{thm2} to $f-f_{cr}$  and  Theorem \ref{thm3} to $f_{cr}$, we immediately get the desired result.

\hfill$\Box$

Suppose $f_{cr}(x)= \sum_{k: 2\mathfrak{R}_k = \lambda_1}
\left(\Phi_k(x)\right)^T v_k$ with $v_k\in \mathbb{C}^{n_k}$ and $ \overline{v}_k=v_{k'}$.
We now rewrite $f_{cr}$ in a different form. In this parapgraph, we always assume $k\in \mathbb{I}$ satisfies $2\mathfrak{R}_k = \lambda_1$.
Recall that $e_k^{(j)}$ is a $\mathbb{C}^{n_k}$-valued vector with $(e_k^{(j)})_i = \delta_{i,j}$ for $1\leq i\leq n_k$ and $\phi_j^{(k)}=(\Phi_k(x))^T e_k^{(j)}$. For each $k$, define $\nu_{k,0}:=0, \nu_{k,i}:=\sum_{m=1}^i d_{k,m}, 1\le i\leq r_k$ and $ d_k:= \max_{1\leq i\leq r_k}d_{k,i}$.
For $1\le i\leq r_k$ and $1\le j\le d_{k, i}$, let $\theta_{k,i}^{(j)}$ be the coefficient of $\phi_{\nu_{k, i-1}+j}^{(k)}$ in $f_{cr}$. Note $\theta_{k,i}^{(\ell)}=\overline{\theta}_{k',i}^{(\ell)}$.
Then $f_{cr}(x)$ can be rewritten as
\[
f_{cr}(x)=
\sum_{k: 2\mathfrak{R}_k=\lambda_1} \sum_{i=1}^{r_k} \sum_{\ell=1}^{d_{k, i}}\theta_{k,i}^{(\ell)} (\Phi_k(x))^T  e_k^{(\nu_{k,i-1}+\ell)}.
\]
Let $\Phi_k= (\Phi_{k,i}, 1\leq i\leq r_k)$ where $\Phi_{k,i}$ is a $\mathbb{C}^{d_{k,i}}$-valued function for each $i=1, \dots, r_k$.
 For each $\ell \leq d_{k,i}$, let
  $e_{k,i}^{(\ell)}$ be the $\mathbb{C}^{d_{k,i}}$-vector with $(e_{k,i}^{(\ell)})_q=\delta_{\ell, q}$.
 Let $d:=\max_{k: 2\mathfrak{R}_k=\lambda_1} d_k$.
For $1\le \ell\le d$,
set
$\mathcal{A}_{\ell}
= \left\{(k, i):\  2\mathfrak{R}_k=\lambda_1, 1\leq i\leq r_k, \ell\leq d_{k,i} \right\}$ and
\begin{align}\label{Def-of-Q}
Q_\ell (x): =
\sum_{(k,i)\in \mathcal{A}_{\ell}}
 \theta_{k,i}^{(\ell)}  (\Phi_{k,i}(x))^T  e_{k,i}^{(\ell)}.
\end{align}
Then
\[
 f_{cr}(x)=\sum_{\ell=1}^{d}Q_\ell(x).
\]
It is easy to see that if $Q_\ell \neq 0$, then $\tau(Q_\ell)=\ell-1$.
For any $t>0$,
\[
T_tQ_\ell (x)=
\sum_{(k,i)\in \mathcal{A}_{\ell} } e^{-\lambda_k t}
\theta_{k,i}^{(\ell)}
(\Phi_{k,i}(x))^T J_{k,i}(t)e_{k,i}^{(\ell)}.
\]
We first consider  $Q_1$.
Set $\theta_{k,i}:= \theta_{k,i}^{(1)}$ and  $\phi_1^{(k,i)}:= \phi_{\nu_{k,i-1}+1}^{(k)}$. For $t>0$, define
$$
T_{-t}Q_1 (x):= \sum_{(k,i)\in\mathcal{A}_{1}} e^{\lambda_k t}
\theta_{k,i}
(\Phi_{k,i}(x))^T J_{k,i}(-t)
e_{k,i}^{(1)} = \sum_{(k,i)\in\mathcal{A}_{1}} e^{\lambda_k t}
\theta_{k,i}
\phi_1^{(k,i)}(x)
$$ and
\[
\mathcal{W}_t= \mathcal{W}_t^{(1)}:=
\sum_{(k,i)\in \mathcal{A}_{1}} e^{\lambda_k t} \theta_{k,i} \langle \phi_1^{(k,i)}, X_t\rangle
= \langle T_{-t}Q_1, X_t\rangle.
\]
Then it is easy to see that $\mathcal{W}_t$ is a martingale.  Let $B_t$ be an independent standard Brownian motion.
We say a sequence $\{a_k: k=0,1,...\}$ of integers is
{\it syndetic}
 if $a_0<a_1<\cdots $ and $\sup_{k\in \mathbb{N}} (a_{k+1}-a_k)<\infty$.
Suppose that $\{a_k: k=0,1,...\}$ is a syndetic sequence such that $a_0=0$ and $a_{k+1}-a_k\in [1, N]$ for any $k\in \mathbb{N}$,
where $N$ is a positive integer.
Then for any $\delta, \varepsilon>0$,
\begin{align}\label{Def-of-Z-n}
  Z_n:= \mathcal{W}_{a_n \delta} - \mathcal{W}_{0} +\varepsilon B_{a_n} = \sum^n_{j=1}(\mathcal{W}_{a_j\delta}-\mathcal{W}_{a_{j-1}\delta})+\varepsilon B_{a_n} ,\quad n\in \mathbb{N},
\end{align}
 is a martingale.  For simplicity, define $\mathcal{G}_j^Z:= \mathcal{F}_{a_{j}\delta}\vee \sigma( B_r, r\leq a_{j}).$

 By the branching property, for each $j\in \mathbb{N}$,
 \[
 \mathcal{W}_{a_j \delta} -\mathcal{W}_{a_{j-1} \delta} = \sum_{i=1}^{M_{a_{j-1}\delta}}\left(\langle T_{-a_j\delta}Q_1, X_{a_j\delta-a_{j-1}\delta}^i\rangle -T_{-a_{j-1}\delta}Q_1(X_{a_{j-1}\delta}(i))\right)=: \sum_{i=1}^{M_{a_{j-1}\delta}} \mathcal{Y}_j^{Q_1,i}.
 \]
 Define $\mathcal{Y}_j^{Q_1}:= \langle T_{-a_j\delta}Q_1, X_{a_j\delta-a_{j-1}\delta}\rangle -\langle T_{-a_{j-1}\delta}Q_1, X_0\rangle$.  Note that  there exists some constant $C(Q_1)$ such that
 \begin{align}\label{bound-f-cr}
 	|\mathcal{Y}_j^{Q_1}| &  \leq C(Q_1)e^{\lambda_1 a_j\delta/2} \sum_{(k,i)\in\mathcal{A}_{1}}    \left|\langle \phi_1^{(k,i)}, X_{a_j\delta-a_{j-1}\delta}\rangle\right|  \nonumber\\
 	&\qquad +  C(Q_1)e^{\lambda_1 a_{j-1}\delta/2} \sum_{(k,i)\in\mathcal{A}_{1}}    \left|\langle \phi_1^{(k,i)}, X_{0}\rangle\right| .
 \end{align}
 Since $a_{j}-a_{j-1} \in \{1,...,N\}$,  there exists a constant $C'= C'(Q_1,\delta)$ such that under $\mathbb{P}_{\delta_x}$,
 \begin{align}\label{Up-of-Y}
 	|\mathcal{Y}_j^{Q_1}| \leq C' e^{\lambda_1 a_{j-1}\delta/2}  \sum_{(k,i)\in\mathcal{A}_{1}}    \sum_{q=0}^N \left| \langle   \phi_1^{(k,i)}, X_{q\delta} \rangle \right|=: e^{\lambda_1 a_{j-1}\delta/2} \Upsilon.
 \end{align}
 Define
 \[
 \mathcal{U}_{a_j \delta}^{Q_1}:= \sum_{i=1}^{M_{a_{j-1}\delta}} \left(\mathcal{Y}_j^{Q_1,i} 1_{\{|\mathcal{Y}_j^{Q_1,i} | \leq 1 \} }-  \mathbb{E}_{\delta_x}\left[\mathcal{Y}_j^{Q_1,i} 1_{\{|\mathcal{Y}_j^{Q_1,i} | \leq 1 \}} |\mathcal{F}_{a_{j-1}\delta}\right]\right)
 \]
 and
 \begin{align}\label{def-of-modified-mar}
 	\mathcal{Z}_n:= \sum_{j=1}^n \mathcal{U}_{a_j \delta}^{Q_1} + \varepsilon
 	 	B_{a_n}
 	=: \sum_{j=1}^n \mathcal{X}_j.
 \end{align}
 We are going to use \cite[Theorem 4.7, p.117]{HH1980} to prove a discrete-time  law of iterated logarithm  in Lemma \ref{LIL-W-cr}.
 Before that, we first prove a limit theorem for  $\mathcal{Z}_j$, which is used to check  the conditions of \cite[Theorem 4.7, p.117]{HH1980}.

By \cite[(2.20)]{RSZ2017},  
for any $f\in L^2(E,\mu;\mathbb{C})\cap L^4(E,\mu;\mathbb{C})$, $t>0$ and $x\in E$,
\begin{equation}\label{Second-moment}
	\mathbb{E}_{\delta_x}\left( \left|\langle f, X_t \rangle\right|^2 \right)= \int_0^t T_s\left[
	A^{(2)}
	|T_{t-s}f|^2 \right](x)\mathrm{d}s + T_t(|f|^2)(x).
\end{equation}

 \begin{lemma}\label{cor2}
	Let $\delta,\varepsilon>0$ be fixed and $\mathcal{X}_j$ be defined as in \eqref{def-of-modified-mar}.
	Define $s_0^{Q_1}=0$ and
	\begin{align}
		\left(s_n^{Q_1}\right)^2& := \sum_{j=1}^n \mathbb{E}_{\delta_x}\left( \mathcal{X}_j^2\big| \mathcal{G}_{j-1}^Z\right)\in \mathcal{G}_{n-1}^Z.
	\end{align}
	$(1)$   $\mathbb{P}_{\delta_x}$-almost surely,
		\begin{align}
		 \lim_{n\to\infty} \frac{1}{a_n} \sum_{j=1}^n (s_n^{Q_1})^2  &= \varepsilon^2+ \lim_{n\to\infty}  \frac{1}{a_n} \sum_{j=1}^n \langle \mbox{Var}_{\cdot}(\langle T_{-a_j\delta}Q_1, X_{a_j\delta-a_{j-1}\delta}\rangle), X_{a_{j-1}\delta}\rangle  \nonumber\\
		 & = \varepsilon^2 + \delta W_\infty \sigma_{cr}^2(Q_1).
	\end{align}
	 Consequently,  $\mathbb{P}_{\delta_x}$-almost surely, $s_n^{Q_1} \to \infty$ and $s_n^{Q_1}/ s_{n+1}^{Q_1}\to 1$.
	\newline
	$(2)$ $\mathbb{P}_{\delta_x}$-almost surely,
	\begin{align}\label{Fini-4-mom}
	\sup_{j\in \mathbb{Z}_+} \mathbb{E}_{\delta_x}\left( \mathcal{X}_j^4\big| \mathcal{G}_{j-1}^Z\right)<\infty \quad \mbox{and} 	\quad \sup_{j\in \mathbb{Z}_+} \mathbb{E}_{\delta_x}\left( \mathcal{X}_j^4\right)<\infty.
	\end{align}
\end{lemma}
\textbf{Proof: } (1)
For simplicity, we denote $d_j:= a_j\delta$.
The first equality follows from
the branching property and the independence of branching Markov process $X$ and the Brownian motion $B$.  Now we prove the second equality. According to \eqref{Second-moment},
\begin{align}
	& \mbox{Var}_{x}(\langle T_{-d_j}Q_1, X_{d_j-d_{j-1}}\rangle) \nonumber\\
	& = \int_0^{d_j-d_{j-1}} T_s\left[
	A^{(2)}
	|T_{-d_{j-1}-s}Q_1|^2 \right](x)\mathrm{d}s + T_{d_j-d_{j-1}}(|T_{-d_j}Q_1|^2)(x) - \left|T_{-d_{j-1}}Q_1\right|^2(x).
\end{align}
Noticing that for any $t>0$,
\begin{align}
	|T_{-t}Q_1|^2 &= e^{\lambda_1 t}
	\sum_{(k,i), (q,p)\in \mathcal{A}_{1} } e^{ \mathrm{i}(\mathfrak{I}_k+\mathfrak{I}_q) t} \theta_{k,i} \theta_{q,p} \phi_1^{(k,i)}\phi_1^{(q,p)},
\end{align}
we obtain that
\begin{align}\label{decom-var}
	& \langle \mbox{Var}_{\cdot}(\langle T_{-d_j}Q_1, X_{d_j-d_{j-1}}\rangle), X_{d_{j-1}}\rangle  \nonumber\\
	& =\sum_{(k,i), (q,p)\in \mathcal{A}_{1} } \theta_{k,i} \theta_{q,p} e^{\mathrm{i} (\mathfrak{I}_k+\mathfrak{I}_q) d_{j-1}}   e^{\lambda_1 d_{j-1}}  \langle \int_0^{d_j-d_{j-1}} e^{\lambda_1 s}e^{\mathrm{i} (\mathfrak{I}_k+\mathfrak{I}_q) s} T_s\left[	A^{(2)}	\phi_1^{(k,i)}\phi_1^{(q,p)} \right] \mathrm{d}s, X_{d_{j-1}}\rangle \nonumber\\	
& \quad +\sum_{(k,i), (q,p)\in \mathcal{A}_{1} } \theta_{k,i}\theta_{q,p} e^{ \mathrm{i}(\mathfrak{I}_k+\mathfrak{I}_q) d_{j}} e^{\lambda_1 (d_j-d_{j-1})}  e^{\lambda_1 d_{j-1}} \langle T_{d_j-d_{j-1}}( \phi_1^{(k,i)}\phi_1^{(q,p)}), X_{d_{j-1}}\rangle \nonumber\\	
&\quad - \sum_{(k,i), (q,p)\in \mathcal{A}_{1} }  \theta_{k,i}\theta_{q,p} e^{ \mathrm{i} (\mathfrak{I}_k+\mathfrak{I}_q) d_{j-1}}   e^{\lambda_1 d_{j-1}} \langle  \phi_1^{(k,i)}\phi_1^{(q,p)}, X_{d_{j-1}}\rangle.
\end{align}
We would like to replace $e^{\lambda_1 d_{j-1}}\langle f, X_{d_{j-1}}\rangle$ by $W_{\infty}\langle f, \widehat{\phi}_1\rangle_\mu$, so we set
\begin{align*}
	G_1(j)& := W_{\infty}\sum_{(k,i), (q,p)\in \mathcal{A}_{1} }   \theta_{k,i}\theta_{q,p} e^{\mathrm{i} (\mathfrak{I}_k+\mathfrak{I}_q) d_{j-1}} \nonumber\\	&\qquad \times   \langle  \int_0^{d_j-d_{j-1}} e^{\lambda_1 s}e^{\mathrm{i} (\mathfrak{I}_k+\mathfrak{I}_q) s} T_s\left[	A^{(2)}	\phi_1^{(k,i)}\phi_1^{(q,p)}\right], \widehat{\phi}_1\rangle_\mu  \mathrm{d}s \nonumber\\	& \quad + W_{\infty} \sum_{(k,i), (q,p)\in \mathcal{A}_{1} }  \theta_{k,i}\theta_{q,p} e^{ \mathrm{i}(\mathfrak{I}_k+\mathfrak{I}_q) d_{j}} e^{\lambda_1 (d_j-d_{j-1})}  \langle T_{d_j-d_{j-1}}( \phi_1^{(k,i)}\phi_1^{(q,p)}), \widehat{\phi}_1\rangle_\mu  \nonumber\\	&\quad - W_{\infty} \sum_{(k,i), (q,p)\in \mathcal{A}_{1} }  \theta_{k,i}\theta_{q,p} e^{ \mathrm{i} (\mathfrak{I}_k+\mathfrak{I}_q) d_{j-1}}  \langle  \phi_1^{(k,i)}\phi_1^{(q,p)}, \widehat{\phi}_1\rangle_\mu.
\end{align*}
Then by Lemma \ref{lemma10}, comparing
the terms in  \eqref{decom-var}  with the corresponding terms in $G_1(j)$
(also noticing that  $d_j-d_{j-1}$ only takes finite values), we have  $\mathbb{P}_{\delta_x}$-almost surely,
\begin{align}\label{e6}
	& \lim_{j\to\infty}  \left|  \langle \mbox{Var}_{\cdot}(\langle T_{-d_j}Q_1, X_{d_j-d_{j-1}}\rangle), X_{d_{j-1}}\rangle  - G_1(j)\right| =0 .
\end{align}
Therefore,  $\mathbb{P}_{\delta_x}$-almost surely,
\begin{align}\label{e7}
	& \lim_{n\to\infty} \frac{1}{a_n} \left| \sum_{j=1}^n \langle \mbox{Var}_{\cdot}(\langle T_{-d_j}Q_1, X_{d_j-d_{j-1}}\rangle), X_{d_{j-1}}\rangle  - \sum_{k=1}^n G_1(j)\right| =0 .
\end{align}
Since $\langle T_t(h), \widehat{\phi}_1\rangle_\mu= e^{-\lambda_1t}\langle h,\widehat{\phi}\rangle_\mu$, we see that
\begin{align}
	G_1 (j)&= W_{\infty} \sum_{(k,i), (q,p)\in \mathcal{A}_{1} } \theta_{k,i}\theta_{q,p}	\langle A^{(2)}	\phi_1^{(k,i)}\phi_1^{(q,p)} , \widehat{\phi}_1\rangle_\mu \int_{d_{j-1}}^{d_j} e^{\mathrm{i} (\mathfrak{I}_k+\mathfrak{I}_q) s} \mathrm{d}s	\nonumber\\	& \quad + W_{\infty} \sum_{(k,i), (q,p)\in \mathcal{A}_{} } \theta_{k,i}\theta_{q,p } \langle \phi_1^{(k,i)}\phi_1^{(q,p)}, \widehat{\phi}_1\rangle_\mu \left(e^{ \mathrm{i}(\mathfrak{I}_k+\mathfrak{I}_q) d_{j}}  -e^{ \mathrm{i}(\mathfrak{I}_k+\mathfrak{I}_q) d_{j-1}}  \right)  .
\end{align}
Note that $\int_{d_{j-1}}^{d_j} e^{\mathrm{i} us}\mathrm{d}s= \frac{1}{\mathrm{i} u}(e^{\mathrm{i} ud_j}- e^{\mathrm{i} ud_{j-1}})$ for $u\neq 0$ and that
$$
\Big| \sum_{j=1}^n \left(e^{ \mathrm{i} u d_{j}}  -e^{ \mathrm{i}ud_{j-1}}  \right)\Big| \leq 2, \quad \forall u\in \mathbb{R}.
$$
Thus the main contribution to $\sum_{k=1}^n G_1(j)$ comes from pairs
$((k,i),(q,p))$  with $q=k'$, which together with \eqref{e7} implies that
\begin{align}\label{lim-of-var}
	& \lim_{n\to\infty}  \frac{1}{a_n} \sum_{j=1}^n \langle \mbox{Var}_{\cdot}(\langle T_{-d_j}Q_1, X_{d_j-d_{j-1}}\rangle), X_{d_{j-1}}\rangle =  \lim_{n\to\infty}  \frac{1}{a_n} \sum_{j=1}^n G_(j) \nonumber\\
	&= W_{\infty} \sum_{(k,i), (k',p)\in \mathcal{A}_{1} } \theta_{k,i} \theta_{k',p} 	\langle A^{(2)}	\phi_1^{(k,i)}\phi_1^{(k',p)} , \widehat{\phi}_1\rangle_\mu  \lim_{n\to\infty}  \frac{1}{a_n} \sum_{j=1}^n (d_j-d_{j-1}) \nonumber\\	& = \delta W_{\infty} \sum_{k: 2\mathfrak{R}_k=\lambda_1} \sum_{i,p=1}^{r_k} \theta_{k,i} \overline{\theta_{k,p} } \langle A^{(2)}	\phi_1^{(k,i)}\overline{\phi_1^{(k,p)}} , \widehat{\phi}_1\rangle_\mu \nonumber\\
	& =\delta W_{\infty} \sum_{k: 2\mathfrak{R}_k=\lambda_1}  \langle A^{(2)}	\left|\sum_{i=1}^{r_k} \theta_{k,i} \phi_1^{(k,i)}\right|^2, \widehat{\phi}_1\rangle_\mu
\end{align}
Using the definitions of $F_{Q_1,k}$ in \eqref{def-F(f,k)} and of $\sigma_{cr}^2(Q_1)$ in \eqref{e:varcr}, it is easy to check
that the limit is equal to $ \delta W_\infty \sigma_{cr}^2(Q_1)$, which implies the result of (1).

(2)
Since $X$ and $B$ are independent and $a_j-a_{j-1}$ is uniformly bounded,  to prove (2), it suffices to show that $\mathbb{P}_{\delta_x}$-almost surely,
\begin{align}\label{step_14-2}
	\sup_{j\in \mathbb{N}} \mathbb{E}_{\delta_x} \left[\left(\mathcal{U}_{a_j \delta}^{Q_1} \right)^4  \Big|  \mathcal{F}_{a_{j-1}\delta}\right]  + 	\sup_{j\in \mathbb{N}} \mathbb{E}_{\delta_x} \left[\left(\mathcal{U}_{a_j \delta}^{Q_1} \right)^4 \right] <\infty.
\end{align}
By the definition of $\mathcal{U}_{a_j \delta}^{Q_1}$,
we see that,  conditioned on $\mathcal{F}_{a_{j-1}\delta}$, $\mathcal{U}_{a_j \delta}^{Q_1}$
is the sum of finitely many independent random variables of mean 0.
For  independent random variables $Y_1,..., Y_n$ of mean 0, we have
$$
\mathbb{E}\Big(\sum_{j=1}^n Y_j\Big)^4 =\sum_{j=1}^n \mathbb{E}(Y_j^4) + 3\sum_{ i\neq j} \mathbb{E}(Y_i^2) \mathbb{E}(Y_j^2) \lesssim \sum_{j=1}^n \mathbb{E}(Y_j^4)+ \Big(\sum_{j=1}^n\mathbb{E}(Y_j^2) \Big)^2.
$$
Therefore,
\begin{align}\label{Step19}
	\mathbb{E}_{\delta_x} \left[\left(\mathcal{U}_{a_j \delta}^{Q_1} \right)^4  \Big|  \mathcal{F}_{a_{j-1}\delta}\right] & \lesssim  \sum_{i=1}^{M_{a_{j-1}\delta }}\mathbb{E}_{\delta_x} \left[\left|\mathcal{Y}_j^{Q_1,i} \right|^4 1_{\{|\mathcal{Y}_j^{Q_1,i} | \leq 1 \}} \Big| \mathcal{F}_{a_{j-1}\delta} \right]\nonumber \\ & \qquad  + \bigg(\sum_{i=1}^{
	M_{a_{j-1}\delta}
	}\mathbb{E}_{\delta_x} \left[\left|\mathcal{Y}_j^{Q_1,i} \right|^2 1_{\{|\mathcal{Y}_j^{Q_1,i} | \leq 1 \}} \Big| \mathcal{F}_{a_{j-1}\delta} \right]\bigg)^2,
\end{align}
where in the inequality
we also used the  inequalities $\mathbb{E}\left[Y-\mathbb{E}[Y]\right]^4 \leq 16\mathbb{E}[Y^4]$ and $
\mathbb{E}\left[Y-\mathbb{E}[Y]\right]^2 \leq \mathbb{E}[Y^2] $ for $Y= \mathcal{Y}_j^{Q_1,i}  1_{\{|\mathcal{Y}_j^{Q_1,i} | \leq 1 \}}$.
By \eqref{Up-of-Y}, we have the following upper bound
\[
\left|\mathcal{Y}_j^{Q_1,i} \right|^4 1_{\{|\mathcal{Y}_j^{Q_1,i} | \leq 1 \}}  \leq\left|\mathcal{Y}_j^{Q_1,i} \right|^2 1_{\{|\mathcal{Y}_j^{Q_1,i} | \leq 1 \}}  \leq  e^{\lambda_1 a_{j-1}\delta}\Upsilon_i^2,
\]
where $\Upsilon_i$ are iid copies of $\Upsilon$.
Thus, by the Markov property, we conclude from \eqref{Step19} that
\begin{align}\label{step_12}
	&	\mathbb{E}_{\delta_x} \left[\left(\mathcal{U}_{a_j \delta}^{Q_1} \right)^4  \Big|  \mathcal{F}_{a_{j-1}\delta}\right] \lesssim    e^{\lambda_1 a_{j-1}\delta} \langle \mathbb{E}_{\delta_{\cdot}} (\Upsilon^2) , X_{a_{j-1}\delta}\rangle   + \left( e^{\lambda_1 a_{j-1}\delta} \langle \mathbb{E}_{\delta_{\cdot}} (\Upsilon^2) , X_{a_{j-1}\delta}\rangle \right)^2 .
\end{align}
Combining Lemma \ref{lemma10} and \eqref{upp-second-phi} (with $R= N\delta$), we get $\sup_{j\geq 0} e^{\lambda_1 a_{j-1}\delta}\langle T_{11t_0}b_{t_0}^{1/2}, X_{a_{j-1}\delta} \rangle <\infty$ and $\mathbb{E}_{\delta_x}(\Upsilon^2)\lesssim T_{11t_0}b_{t_0}^{1/2}$.
Moreover, by Lemma \ref{lemma2} (3) and
the inequality $\mathbb{E}_{\delta_x}\left(\langle T_{11t_0}b_{t_0}^{1/2}, X_{t} \rangle ^2\right)\leq \mathbb{E}_{\delta_x}\left(\langle b_{t_0}^{1/2}, X_{t+11t_0} \rangle ^2\right)$,
we know that $e^{\lambda_1 t}\langle T_{11t_0}b_{t_0}^{1/2}, X_{t} \rangle $ is $L^2$ bounded. Thus  \eqref{step_14-2} is valid. The proof is complete.

\hfill$\Box$.

\begin{lemma}\label{LIL-W-cr}
	Let $\delta>0$ and $\ell \in \mathbb{N}$ be fixed, and let $Q_\ell$ be defined in \eqref{Def-of-Q}.
	If $Q_\ell\neq 0$ and  $a_k$ is a syndetic sequence with $a_0=0$ and $a_{k+1}-a_k\in \{1,..,N\}$ for some $N\in\mathbb{N}$, then
	$\mathbb{P}_{\delta_x}$-almost surely,
	\begin{align}\label{LIL-level-l}
			\limsup_{n\to\infty} \frac{	\sum_{(k,i)\in \mathcal{A}_{\ell}}
			\theta_{k,i}^{(\ell)} e^{\lambda_k a_n\delta } \langle  (\Phi_{k,i})^T  e_{k,i}^{(\ell)}, X_{a_n\delta}\rangle }{\sqrt{2 (a_n\delta)^{1+2\tau(Q_\ell)} \log \log (a_n\delta)}} =\sqrt{\sigma_{cr}^2 (Q_\ell) W_\infty }.
	\end{align}
\end{lemma}
\textbf{Proof: }
We first prove \eqref{LIL-level-l} for $\ell=1$. Combining Lemma \ref{cor2} and Markov inequality, for each $\varepsilon,\delta>0$,  it is easy to see that $\mathbb{P}_{\delta_x}$-almost surely,
 \begin{align}
 	& \sum_{j=1}^\infty \frac{1}{(s_j^{Q_1})^4} \mathbb{E}_{\delta_x} \left( \mathcal{X}_j^4\big| \mathcal{G}_{j-1}^Z\right) +\sum_{j=1}^\infty \mathbb{E}_{\delta_x} \left( |\mathcal{X}_j| 1_{\{|\mathcal{X}_j| >\sqrt{j}\}}\right) <\infty
 \end{align}
 and
 \[
 \lim_{n\to\infty} \frac{1}{(s_n^{Q_1})^2}\sum_{j=1}^n \mathbb{E}_{\delta_x} \left( \mathcal{X}_j^2 1_{\{|\mathcal{X}_j| >\sqrt{j}\}}\big| \mathcal{G}_{j-1}^Z\right)=0.
 \]
 Therefore, the martingale $\mathcal{Z}_j$ satisfies the condition of \cite[Theorem 4.7]{HH1980}.
 For $t\in \left[\frac{(s_i^{Q_1})^2}{(s_n^{Q_1})^2}, \frac{(s_{i+1}^{Q_1})^2}{(s_n^{Q_1})^2}\right)$ and $i\leq n-1$, define
 \begin{align}
 	\beta_{\mathcal{Z}, n}(t):= \frac{1}{\sqrt{2 (s_n^{Q_1})^2 \log \log (s_n^{Q_1})^2}}\left(\mathcal{Z}_i + \frac{( t(s_{n}^{Q_1})^2 -(s_i^{Q_1})^2)(\mathcal{Z}_{i+1}-\mathcal{Z}_i  )}{(s_{i+1}^{Q_1})^2 -(s_i^{Q_1})^2} \right).
 \end{align}
 We define $\beta_{Z,n}$ in  the same way, with $\mathcal{Z}_i$ replaced by $Z_i$.
 Combining $|\mathbb{E}(X|\mathcal{F})| \leq \mathbb{E}(|X| |\mathcal{F})$ and the fact that $\mathbb{E}_{\delta_x}(\mathcal{Y}_j^{Q_1, i}| \mathcal{F}_{a_{j-1}\delta})=0$, we get
 \begin{align}
 	& \sum_{j=1}^\infty \mathbb{E}_{\delta_x}\left( |Z_j -Z_{j-1}- \mathcal{X}_j| \right)  \leq 2 \sum_{j=1}^\infty \mathbb{E}_{\delta_x}\bigg(   \sum_{i=1}^{M_{a_{j-1}\delta}} |\mathcal{Y}_j^{Q_1,i} |1_{\{|\mathcal{Y}_j^{Q_1,i} | >1 \} }  \bigg) \nonumber\\
 	&  = \sum_{j=1}^\infty \mathbb{E}_{\delta_x} \left( \langle \mathbb{E}_{\delta_{\cdot}} ( |\mathcal{Y}_j^{Q_1} |1_{\{|\mathcal{Y}_j^{Q_1} | >1 \} } ), X_{a_{j-1}\delta}\rangle \right) \leq \sum_{j=1}^\infty e^{\lambda_1 a_{j-1}\delta}\mathbb{E}_{\delta_x} \left( \langle \mathbb{E}_{\delta_{\cdot}} ( \Upsilon 1_{\{ \Upsilon>e^{-\lambda_1 a_{j-1}\delta} \} } ), X_{a_{j-1}\delta}\rangle \right),
 \end{align}
 where in the last equality we used \eqref{Up-of-Y}. Therefore, repeating the same argument  for \eqref{step_17}
 with $Y^f(s,r)$ replaced by $\Upsilon$, we conclude that $\mathbb{P}_{\delta_x}$-almost surely,
 \begin{align}
 	\sup_{n\in\mathbb{N}}|Z_n-\mathcal{Z}_n|\leq \sum_{j=1}^\infty  |Z_j -Z_{j-1}- \mathcal{X}_j| <\infty,
 \end{align}
 which implies that
  \begin{align} \label{bounds-of-Zs}
 \lim_{n\to\infty}\sup_{t\in [0,1]} |\beta_{\mathcal{Z},n}(t)- \beta_{Z,n}(t)|=0 .
 \end{align}
 Combining \eqref{bounds-of-Zs} and \cite[Theorems 4.7 and 4.8]{HH1980}, we get
 \begin{align}
 	\limsup_{n\to\infty} \frac{  \mathcal{W}_{a_n\delta}+\varepsilon B_{a_n} } {\sqrt{2 (a_n\delta)^{1+2\tau(Q_1)} \log \log (a_n\delta)}} =\sqrt{\sigma_{cr}^2 (Q_1) W_\infty  +\varepsilon^2 },\quad \mathbb{P}_{\delta_x}\textup{-a.s.},
 \end{align}
 and
 \begin{align}\label{Useful-result}
 	\mbox{$\{\beta_{Z,n}\}_{n>3}$ is relatively compact in $C[0,1]$ with
	closure equal to $\mathcal{K}$ a.s.,
	}
 \end{align}
 where $\mathcal{K}$ is the set of absolutely continuous function $z(t)\in C[0,1]$ with $z(0)=0$ and $\int_0^1 (z'(t))^2 \mathrm{d}t\leq 1$.
Now letting  $\varepsilon\to 0$, we get that \eqref{LIL-level-l} holds for $\ell=1$.

For any $(k,i)\in \mathcal{A}_{2}$ and  $2\leq \ell \leq d_{k,i}$,  set $\phi_\ell^{(k,i)}:= (\Phi_{k,i})_\ell$.
 Taking $g_k= T_{-p\delta}\phi_\ell^{(k,i)}$
  for all $p \in \{1,...,N\}$ in  \eqref{e2}, we see that there exists a random variable  $\mathcal{U}$  such that for large $n$,
 \[
 \sup_{2\leq q\leq \ell}\sup_{p\in \{0, 1,...,N\}} \left|e^{\lambda_k n\delta}\langle T_{-p\delta}
 \phi_q^{(k,i)},
  X_{(n+1)\delta} \rangle - e^{\lambda_k n\delta}\langle T_{(-p+1)\delta}
  \phi_q^{(k,i)},
   X_{n\delta} \rangle \right|\leq   \mathcal{U}   \sqrt{\log n}.
 \]
 Therefore, when $n$ is large enough (say, $n\geq R$), for all $2\leq q\leq \ell$,
 \begin{align}\label{Ite-of-W}
 	& \left| e^{\lambda_k a_n\delta}\langle T_{-a_n\delta+a_{n-1}\delta}
 	\phi_q^{(k,i)},
 	X_{a_n\delta} \rangle - e^{\lambda_k a_n\delta}\langle
 	\phi_q^{(k,i)},
 	X_{a_{n-1}\delta} \rangle \right|\nonumber\\
 	&\leq \sum_{j=a_{n-1}+1}^{a_n}  \left| e^{\lambda_k a_n\delta}\langle T_{(a_{n-1}-j)\delta}
 	\phi_q^{(k,i)},
 	X_{j\delta} \rangle - e^{\lambda_k a_n\delta}\langle T_{(a_{n-1}-j+1)\delta}
 	\phi_q^{(k,i)},
 	X_{(j-1)\delta} \rangle \right|\nonumber\\
 	& \leq N \mathcal{U}  \sqrt{\log n}.
 \end{align}
 Write ${\mathcal{W}}_t^{(q,k,i)}:= e^{\lambda_k t} \langle \phi_q^{(k,i)}, X_t\rangle $ for simplicity. It is routine to check that
 \begin{align}
 	& e^{\lambda_k a_n\delta}\langle T_{-(a_n-a_{n-1})\delta}
 	\phi_q^{(k,i)}, X_{a_n\delta} \rangle - e^{\lambda_k a_n\delta}\langle
 	\phi_q^{(k,i)},	X_{a_{n-1}\delta} \rangle\nonumber\\
 & = e^{\lambda_k a_{n}\delta}\langle (\Phi_k)^T
 	J_{k,i}((a_n-a_{n-1})\delta)^{-1}e_{k,i}^{(q)}, X_{a_{n}\delta} \rangle  -
 	{\mathcal{W}}_{a_{n-1}\delta}^{(q,k,i)}
 	 \nonumber\\
 	& =
 	{\mathcal{W}}_{a_n\delta}^{(q,k,i)}- {\mathcal{W}}_{a_{n-1}\delta}^{(q, k,i)} +\sum_{u=1}^{q-1} \frac{(a_{n-1}\delta-a_{n}\delta)^u}{u!} {\mathcal{W}}_{a_{n}\delta}^{(q-u,k,i)}.
 \end{align}
 Plugging the above equation back to \eqref{Ite-of-W},   we see that for $n\geq R+1$,
 \begin{align}\label{dif-w-2}
 	&
 	\bigg|{\mathcal{W}}_{a_n\delta}^{(q,k,i)}-{\mathcal{W}}_{a_{R}\delta}^{(q, k,i)} + \sum_{j=R+1}^n \sum_{u=1}^{q-1} \frac{(a_{j-1}\delta-a_{j}\delta)^u}{u!} {\mathcal{W}}_{a_{j}\delta}^{(q-u,k,i)} \bigg| \leq N
	\mathcal{U}
	n \sqrt{\log n}.
 \end{align}
 Recall the definition of $Z_n$ in  \eqref{Def-of-Z-n}.
 Define $ \mathcal{S}_{a_n\delta}^{(1, 1)}: =\mathcal{W}_{a_n\delta}, \mathcal{S}_{a_n\delta}^{(2,1)}: =Z_n$,  $\mathcal{S}_{a_n\delta}^{(3,1)}: =\varepsilon B_{a_n}$ and for $j\in \{1,2,3\}$,
 \begin{align}\label{Def-of-S-q}
 	\mathcal{S}_{a_n\delta}^{(j, q)}:=- \sum_{k=1}^{n} \sum_{u=1}^{q-1}\frac{(a_{k-1}\delta-a_{k}\delta)^u}{u!} \mathcal{S}_{a_{k}\delta}^{(j, q-u)}, \quad 2\leq q\leq \ell.
 \end{align}
 We claim that,  $\mathbb{P}_{\delta_x}$-almost surely,
 \begin{align}\label{coupling4}
 	\limsup_{n\to\infty} \frac{|
 		\sum_{(k,i)\in \mathcal{A}_{\ell} }
 		\theta_{k,i}^{(\ell)} e^{\lambda_k a_n\delta } \langle  (\Phi_{k,i})^T  e_{k,i}^{(\ell)},
 		X_{a_n\delta}\rangle
 		-\mathcal{S}_{a_n\delta}^{(1,\ell)} |}{n^{\ell-1}\sqrt{\log n}}<\infty.
 \end{align}
Note that
 $\mathcal{S}_{a_n\delta}^{(1,1)}=\mathcal{W}_{a_n\delta}=\sum_{(k,i)\in \mathcal{A}_{\ell} } e^{\lambda_k t} \theta^{(\ell)}_{k,i} \langle \phi_{\ell}^{(k,i)}, X_{a_n\delta}\rangle =: \sum_{(k,i)\in \mathcal{A}_{\ell} } \theta_{k,i}^{(\ell)} \mathcal{S}_{a_n\delta}^{(1,1,k, i)}$.
To prove \eqref{coupling4}, it suffices to show that  for each pair $(k,i)\in \mathcal{A}_{\ell} $,
 \begin{align}\label{coupling5}
 	 		\limsup_{n\to\infty} \frac{| \mathcal{W}_{a_n\delta}^{(\ell, k,i)} -	\mathcal{S}_{a_n\delta}^{(1,\ell,k, i)}|}{n^{\ell-1}\sqrt{\log n}}=
 		 			\limsup_{n\to\infty} \frac{| e^{\lambda_k a_n\delta } \langle  (\Phi_{k,i})^T  e_{k,i}^{(\ell)},	X_{a_n\delta}\rangle	-
 			\mathcal{S}_{a_n\delta}^{(1,\ell,k, i)}
 			|}{n^{\ell-1}\sqrt{\log n}}<\infty,
 \end{align}
where
$\mathcal{S}_{a_n\delta}^{(1, q,k, i)}$ is
defined in the same way as \eqref{Def-of-S-q} with  $\mathcal{S}_{a_{k}\delta}^{(j, q-u)}$ replaced by
$\mathcal{S}_{a_{k}\delta}^{(j, q-u,k, i)}$.
 If $\ell=1$, then $ e^{\lambda_k a_n\delta } \langle  (\Phi_{k,i})^T  e_{k,i}^{(\ell)},
 X_{a_n\delta}\rangle
 -\mathcal{S}_{a_n\delta}^{(1,\ell,k, i)} =0$.  Suppose that \eqref{coupling5} holds for
   $\ell=1,...,m$, then for $\ell=m+1$ and $(k,i)\in \mathcal{A}_{m+1}
   \subset \mathcal{A}_{m}$,
  by \eqref{dif-w-2} and the definition of $\mathcal{S}_{a_n\delta}^{(1,q,k,i)}$, we have
 \begin{align}
 \limsup_{n\to\infty} \frac{|
 		{\mathcal{W}}_{a_n\delta}^{(m+1,k,i)}
 		-\mathcal{S}_{a_n\delta}^{(1,m+1, k,i)} |}{n^{\ell}\sqrt{\log n}}
 		 \leq &	\limsup_{n\to\infty} \frac{|
 		{\mathcal{W}}_{a_n\delta}^{(m+1,k,i)}+
 		\sum_{j=1}^n \sum_{u=1}^{m} \frac{(a_{j-1}\delta-a_{j}\delta)^u}{u!} {\mathcal{W}}_{a_{j}\delta}^{(\ell+1-u,k,i)} |}{n^{\ell}\sqrt{\log n}} \nonumber\\
 		& + \limsup_{n\to\infty}\sum_{j=1}^n \sum_{u=1}^m \frac{(a_{j}\delta-a_{j-1}\delta)^u}{u!}  \frac{|
 			{\mathcal{W}}_{a_{j}\delta}^{(\ell+1-u,k,i)} - \mathcal{S}_{a_j\delta}^{(1,\ell+1-u,k, i)} |}{n^{\ell}\sqrt{\log n}}\nonumber\\
 			 \leq & N\mathcal{U}
			 + \limsup_{n\to\infty} \sum_{u=1}^m \frac{(N\delta)^u}{u!}  \sum_{j=1}^n\frac{|	{\mathcal{W}}_{a_{j}\delta}^{(\ell+1-u,k,i)} - \mathcal{S}_{a_j\delta}^{(1,\ell+1-u,k, i)} |}{n^{\ell}\sqrt{\log n}}.
 \end{align}
By induction, there is a random variable
$\mathcal{U}'$ such that
 $\sup_{1\leq u\leq m} |
 \widehat{\mathcal{W}}_{a_{j}\delta}^{(\ell+1-u,k,i)} - \mathcal{S}_{a_j\delta}^{(1,\ell+1-u,k, i)} |\leq   \mathcal{U}'
 j^{\ell-1} \sqrt{\log j}\leq  \mathcal{U}'  n^{\ell-1}\sqrt{\log n}$. Plugging this back to the above display, we  obtain
 \begin{align}
 	& \limsup_{n\to\infty} \frac{|
 		{\mathcal{W}}_{a_n\delta}^{(m+1,k,i)}
 		-\mathcal{S}_{a_n\delta}^{(1,m+1,k, i)} |}{n^{m}\sqrt{\log n}}  \leq N	\mathcal{U}+	\mathcal{U}'	\sum_{u=1}^m \frac{(N\delta)^u}{u!} <\infty,
 \end{align}
 which implies \eqref{coupling5} for $m+1$. Therefore,  \eqref{coupling5} holds by induction.

 Let $\widehat{J}_\ell(t)$ be the $\ell\times \ell$ matrix with $(\widehat{J}_\ell(t))_{a,b}=1_{\{b\geq a\}} t^{b-a}/ (b-a)!$.
 Then from the definition of ${\mathcal{S}}_{a_n\delta}^{(j, \ell)} $, for any $j\in \{1,2,3\}$,
 \begin{align}
 	& \left(\mathcal{S}_{a_n\delta}^{(j, q)} ,1\leq q\leq \ell \right)= \left(\mathcal{S}_{a_{n-1}\delta}^{(j, q)} ,1\leq q\leq \ell \right) \widehat{J}_\ell((a_n-a_{n-1})\delta)  	+ (  \mathcal{S}_{a_n\delta}^{(j,1)}- \mathcal{S}_{a_{n-1}\delta}^{(j,1)}, 0,...,0 )\nonumber\\	&=\cdots= \sum_{k=1}^n ( \mathcal{S}_{a_k\delta}^{(j,1)}-\mathcal{S}_{a_{k-1}\delta}^{(j,1)}, 0,...,0 ) \widehat{J}_\ell((a_n-a_k)\delta),
 \end{align}
 which implies that
 \begin{align}\label{iteration}
 	& \mathcal{S}_{a_{n}\delta}^{(j, \ell)} =  \sum_{k=1}^n \frac{(a_n\delta-a_k\delta)^{\ell-1}}{(\ell-1)!} \left(\mathcal{S}_{a_k\delta}^{(j,1)}-\mathcal{S}_{a_{k-1}\delta}^{(j,1)}\right) .
 \end{align}
 Taking $j=3$ in the above inequality, we get
 \begin{align}
    \left| \mathcal{S}_{a_{n}\delta}^{(2, \ell)} -\mathcal{S}_{a_{n}\delta}^{(1, \ell)}  \right| & =\left| \mathcal{S}_{a_{n}\delta}^{(3, \ell)}\right| = \left| \frac{\varepsilon \delta^{\ell-1}}{(\ell-1)!} \sum_{k=1}^{n-1} B_{a_k}(a_n-a_k)^{\ell-1}\right|\leq \frac{\varepsilon \delta^{\ell-1} a_n^{\ell-1}}{(\ell-1)!}\sum_{k=1}^{n-1} |B_{a_k}|.
 \end{align}
 According to the LIL for Brownian motion,
 there exists a random variable  $\mathcal{U}$
 such that $|B_{a_k}|\leq   \mathcal{U}
 \sqrt{a_k \log \log a_k}\leq    \mathcal{U} \sqrt{a_n \log \log a_n}$ almost surely. Combining  this with $a_n\in [n, Nn]$, we conclude that almost surely,
 \begin{align}\label{coupling2}
 	\lim_{\varepsilon\to 0} \limsup_{n\to\infty} \frac{| \mathcal{S}_{a_n\delta}^{(1,\ell)}-{\mathcal{S}}_{a_n\delta}^{(2, \ell)} |}{\sqrt{2n^{1+2(\ell-1)}\log\log n}} =0.
 \end{align}
 For ${\mathcal{S}}_{a_n\delta}^{(2, \ell)}$, by \eqref{iteration}, we have
 \begin{align}
 		{\mathcal{S}}_{a_n\delta}^{(2, \ell)} 	= (a_n\delta)^{\ell-1} \frac{\sqrt{2(s_n^{Q_1})^2\log\log (s_n^{Q_1})^2}}{(\ell-1)!}\sum_{j=1}^{n-1} \beta_{Z,n}\left(\frac{(s_j^{Q_1})^2 }{(s_n^{Q_1})^2}\right) \left(\left(1-\frac{a_j}{a_n}\right)^{\ell-1}- \left(1-\frac{a_{j+1}}{a_n}\right)^{\ell-1}\right)  .
 \end{align}
 According to \eqref{Useful-result}, for any $\gamma>0$ and $\zeta \in \mathcal{K}$, almost surely $\sup_{t\in [0,1]}|\beta_{Z,n}(t)-\zeta (t)|<\gamma$
 for  infinitely many $n$. Now
 we assume that $n$ is large enough such that $\sup_{1\leq j\leq n}|(s_j^{Q_1})^2/ (s_n^{Q_1})^2 - a_j/ a_n|<\gamma$. Therefore,
 since
 $\sum_{j=1}^{n-1}  \left|\left(1-\frac{a_j}{a_n}\right)^{\ell-1}- \left(1-\frac{a_{j+1}}{a_n}\right)^{\ell-1}\right| <1 $, when $n$ is large enough, we may replace $\beta_{Z,n}$ by $\zeta $ and $(s_j^{Q_1})^2/
 (s_n^{Q_1})^2
 $ by $ a_j/ a_n$,
and the resulting  error is at most $2\gamma$.
 Since $\gamma$ is arbitrary, by Lemma \ref{cor2}(1),
 \begin{align}\label{coupling3}
 	& \lim_{\varepsilon\to 0} \limsup_{n\to\infty} \frac{{\mathcal{S}}_{a_n	\delta}^{(2,\ell)}}{\sqrt{2 (a_n\delta)^{1+2(\ell-1)} \log \log (a_n\delta)}} \nonumber\\
 	& =\lim_{\varepsilon\to 0}  \sqrt{\frac{\varepsilon^2+\delta W_\infty \sigma_{cr}^2(Q_1)}{\delta}} \frac{1}{(\ell-1)!} \sup_{\zeta \in \mathcal{K}} \lim_{n\to\infty} \sum_{j=1}^{n-1} \zeta \left(a_j/a_n\right) \left(\left(1-\frac{a_j}{a_n}\right)^{\ell-1}- \left(1-\frac{a_{j+1}}{a_n}\right)^{\ell-1}\right) \nonumber\\
 	& =  \sqrt{W_\infty \sigma_{cr}^2(Q_1)} \frac{1}{(\ell-1)!} \sup_{\zeta \in \mathcal{K}}  \int_0^1 (1-t)^{\ell-1} \zeta '(t)\mathrm{d}t.
 \end{align}
 According to \cite[p219]{Strassen}, $\sup_{\zeta \in \mathcal{K}}  \int_0^1 (1-t)^{\ell-1} \zeta '(t)\mathrm{d}t= \sqrt{\int_0^1 (1-t)^{2\ell-2}\mathrm{d} t  }= (2\ell-1)^{-1/2} $. Therefore,
 by \eqref{coupling4}, \eqref{coupling2} and \eqref{coupling3},  we obtain
 \begin{align}
 	\limsup_{n\to\infty} \frac{
 		\sum_{(k,i)\in \mathcal{A}_{\ell}}
 		\theta_{k,i}^{(\ell)} e^{\lambda_k a_n\delta } \langle  (\Phi_{k,i})^T  e_{k,i}^{(\ell)},
 		X_{a_n\delta}\rangle
 	}{\sqrt{2 (a_n\delta)^{1+2\tau(Q_\ell)} \log \log (a_n\delta)}} =\sqrt{\frac{(2\ell-1)^{-1}}{((\ell-1)!)^2}\sigma_{cr}^2 (Q_1) W_\infty }.
 \end{align}
 An elementary calculation yields that $\frac{(2\ell-1)^{-1}}{((\ell-1)!)^2}\sigma_{cr}^2 (Q_1)=\sigma_{cr}^2 (Q_\ell)$, which completes the proof of the lemma.

\hfill$\Box$

\begin{cor}\label{LIL-conti-mar}
	$\mathbb{P}_{\delta_x}(\cdot| \mathcal{E}^c)$-almost surely,
	\begin{align}
		\limsup_{t\to\infty} \frac{\sum_{(k,i)\in \mathcal{A}_{\ell} }
			\theta_{k,i}^{(\ell)} e^{\lambda_k t } \langle  (\Phi_{k,i})^T  e_{k,i}^{(\ell)},	X_{t}\rangle
			}{\sqrt{2 t^{1+2\tau(Q_\ell)} \log \log t}} =\sqrt{\sigma_{cr}^2 (Q_\ell) W_\infty }.
	\end{align}
\end{cor}
\textbf{Proof: }
By Lemma \ref{LIL-W-cr}, it suffices to prove that
\begin{align}\label{goal2}
	\lim_{\delta\to 0} \limsup_{n\to\infty} \sup_{t\in [n\delta, (n+1)\delta)} \frac{
		\sum_{(k,i)\in \mathcal{A}_{\ell} }
		|\theta_{k,i}^{(\ell)} | \left|e^{\lambda_k t } \langle (\Phi_{k,i})^T  e_{k,i}^{(\ell)} , X_{t}\rangle -e^{\lambda_k n\delta } \langle  (\Phi_{k,i})^T  e_{k,i}^{(\ell)}, X_{n\delta}\rangle \right|
		 }{\sqrt{2 (n\delta)^{1+2\tau(Q_\ell)} \log \log (n\delta)}}  =0.
\end{align}
By  Lemma \ref{lem:conti-mart},
we only need to show that for each
$(k,i)\in \mathcal{A}_{d}$ and $\ell \leq d_{k,i}$,
\begin{align}
		\lim_{\delta\to 0} \limsup_{n\to\infty} \sup_{t\in [n\delta, (n+1)\delta)} \frac{  \left|e^{\lambda_k t } \langle  T_{t-n\delta}
			(\Phi_{k,i})^T  e_{k,i}^{(\ell)}, X_{n\delta}\rangle -e^{\lambda_k n\delta } \langle
			(\Phi_{k,i})^T  e_{k,i}^{(\ell)},
			 X_{n\delta}\rangle \right|  }{\sqrt{2 (n\delta)^{1+2\tau(Q_\ell)} \log \log (n\delta)}}  =0.
\end{align}
Note that
\begin{align}
	&  \frac{  \left|e^{\lambda_k t } \langle  T_{t-n\delta}	(\Phi_{k,i})^T  e_{k,i}^{(\ell)}, X_{n\delta}\rangle -e^{\lambda_k n\delta } \langle  	(\Phi_{k,i})^T  e_{k,i}^{(\ell)}, X_{n\delta}\rangle \right|  }{\sqrt{2 (n\delta)^{1+2\tau(Q_\ell)} \log \log (n\delta)}}  \nonumber\\	&= \frac{ e^{\lambda_1 n\delta/2} \left| \langle  \Phi_{k,i}^T(J_{k,i}(t-n\delta)-I) e_{k,i}^{(\ell)} , X_{n\delta}\rangle  \right|  }{\sqrt{2 (n\delta)^{1+2(\ell-1)} \log \log (n\delta)}}= \frac{ e^{\lambda_1 n\delta/2} \left| \sum_{j=1}^{\ell-1} \langle \frac{(t-n\delta)^{\ell-j}}{(\ell-j)!}  	(\Phi_{k,i})^T  e_{k,i}^{(j)} , X_{n\delta}\rangle  \right|  }{\sqrt{2 (n\delta)^{1+2(\ell-1)} \log \log (n\delta)}}\nonumber\\	& \leq \sum_{j=1}^{\ell-1}  \frac{\delta^{\ell-j}}{(\ell-j)!}  \frac{ e^{\lambda_1 n\delta/2} \left|  \langle 	(\Phi_{k,i})^T  e_{k,i}^{(j)} , X_{n\delta}\rangle  \right|  }{\sqrt{2 (n\delta)^{1+2(\ell-1)} \log \log (n\delta)}} = \sum_{j=1}^{\ell-1}  \frac{ n^{- (\ell-j)}}{(\ell-j)!}  \frac{ e^{\lambda_1 n\delta/2} \left|  \langle 	(\Phi_{k,i})^T  e_{k,i}^{(j)}, X_{n\delta}\rangle  \right|  }{\sqrt{2 (n\delta)^{1+2(j-1)} \log \log (n\delta)}}.
\end{align}
Combining the above with Lemma \ref{LIL-W-cr}, we get the desired result.
	
\hfill$\Box$

\begin{lemma}\label{lem:LIL-without-Im}
		$\mathbb{P}_{\delta_x}(\cdot| \mathcal{E}^c)$-almost surely,
	\begin{align}
		\limsup_{t\to\infty} \frac{  e^{\lambda_1 t/2}\langle Q_\ell,X_t\rangle  }{\sqrt{2 t^{1+2\tau(Q_\ell)} \log \log t}} =\sqrt{\sigma_{cr}^2 (Q_\ell) W_\infty }.
	\end{align}
\end{lemma}
\textbf{Proof: }
Define $\mathcal{S}:=\left\{ \Theta= (\Theta_{k}: \exists i \mbox{ s. t. }  (k,i)\in \mathcal{A}_{\ell} )^T :  \Theta_k\in \mathbb{C}, |\Theta_k|=1\ \mbox{and}\ \Theta_k= \Theta_{k'}\right\}$.
For any $\Theta\in \mathcal{S}$, define
\[
\Theta \star Q_\ell := \sum_{(k,i)\in \mathcal{A}_{\ell}}
\theta_{k,i}^{(\ell)}  (\Phi_{k,i})^T  (\Theta_{k} e_{k,i}^{(\ell)} ) =  \sum_{k: 2\mathfrak{R}_k=\lambda_1} \Theta_k \sum_{i=1}^{r_k} 1_{\{\ell \leq d_{k,i}\}}\theta_{k,i}^{(\ell)} (\Phi_k(x))^T  e_k^{(\nu_{k,i-1}+\ell)}.
\]
 Then it is easily seen
 from \eqref{lim-of-var} and the identity $\frac{(2\ell-1)^{-1}}{((\ell-1)!)^2}\sigma_{cr}^2 (Q_1)=\sigma_{cr}^2 (Q_\ell)$ that
  $\sigma_{cr}^2(Q_\ell)= \sigma_{cr}^2(\Theta\star Q_\ell)$ for any $\Theta\in \mathcal{S}$.

 Set $\mathcal{P}_{k,t}:= \sum_{i=1}^{r_k} 1_{\{\ell \leq d_{k,i}\}}\theta_{k,i}^{(\ell)} \langle (\Phi_k(x))^T  e_k^{(\nu_{k,i-1}+\ell)}, X_t\rangle$. For each fixed pair $(k_0, k_0')$,
applying Corollary \ref{LIL-conti-mar} (with $\Theta_k=\Theta_{k'}= \pm 1 $ for $k\neq k_0, k_0'$ and $Q_\ell$ replaced by $\Theta \star Q_\ell$),
we obtain that
  there exists some constant $\Gamma_1= \Gamma_1(\#\{k: 2\mathfrak{R}_k=\lambda_1 \})$ such that almost surely,
\begin{align}\label{e10}
	\limsup_{t\to\infty} \frac{ \sup_{k: 2\mathfrak{R}_k=\lambda_1}
		| e^{\lambda_k t}  \mathcal{P}_{k,t} |
		 }{\sqrt{2 t^{1+2(\ell-1)} \log \log t}}
		=\limsup_{t\to\infty} \frac{ \sup_{k: 2\mathfrak{R}_k=\lambda_1}
			e^{\lambda_1 t/2}|   \mathcal{P}_{k,t} |
		 }{\sqrt{2 t^{1+2(\ell-1)} \log \log t}}  \leq \Gamma_1\sqrt{\sigma_{cr}^2 (Q_\ell) W_\infty }.
\end{align}
For any $\varepsilon>0$, since $\mathcal{S}$ is compact, we may find a finite subset
$\mathcal{R}$ of $\mathcal{S}$ such that for any $\Theta\in \mathcal{S}$, there exists
$\mathcal{R}^u= (\mathcal{R}^u_k: \exists\ i\mbox{ s. t.} (k,i)\in \mathcal{A}_{\ell} )\in \mathcal{R}$ such that $|\Theta -\mathcal{R}^u|<\varepsilon$. Taking $Q_\ell=\mathcal{R}^u\star Q_\ell$ in
Corollary \ref{LIL-conti-mar}, we obtain that
\begin{align}\label{e11}
	& \limsup_{t\to\infty}\sup_{\mathcal{R}^u\in \mathcal{R}} \frac{ \sum_{k: 2\mathfrak{R}_k=\lambda_1}  e^{\lambda_1 t/2} \left( e^{\mathrm{i}\mathfrak{I}_k t}\mathcal{R}^u_k\right)  \mathcal{P}_{k,t} }{\sqrt{2 t^{1+2(\ell-1)} \log \log t}} = \sqrt{\sigma_{cr}^2 (Q_\ell) W_\infty }.
\end{align}
Suppose that $\mathcal{R}^{t}\in \mathcal{R}$ satisfies $| \Theta^*(t)-\mathcal{R}^t|<\varepsilon$ where
$\Theta^*(t)=( e^{-\mathrm{i}\mathfrak{I}_k t}:  \exists i \mbox{ s. t. } (k,i)\in \mathcal{A}_{\ell} )^T$.
By \eqref{e10} and \eqref{e11}, we have
\begin{align}
		& \limsup_{t\to\infty} \frac{  e^{\lambda_1 t/2}\langle Q_\ell,X_t\rangle  }{\sqrt{2 t^{1+2\tau(Q_\ell)} \log \log t}} \nonumber\\
		&  \leq \limsup_{t\to\infty}\sup_{\mathcal{R}^u\in \mathcal{R}} \frac{ \sum_{k: 2\mathfrak{R}_k=\lambda_1}  e^{\lambda_1 t/2} \left( e^{\mathrm{i}\mathfrak{I}_k t}\mathcal{R}^u_k\right)  \mathcal{P}_{k,t}  }{\sqrt{2 t^{1+2(\ell-1)} \log \log t}}  + \varepsilon \#\{k: 2\mathfrak{R}_k=\lambda_1 \} \Gamma_1\sqrt{\sigma_{cr}^2 (f) W_\infty }		\nonumber\\		& = \varepsilon\#\{k: 2\mathfrak{R}_k=\lambda_1 \} \Gamma_1\sqrt{\sigma_{cr}^2 (Q_\ell) W_\infty } + \sqrt{\sigma_{cr}^2 (Q_\ell) W_\infty }.
\end{align}
Taking $\varepsilon \to 0$, we arrive at the upper bound.

Now we prove the lower bound. For any $\varepsilon>0$, by \cite[Theorem 1.21]{Furstenberg}, there exists a syndetic sequence $\{a_{n}: n\in \mathbb{N}\}$ such that
$
\sup_{ k: 2\mathfrak{R}_k=\lambda_1} |e^{\mathrm{i}\mathfrak{I}_k a_n}-1|<\varepsilon.
$
Thus, together with Lemma \ref{LIL-W-cr}, we obtain the  lower bound
\begin{align}\label{lo1}
		& \limsup_{t\to\infty} \frac{ e^{\lambda_1 t/2}\langle Q_\ell,X_t\rangle  }{\sqrt{2 t^{1+2\tau(Q_{\ell})} \log \log t}} \geq \limsup_{n\to\infty} \frac{ e^{\lambda_1 a_n/2}\langle Q_\ell,X_{a_n}\rangle  }{\sqrt{2 {a_n}^{1+2\tau(Q_{\ell})} \log \log a_n}} \nonumber\\
		& \geq \sqrt{\sigma_{cr}^2(Q_\ell)
			 W_\infty } - \varepsilon\#\{k: 2\mathfrak{R}_k=\lambda_1 \} \Gamma_1\sqrt{\sigma_{cr}^2 (Q_\ell) W_\infty }.
\end{align}
Taking $\varepsilon\to 0$, we arrive at the lower bound. The proof is complete.

\hfill$\Box$

\textbf{Proof of Theorem \ref{thm3}:}	
Define $\ell= 1+\tau(f_{cr})$, and,  for $1\leq q\leq \ell$, let $Q_q$ be defined as in \eqref{Def-of-Q}, then $f_{cr}=\sum_{q=1}^\ell Q_q$.
Applying Lemma \ref{lem:LIL-without-Im}  to each $Q_q$ and using the fact that $\tau(f_{cr})=\tau(Q_\ell)$, we get the desired result.

\hfill$\Box$

\section{Proof of Theorem \ref{thm5}}\label{ss: 5}

In this section, we always assume that {\bf(H1)}--{\bf(H3)} hold.

\subsection{Proof of Theorem \ref{thm5}}\label{ss:5.2}

In this subsection, we first
prove Theorem \ref{thm5} using the following Proposition \ref{prop1},
and then give  the proof of Proposition \ref{prop1}.

\begin{prop}\label{prop1}
	Let $f\in {\cal T}$ with  $\mathfrak{R}_{\gamma(f)}>0$.
	Suppose in addition that {\bf(H4)} holds. Then

	\begin{align}
		\limsup_{t\to\infty} \frac{e^{\lambda_1 t/2} \left|\langle f, X_t \rangle\right|}{\sqrt{2\log t}}\leq  18 \sqrt{\sigma_{sm}^2(f) W_\infty} ,\quad \mathbb{P}_{\delta_x}\left(\cdot| \mathcal{E}^c\right)
			\mbox{-a.s.}
	\end{align}
\end{prop}

Note that $\mathfrak{R}_{\gamma(f)}>0$ implies
 $\lambda_1< 2\mathfrak{R}_{\gamma(f)}$, which corresponds to the small branching rate case.

\bigskip

\textbf{Proof of Theorem \ref{thm5}:} We only give the proof of  Theorem \ref{thm2} here, the proof for Theorem \ref{thm4} is similar.
Combining Lemma \ref{lemma2} (1) and {\bf(H4)}{\bf(a)},
we have
$\sigma_{sm}^2(f)\lesssim\Vert f\Vert^2_2 +  \langle |f|^2, \widehat{\phi}_1\rangle_\mu\lesssim \Vert f\Vert_2^2.$
Therefore, by Proposition \ref{prop1}, for any $f\in {\cal T}$ with $\mathfrak{R}_{\gamma(f)}>0$, there exists a constant $C$ independent of $f$ such that
\begin{align}\label{Step28}
	\limsup_{t\to\infty} \frac{e^{\lambda_1 t/2} \left|\langle f, X_t \rangle\right|}{\sqrt{2\log t}}\leq  C\sqrt{W_\infty} \Vert f\Vert_2 ,\quad \mathbb{P}_{\delta_x}\left(\cdot| \mathcal{E}^c\right)\textup{-a.s.}
\end{align}
Now for any $f\in {\cal T}$,  we write $f= f_{main}+ f_{rest}$, where
\[
f_{main} = \sum_{k\in \mathbb{I}: k\leq N} e^{-\lambda_k r} (\Phi_k)^TD_k(r) v_k, \quad	f_{rest}= f- f_{main}\in {\cal T}
\]
and $N$ is a large integer such that $\mathfrak{R}_k>0$ for all $k>N$.
Applying Theorem \ref{thm2} to $f_{main}$ and  \eqref{Step28} to $f_{rest}$, we see that $\mathbb{P}_{\delta_x}\left(\cdot| \mathcal{E}^c\right)$ almost surely,
\begin{align}\label{Step31}
	&	 \limsup_{t\to\infty}  \frac{e^{\lambda_1 t/2} (\langle f,  X_t\rangle- E_t(f_{la}))}{\sqrt{2\log t}} \nonumber\\
	& \leq   \limsup_{t\to\infty}  \frac{e^{\lambda_1 t/2} (\langle f_{main},  X_t\rangle- E_t(f_{la})) }{\sqrt{2\log t}}+ \limsup_{t\to\infty}  \frac{e^{\lambda_1 t/2} \left|\langle f_{rest},  X_t\rangle\right|}{\sqrt{2\log t}} \nonumber\\
	& \leq \sqrt{\left(\sigma_{sm}^2(f_{main})+\sigma_{la}^2(f)\right) W_\infty}+C\sqrt{W_\infty} \Vert  f_{rest}\Vert_2 ,
\end{align}
and similarly
\begin{align}\label{Step32}
	&	\limsup_{t\to\infty}  \frac{e^{\lambda_1 t/2} (\langle f,  X_t\rangle- E_t(f_{la}))}{\sqrt{2\log t}} \nonumber\\
	& \geq \sqrt{\left(\sigma_{sm}^2(f_{main})+\sigma_{la}^2(f)\right) W_\infty}-C\sqrt{W_\infty} \Vert  f_{rest}\Vert_2 .
\end{align}
By the dominated convergence theorem,   as $N\to\infty$,
\[
\sigma_{sm}(f_{main})\to \sigma_{sm}(f),\quad \Vert  f_{rest}\Vert_2  \to 0.
\]
Therefore, letting $N\to \infty$ in \eqref{Step31} and \eqref{Step32},
we get
\begin{align}
	\limsup_{t\to\infty}  \frac{e^{\lambda_1 t/2} (\langle f,  X_t\rangle- E_t(f_{la}))}{\sqrt{2\log t}} = \sqrt{\left(\sigma_{sm}^2(f)+\sigma_{la}^2(f)\right) W_\infty},\quad \mathbb{P}_{\delta_x}\left(\cdot| \mathcal{E}^c\right)\textup{-a.s.}
\end{align}
The proof for the liminf is similar and we complete the proof of the theorem.

\hfill$\Box$

The rest of the subsection is devoted to the proof of Proposition \ref{prop1}.
Recall  that \eqref{Semigroup} holds by the definition of $\mathcal{T}$.
We will use a different discretization scheme.
 For any $n\in \mathbb{N}$, define
\begin{align}\label{Defs}
	t_n:= n^{1/10}.
\end{align}
The following lemma shows that $\langle T_{t_{n+1}-t}f, X_t\rangle \approx \langle f, X_t\rangle$ for any $t\in [t_n, t_{n+1})$ as $n\to\infty$.

\begin{lemma}\label{Tec-cor}
		Let $f\in {\cal T}$ with $\mathfrak{R}_{\gamma(f)}>0$.
		Then under $\mathbb{P}_{\delta_x}(\cdot| \mathcal{E}^c)$, almost surely,
		\begin{align}
			\lim_{n\to\infty} \sup_{t_n \leq t<t_{n+1}} e^{\lambda_1 t/2}\left|\langle T_{t_{n+1}-t}f- f, X_t \rangle\right|=0.
		\end{align}
\end{lemma}
\textbf{Proof: }
Set $h:= \mathcal{L}f$. By \eqref{Semigroup},  for $t_n\leq t<t_{n+1}$,
\begin{align}
\left| \langle T_{t_{n+1}-t}f- f, X_t\rangle \right|& =\left| \int_0^{t_{n+1}-t} \langle  T_s h, X_t\rangle \mathrm{d}s\right| =\left|  \int_0^{t_{n+1}-t} \mathbb{E}_{\delta_x}\left(\langle h, X_{t+s}\rangle \Big|\mathcal{F}_t\right) \mathrm{d}s\right|  \nonumber\\
&\leq   \int_{t_n}^{t_{n+1}} \left|  \mathbb{E}_{\delta_x}\left(\langle h, X_{s}\rangle \Big|\mathcal{F}_t\right) 1_{\{s\geq t\geq t_n\}}\right| \mathrm{d}s.
\end{align}
Since $\mathcal{M}_t^{(s)}:= \mathbb{E}_{\delta_x}\left(\langle h, X_{s}\rangle \Big|\mathcal{F}_t\right) $ is a martingale for $t\in [t_n, s]$, it follows from Jensen's inequality and the $L^2$-maximal inequality that
\begin{align}
	\mathbb{E}_{\delta_x}\left(\sup_{t_n\leq t<t_{n+1}}\left| \langle T_{t_{n+1}-t}f- f, X_t\rangle \right|^2\right)&\leq (t_{n+1}-t_n) \int_{t_n}^{t_{n+1}} \mathbb{E}_{\delta_x}\left(\sup_{t_n\leq t\leq s}   (\mathcal{M}_t^{(s)})^2 \right)\mathrm{d}s\nonumber\\
	&\leq 4(t_{n+1}-t_n) \int_{t_n}^{t_{n+1}} \mathbb{E}_{\delta_x}\left( (\mathcal{M}_s^{(s)})^2 \right)\mathrm{d}s\nonumber\\
	&= 4(t_{n+1}-t_n) \int_{t_n}^{t_{n+1}} \mathbb{E}_{\delta_x}\left(  \langle h, X_s\rangle^2 \right)\mathrm{d}s.
\end{align}
Combining Lemma \ref{lemma2} (1) and the fact that $\gamma(f) = \gamma(h)$, we finally conclude that
\begin{align}
  & \sum_{n> (10t_0)^{10}} e^{\lambda_1 t_n}	\mathbb{E}_{\delta_x}\left(\sup_{t_n\leq t<t_{n+1}}\left| \langle T_{t_{n+1}-t}f- f, X_t\rangle \right|^2\right) \nonumber\\
  &\lesssim_{h, t_0}(b_{t_0}^{1/2}(x)+b_{t_0}(x))\sum_{n> (10 t_0)^{10}} (t_{n+1}-t_n)^2\lesssim (b_{t_0}^{1/2}(x)+b_{t_0}(x))  \sum_{n> (10t_0)^{10}} n^{-9/5}<\infty,
\end{align}
which implies the desired result by Markov's inequality and the Borel-Cantelli lemma.

\hfill$\Box$

Define
\begin{align}\label{New-def-of-Y-V}
	J_{t}^f&:=T_{t}f\left(x\right)  - \langle f, X_{t}\rangle, \quad
R_{t}^f(x):=
\mathbb{E}_{\delta_x} \left(\left(J_t^f\right) ^2\right).
\end{align}
The following lemma is a modification of Lemma \ref{lemma5}. We give a rough bound
for the conditional variance of $\langle T_{s_n} f, X_{t_n}\rangle$,  where  either $s_n=0$ or $s_n =t_{n+1}-t_n$.
\begin{lemma}\label{lemma5-2}
	Let $f\in L^2(E,\mu)\cap L^4(E,\mu)$ with $\mathfrak{R}_{\gamma(f)}>0$.
	Assume either $s_n=0$ or $s_n=t_{n+1}-t_n$ for all $n\in \mathbb{N}$. Then it holds that
	\begin{align}\label{e:lbl5.2}
		\liminf_{n \to \infty} e^{\lambda_1 t_n }\textup{Var}_x\left[\langle T_{s_n}f, X_{t_n} \rangle \bigg| \mathcal{F}_{t_n/2}\right] \geq \langle f^2, \widehat{\phi}_1 \rangle_\mu W_\infty,\quad \mathbb{P}_{\delta_x}\mbox{-a.s.}
	\end{align}
	and
	\begin{align}\label{e:ubl5.2}
		\limsup_{n \to \infty} e^{\lambda_1 t_n }\textup{Var}_x\left[\langle T_{s_n}f, X_{t_n} \rangle \bigg| \mathcal{F}_{t_n/2}\right] \leq
		\sigma_{sm}^2(f)
		W_\infty,\quad \mathbb{P}_{\delta_x}\mbox{-a.s.}
	\end{align}
\end{lemma}
\textbf{Proof:}
Using conditional independence, we get
\begin{align}\label{step_30}
	e^{\lambda_1 t_n }\textup{Var}_x\left[\langle T_{s_n}f, X_{t_n} \rangle \bigg| \mathcal{F}_{t_n/2}\right] =  e^{\lambda_1 t_n} \sum_{i=1}^{M_{t_n/2}} R_{t_n/2}^{T_{s_n}f}\left(X_{t_n/2}(i)\right) = e^{\lambda_1 t_n} \langle R_{t_n/2}^{T_{s_n}f}, X_{t_n/2}\rangle.
\end{align}
It follows from \eqref{Second-moment} that
\begin{align}\label{step_38}
	R_t^{T_{s_n}f} = \int_0^t T_s\left[A^{(2)}\cdot(T_{t-s+s_n}f)^2 \right]\mathrm{d}s + T_t((T_{s_n}f)^2)-(T_{t+s_n} f)^2  \geq T_t((T_{s_n}f)^2)- (T_{t+s_n} f)^2.\qquad
\end{align}
Combining \eqref{T-t-upp} with the fact $\Vert T_s \Vert_2\leq e^{\Vert A^{(1)}\Vert_\infty s}$, we get
$\Vert (T_{s_n}f)^2\Vert_2\lesssim \Vert T_{s_n}(f^2)\Vert_2 \leq \Vert T_{s_n}\Vert \Vert f^2 \Vert_2 \lesssim \Vert  f\Vert_4^2.$ Therefore,
applying Lemma \ref{lemma3} (1) with $a=\frac{\lambda_1+\mathfrak{R}_2}{2}$,  we get for any $t>2t_0$,
\begin{align}\label{step_39}
	 \left|T_t\left( \widetilde{(T_{s_n}f)^2}\right)\right| & =\left|T_t\left( |T_{s_n}f|^2\right) -  e^{-\lambda_1 t} \langle |T_{s_n}f|^2, \widehat{\phi}\rangle_\mu \phi_1 \right|
	\nonumber\\
	& \lesssim_{t_0} e^{-at}\Vert (T_{s_n}f)^2 \Vert_2
	b_{t_0}^{1/2} \lesssim_{f, t_0} e^{-at} b_{t_0}^{1/2}.
\end{align}
Therefore,
combining \eqref{step_38},  \eqref{step_39} and Lemma \ref{lemma3} (2) for $T_{t+s_n} f$,
we see that there exists a constant $C(f)>0$ such that for any $t>2t_0$ and $x\in E$,
\begin{align}
	e^{\lambda_1 t}R_t^{T_{s_n}f} & \geq  e^{\lambda_1 t}T_t((T_{s_n}f)^2)- e^{\lambda_1 t}(T_{t+s_n} f)^2 \nonumber\\
	& \geq \langle (T_{s_n}f)^2, \widehat{\phi}_1 \rangle_\mu \phi_1 -e^{\lambda_1 t} \left|T_t\left( \widetilde{(T_{s_n}f)^2}\right)\right|  -e^{\lambda_1 t}(T_{t+s_n} f)^2  \nonumber\\
	&\geq 	 \langle (T_{s_n}f)^2, \widehat{\phi}_1 \rangle_\mu \phi_1  - C(f)\left( e^{(\lambda_1 -a)t}+(t+s_n)^{2\tau(f)}e^{-\left(2\mathfrak{R}_{\gamma(f)}- \lambda_1\right)(t+s_n)} \right)
	 b_{t_0}^{1/2},
\end{align}
which together with \eqref{step_30} implies that
\begin{align}\label{low-var}
	& e^{\lambda_1 t_n }\textup{Var}_x\left[\langle T_{s_n}f, X_{t_n} \rangle \bigg| \mathcal{F}_{t_n/2}\right]  \geq  \langle (T_{s_n}f)^2, \widehat{\phi}_1 \rangle_\mu W_{t_n/2} \nonumber\\
	&\qquad - C(f) \left(e^{(\lambda_1 -a)t_n/2}+(t_n/2+s_n)^{2\tau(f)}e^{-\left(2\mathfrak{R}_{\gamma(f)}- \lambda_1\right)(t_n/2+s_n)} \right)  e^{\lambda_1 t_n/2}\langle b_{t_0}^{1/2}, X_{t_n/2}\rangle.
\end{align}
Since
\begin{align}
	\sum_{n> (4t_0)^{10} } \left(e^{(\lambda_1 -a)t_n/2}+(t_n/2+s_n)^{2\tau(f)}e^{-\left(2\mathfrak{R}_{\gamma(f)}- \lambda_1\right)(t_n/2+s_n)} \right)  e^{\lambda_1 t_n/2} \mathbb{E}_{\delta_x} \left(\langle b_{t_0}^{1/2}, X_{t_n/2}\rangle\right)<\infty,
\end{align}
the last term on the right hand side of \eqref{low-var} converges to $0$ almost surely as $n\to \infty$.
Letting $n\to\infty$ in \eqref{low-var} yields \eqref{e:lbl5.2}.

For the upper bound, combining
\eqref{Second-moment} and Jensen's inequality, we get
\begin{align}
	e^{\lambda_1 t}R_t^{T_{s_n}f }(x)& \leq e^{\lambda_1 t} R_{t+s_n}^f(x)	
	\leq 	e^{\lambda_1 t}  \mathbb{E}_{\delta_x}\left(\langle f, X_{t+s_n} \rangle^2 \right)\nonumber\\
	&= e^{\lambda_1 t}\int_0^{t+s_n} T_{t+s_n-s}\left[A^{(2)}\cdot(T_{s}f)^2 \right](x)\mathrm{d}s + e^{\lambda_1 t} T_{t+s_n}(f^2)(x).
\end{align}
Using the fact $s_n\in [0,1]$ and
an argument similar to that  for \eqref{step_39},
we get that there exists a constant $C(f)$ such that for $a=\frac{\lambda_1+\mathfrak{R}_2}{2}$ and $t>2t_0$,
\begin{align}\label{Step23}
	e^{\lambda_1 t} T_{t+s_n}(f^2)(x) \leq e^{-\lambda_1 s_n}\langle f^2,\widehat{\phi}_1
	\rangle_\mu \phi_1 (x)
	+ C(f) e^{-(a-\lambda_1)t} b_{t_0}^{1/2}(x).
\end{align}
Thus for any $2t_0< N<t/2$, combing  Lemma \ref{lemma3} (2) and \cite[(2.25)]{RSZ2017}, we get
\begin{align}\label{Step24}
	&e^{\lambda_1 t}\!\int_N^{t+s_n} T_{t+s_n-s}\left[A^{(2)}\cdot(T_{s}f)^2 \right](x)\mathrm{d}s \nonumber\\
	& \lesssim_{f,t_0} e^{\lambda_1 t}\!\int_N^{t+s_n-2t_0} s^{2\tau(f)}e^{-2\mathfrak{R}_{\gamma(f)}s} T_{t+s_n-s}(b_{t_0}) (x)\mathrm{d}s + t^{2\tau(f)} e^{(\lambda_1-2\mathfrak{R}_\gamma(f))t}b_{t_0}^{1/2}(x)\nonumber\\
	&
	\lesssim_{f, t_0}  \left( \int_N^\infty s^{2\tau(f)} e^{-(2\mathfrak{R}_{\gamma(f)}-\lambda_1)s}\mathrm{d}s +t^{2\tau(f)} e^{(\lambda_1-2\mathfrak{R}_\gamma(f))t} \right)   b_{t_0}^{1/2}(x)\nonumber\\
	&\lesssim_{f,t_0}  N^{2\tau(f)} e^{(\lambda_1-2\mathfrak{R}_\gamma(f))N}   b_{t_0}^{1/2}(x).
\end{align}
By Lemma \ref{lemma3} (1),
there exists a constant $C>0$ independent of $N$ such that
\begin{align}\label{Step25}
	& e^{\lambda_1 t}\int_0^N T_{t+s_n-s}\left[A^{(2)}\cdot(T_{s}f)^2 \right](x)\mathrm{d}s \nonumber\\
	& \leq 	e^{\lambda_1 t}\int_0^N e^{- \lambda_1(t+s_n-s)} \langle A^{(2)}\cdot(T_{s}f)^2 ,\widehat{\phi}_1 \rangle_\mu  \mathrm{d}s  \phi_1(x) +  C e^{\lambda_1 t}\int_0^N e^{- a(t+s_n-s)} \Vert  (T_{s}f)^2 \Vert_2  \mathrm{d}s  b_{t_0}^{1/2}(x) \nonumber\\
	& \leq  \phi _1(x) \left( e^{-\lambda_1 s_n}
	\int_0^\infty
	 e^{ \lambda_1 s} \langle A^{(2)}\cdot(T_{s}f)^2 ,\widehat{\phi}_1 \rangle_\mu  \mathrm{d}s \right) \nonumber\\
	&\qquad +  C \left(e^{-a s_n}e^{(\lambda_1-a) t} \Vert  f \Vert_4^2 \int_0^N e^{ (a+2\Vert A^{(1)}\Vert_\infty)s}   \mathrm{d}s   \right)b_{t_0}^{1/2}(x).\qquad
\end{align}
Therefore, combining \eqref{Step23}, \eqref{Step24} and \eqref{Step25}, there exists a constant $C'(f)= C'(f, t_0)>0$ such that for all $t>4t_0, x\in E$ and $2t_0<N <t/2$,
\begin{align}
	& e^{\lambda_1 t}R_t^{T_{s_n}f}
     \nonumber\\  & \leq e^{-\lambda_1 s_n} \sigma_{sm}^2(f)\phi_1
        + C'(f)\bigg(  e^{-(a-\lambda_1)t} +  N^{2\tau(f)} e^{(\lambda_1-2\mathfrak{R}_\gamma(f))N}  +   e^{(\lambda_1-a) t} e^{(|a|+2\Vert A^{(1)}\Vert_\infty )N } \bigg)b_{t_0}^{1/2}.
\end{align}
Taking $N= \varepsilon t_n$ such that  $-(\lambda_1-a)t_n/2  = (|a|+2\Vert A^{(1)}\Vert_\infty )N$,  we conclude that
\begin{align}
	& e^{\lambda_1 t_n }\textup{Var}_x\left[\langle
	T_{s_n}f,
	X_{t_n} \rangle \bigg| \mathcal{F}_{t_n/2}\right] \leq
		e^{-\lambda_1 s_n} \sigma_{sm}^2(f)
	W_{t_n/2}\nonumber\\
	&\qquad   + C'(f)\bigg( 2 e^{-(a-\lambda_1)t_n/2} +  (\varepsilon t_n)^{2\tau(f)} e^{(\lambda_1-2\mathfrak{R}_\gamma(f))\varepsilon t_n}    \bigg) e^{\lambda_1 t_n/2} \langle b_{t_0}^{1/2}, X_{t_n/2}\rangle.
\end{align}
Similar to the argument in proof of the lower bound,
the last term  of inequality converges to $0$ almost surely as $n\to \infty$.
Letting $n\to\infty$ in above inequality, we get \eqref{e:lbl5.2}.

\hfill$\Box$

Under the Assumption {\bf(H4)}, we have the following useful
lemma whose proof is postponed to Section \ref{ss:5.3}.

\begin{lemma}\label{Useful-lemma}
	Suppose in addition that {\bf(H4)} holds. If $f$ satisfies
	$|f|\lesssim_f b_{4t_0}^{1/2}$	and  $\mathfrak{R}_{\gamma(f)}>0$,
	then
	\begin{align}
		e^{2\lambda_1 t}
		\mathbb{E}_{\delta_x}\left(\langle f, X_t\rangle^4\right)\lesssim_f b_{t_0}^{1/2}(x),\quad t>T_0:=164 t_0, x\in E.
	\end{align}
\end{lemma}

Recall the definition of $J_t^f$ in \eqref{New-def-of-Y-V}. Combining
Lemmas \ref{lemma2} (1), \ref{Useful-lemma} and
inequalities $x^3\lesssim x^2+ x^4$ (for $e^{\lambda_1 t/2}|J_t^f|)$ and  $\mathbb{E}(|X- \mathbb{E}X|^4)\lesssim \mathbb{E}\left(X^4\right)$ (for $X= \langle f, X_t\rangle$), it is easy to get
that, for any $t>T_0$ and $x\in E$,
\begin{align}
	e^{3\lambda_1 t/2} \mathbb{E}_{\delta_x} \left(|J_t^f|^3\right) =  \mathbb{E}_{\delta_x} \left(|e^{\lambda_1 t/2} J_t^f|^3\right) \lesssim  e^{\lambda_1 t} R_t^f(x)+  e^{2\lambda_1 t}\mathbb{E}_{\delta_x} \left(|J_t^f|^4\right) \lesssim_{f,t_0} b_{t_0}^{1/2}(x)+b_{t_0}(x).
\end{align}
By Jensen's inequality, we deduce from the inequality above that for all $t> T_0$,
\begin{align}\label{Useful-ineq}
	e^{3\lambda_1 t/2} \sup_{s\in [0,1]} \mathbb{E}_{\delta_x} \left( |J_t^{T_sf}|^3\right) \leq \sup_{s\in [0,1]} e^{3\lambda_1 t/2} \mathbb{E}_{\delta_x} \left( |J_{t+s}^{f}|^3\right) \lesssim_{f, t_0} b_{t_0}^{1/2}(x)+ b_{t_0}(x).
\end{align}

The following result is a modification of Lemmas \ref{lemma7} and \ref{lemma12}.
With the help of Lemma \ref{Useful-lemma}, we are ready to give an upper bound for the limsup of the discrete-time version of the quantity in Proposition \ref{prop1}, as stated in the following lemma.

\begin{lemma}\label{lemma12-2}
	Suppose in addition that {\bf(H4)} holds.	If $f$ satisfies $|f|\lesssim_f b_{4t_0}^{1/2}$	and  $\mathfrak{R}_{\gamma(f)}>0$, then
	\begin{align}
		\limsup_{n\to\infty} \frac{e^{\lambda_1 t_n/2} \left(\left| \langle f, X_{t_n}\rangle\right| +\left| \langle T_{t_{n+1}-t_n}f, X_{t_n}\rangle\right| \right)}{\sqrt{2\log (t_n)}} \leq  8\sqrt{\sigma_{sm}^2(f) W_\infty},\quad \mathbb{P}_{\delta_x}\left(\cdot| \mathcal{E}^c\right)
		\mbox{-a.s.}
	\end{align}
\end{lemma}
\textbf{Proof: }
In the following  $s_n= 0$ or $s_n= t_{n+1}- t_n$ for all $n\in \mathbb{N}$.
Define
\begin{equation*}
	\Delta_{n}^{T_{s_n}f }: = \sup_{y\in\mathbb{R}} \left\vert \mathbb{P}_{\delta_x}\left[\frac{ \langle T_{s_n}f, X_{t_n} \rangle-\langle T_{t_n/2 + s_n}f, X_{t_n/2} \rangle}{\sqrt{\textup{Var}_x\left[\langle T_{s_n}f, X_{t_n} \rangle \big| \mathcal{F}_{t_n/2}\right]}} \leq y \bigg\vert
	\mathcal{F}_{t_n/2}
	\right]   -\Phi(y) \right\vert.
\end{equation*}
We claim that, $\mathbb{P}_{\delta_x}$-almost surely,
\begin{equation}\label{step_23-2}
	1_{\mathcal{E}^c}\sum_{n\geq 0} \Delta_{n}^{T_{s_n}f }< \infty.
\end{equation}
 Indeed, combining
 the branching property and \eqref{Useful-ineq}, we get
\begin{align}
	&\mathbb{E}_{\delta_x}\bigg( \sum_{n  > (2T_0)^{10}} e^{3\lambda_1 t_n/2} \sum_{i=1}^{M_{t_n/2}}\mathbb{E}_{\delta_x} \left[\left|T_{t_n/2 +s_n}f\left(X_{t_n/2}(i)\right)- \langle T_{s_n} f, X_{t_n/2}^i\rangle \right|^3 \Big| \mathcal{F}_{t_n/2}\right]\bigg)  \nonumber\\
	&= \sum_{n  > (2T_0)^{10}}  e^{3\lambda_1 t_n/2}  \mathbb{E}_{\delta_x}\left( \langle \mathbb{E}_{\delta_\cdot}\left(\left|J_{t_n/2}^{T_{s_n}f}\right|^3\right),X_{t_n/2}\rangle \right) \nonumber\\
	& \lesssim_{f,t_0} \sum_{n  > (2T_0)^{10}}  e^{3\lambda_1 t_n/2} e^{-3\lambda_1 t_n /4} \mathbb{E}_{\delta_x}\left(\langle b_{t_0}^{1/2}+b_{t_0}, X_{t_n/2}\rangle\right) \lesssim_{f,t_0} \sum_{n  > (2T_0)^{10}} e^{\lambda_1 t_n/4}<\infty,
\end{align}
where in the second inequality we also used  Lemma \ref{lemma3} (2). Therefore, almost surely,
\begin{align}\label{step_8-2}
	\sum_{n  \geq 1} e^{3\lambda_1 t_n/2} \sum_{i=1}^{M_{t_n/2}}\mathbb{E}_{\delta_x} \left[\left|T_{t_n/2 +s_n}f\left(X_{t_n/2}(i)\right)- \langle T_{s_n} f, X_{t_n/2}^i\rangle \right|^3 \Big| \mathcal{F}_{t_n/2}\right]<\infty.
\end{align}
It is trivial that $\Delta_{n}^{T_{s_n}f} \leq 2$. Combining $\{M_{t_n/2} > 0 \}\in \mathcal{F}_{t_n/2}$ and Lemma \ref{lemma6}, we see that there exists a constant $C_1$ such that under $\P_{\delta_x},$ on the event $\{M_{t_n/2}>0 \}$,
\begin{align}\label{step_13-2}
	\Delta_{n}^{T_{s_n}f} &\leq C_1 \frac{\sum_{i=1}^{M_{t_n/2}}\mathbb{E}_{\delta_x} \left[\left|T_{t_n/2+s_n}f\left(X_{t_n/2}(i)\right)- \langle T_{s_n}f, X_{t_n/2}^i\rangle \right|^3 \Big| \mathcal{F}_{t_n/2}\right]  }{\sqrt{\left(\textup{Var}_x\left[\langle T_{s_n}f, X_{t_n} \rangle \big| \mathcal{F}_{t_n/2} \right]\right)^3}} .
\end{align}
Since $\mathcal{E}^c \subset \{M_{t_n/2}>0 \},$ we see that  \eqref{step_13-2} holds on
the event $\mathcal{E}^c$.
Now suppose $\Omega_0$ is an event with
$\P_{\delta_x}\left(\Omega_0\right) = 1$
such that, for any $\omega \in \Omega_0$, the conclusion
of Lemma \ref{lemma5-2}, \eqref{step_8-2} and \eqref{step_13-2} hold.
Then for $\omega \in \Omega_0\cap \mathcal{E}^c$,
there exists a large $N = N(\omega)$ such that for $n \geq N,$
$$ \textup{Var}_x\left[\langle T_{s_n}f, X_{t_n} \rangle \big| \mathcal{F}_{t_n/2}\right](\omega) \geq \frac{e^{-\lambda_1 t_n}}{2}\langle f^2, \widehat{\phi}_1 \rangle_\mu W_\infty(\omega)>0.$$
Together with \eqref{step_13-2}, we have that on $\Omega_0\cap\mathcal{E}^c$,
\begin{align}
	&\sum_{n\geq 0} \Delta_{n}^{T_{s_n}f} \leq 2\left(1+N\right) \nonumber\\
	&\qquad+  \frac{ C_1\sqrt{8}}{\sqrt{\left[\langle f^2, \widehat{\phi}_1 \rangle_\mu W_\infty\right]^3} }\sum_{n \geq N} e^{3\lambda_1 t_n/2} \sum_{i=1}^{M_{t_n/2}}\mathbb{E}_{\delta_x} \left[\left|T_{t_n/2+s_n}f\left(X_{t_n/2}(i)\right)- \langle T_{s_n}f, X_{t_n/2}^i\rangle \right|^3 \Big| \mathcal{F}_{t_n/2}\right] .
\end{align}
Combining \eqref{step_8-2} with the inequality above, we get \eqref{step_23-2}.

Combining Lemma \ref{lemma1} (with $B = \mathcal{E}^c$) and \eqref{step_23-2}, we get
\begin{equation*}
	\limsup_{n\to\infty} \frac{\langle T_{s_n}f, X_{t_n} \rangle-\langle T_{t_n/2+s_n}f, X_{t_n/2}\rangle}{\sqrt{2\log n \textup{Var}_x\left[\langle T_{s_n}f, X_{t_n} \rangle\big| \mathcal{F}_{t_n/2}\right]}} \leq  1,\quad \mathbb{P}_{\delta_x}\left(\cdot \vert \mathcal{E}^c\right)\textup{-a.s.}
\end{equation*}
Recall that $t_n = n^{1/10}$. It follows from Lemma \ref{lemma5-2} and $\sqrt{10}<4$ that
\begin{equation}
	\limsup_{n\to\infty} \frac{e^{\lambda_1 t_n/2}\left(\langle T_{s_n}f, X_{t_n} \rangle-\langle T_{t_n/2+s_n}f, X_{t_n/2}\rangle\right)}{\sqrt{2\log (t_n)}} \leq  4\sqrt{\sigma_{sm}^2(f)  W_\infty},\quad \mathbb{P}_{\delta_x}\left(\cdot \vert \mathcal{E}^c\right)\textup{-a.s.}
\end{equation}
Since $\mathfrak{R}_{\gamma(f)}>0$, by Lemma \ref{lemma3} (2),
we have,
\begin{align}
 & \sum_{n> (2  t_0)^{10}} e^{\lambda_1 t_n/2} \mathbb{E}_{\delta_x} \left(\left| \langle T_{t_n/2+s_n} f,X_{t_n/2}\rangle \right| \right) \nonumber\\
 &\lesssim_{f, t_0} \sum_{n> (2  t_0)^{10}}  (t_n/2 +s_n)^{\tau(f)}e^{ - \mathfrak{R}_{\gamma(f)}(t_n/2+s_n)}  e^{\lambda_1 t_n/2}  T_{t_n/2}(b_{t_0}^{1/2})(x) \nonumber\\
 & \lesssim_{f,t_0} b_{t_0}^{1/2}(x) \sum_{n> (2 t_0)^{10}} t_n^{\tau(f)}e^{ - \mathfrak{R}_{\gamma(f)}t_n/2} <\infty,
\end{align}
which implies that $e^{\lambda_1 t_n/2}  \left| \langle T_{t_n/2+s_n} f,X_{t_n/2}\rangle \right|\to 0$ almost surely as $n\to\infty$. Therefore,
\begin{equation}
	\limsup_{n\to\infty} \frac{e^{\lambda_1 t_n/2}\langle T_{s_n}f, X_{t_n} \rangle}{\sqrt{2\log (t_n)}} \leq  4\sqrt{\sigma_{sm}^2(f)  W_\infty},\quad \mathbb{P}_{\delta_x}\left(\cdot \vert \mathcal{E}^c\right)\textup{-a.s.}
\end{equation}
Repeating the argument above with $f$ replaced by $-f$, we arrive at the desired assertion.

\hfill$\Box$

Now we treat the continuous-time setting in the following lemma using an idea roughly similar to that used in Lemma \ref{lemma9}.
\begin{lemma}\label{lemma8}
	Suppose in addition that {\bf(H4)} holds.	If $f$ satisfies 	$|f|\lesssim_f b_{4t_0}^{1/2}$  and  $\mathfrak{R}_{\gamma(f)}>0$, then
	\begin{align}
		\limsup_{n\to\infty} \sup_{t\in [t_n, t_{n+1})} \frac{e^{\lambda_1 t_n/2}\left|\langle T_{t_{n+1}-t}f, X_t\rangle\right|}{\sqrt{2\log t_n}}\leq 18 \sqrt{ \sigma_{sm}^2(f) W_\infty},\quad \mathbb{P}_{\delta_x}\left(\cdot| \mathcal{E}^c\right)
				\mbox{-a.s.}
	\end{align}
\end{lemma}
\textbf{Proof: }
From Lemma \ref{lemma12-2}, we see that
\begin{align}\label{Goal-4}
	\limsup_{n\to\infty} \frac{ e^{\lambda_1 t_n/2} \left|\langle T_{t_{n+1}-t_n}f, X_{t_n}\rangle -\langle f, X_{t_{n+1}}\rangle \right|}{\sqrt{2\log (t_n)}} \leq  8\sqrt{\sigma_{sm}^2(f) W_\infty},\quad \mathbb{P}_{\delta_x}\left(\cdot| \mathcal{E}^c\right)\textup{-a.s.}
\end{align}
Define
\begin{align}
	\varepsilon_n(f):= 10\sqrt{2\sigma_{sm}^2(f) e^{-\lambda_1 t_n}\log (t_n) W_{t_n}}.
\end{align}
Set $\mathcal{G}_n=\mathcal{F}_{t_n} $ and $B_n:= \left\{\langle T_{t_{n+1}-t_n}f, X_{t_n}\rangle -\langle f, X_{t_{n+1}}\rangle> \varepsilon_n(f)\right\}$, then  $B_n\in \mathcal{G}_{n+1}$ for all $n$. From the second Borel-Cantelli lemma, we get that
\begin{align}
	& \left\{\langle T_{t_{n+1}-t_n}f, X_{t_n}\rangle -\langle f, X_{t_{n+1}}\rangle > \varepsilon_n(f), \mbox{ i.o. }\right\}\nonumber\\
	&=\left\{\sum_{n=1}^\infty \mathbb{P}_{\delta_x}\left( \langle T_{t_{n+1}-t_n}f, X_{t_n}\rangle -\langle f, X_{t_{n+1}}\rangle > \varepsilon_n(f) \big| \mathcal{F}_{t_n}\right)=\infty \right\}.
\end{align}
By \eqref{Goal-4}, on $\mathcal{E}^c$, $\mathbb{P}_{\delta_x}$-almost surely,
\begin{align}\label{Step9}
	\sum_{n=1}^\infty \mathbb{P}_{\delta_x}\left( \langle T_{t_{n+1}-t_n}f, X_{t_n}\rangle -\langle f, X_{t_{n+1}}\rangle > \varepsilon_n(f) \big| \mathcal{F}_{t_n}\right)<\infty.
\end{align}
Define
\begin{align*}
	Z_t(f)&:= \mathbb{E}_{\delta_x}\left[\left(\langle f, X_{t_{n+1}}\rangle  -\langle T_{t_{n+1}-t}f, X_{t}\rangle \right)^2 \Big| \mathcal{F}_t\right],
	\quad t\in[t_n, t_{n+1}),
	\nonumber\\
	B_n(f)& := \sup_{t\in[t_n, t_{n+1})} \left[  \langle T_{t_{n+1}-t_n}f, X_{t_n}\rangle  - \langle T_{t_{n+1}-t}f, X_{t}\rangle - \sqrt{2Z_t(f)}\right],
	\\	T_n(f)&:= \inf\left\{
	s\in [t_n, t_{n+1}):
	\langle T_{t_{n+1}-t_n}f, X_{t_n}\rangle  - \langle T_{t_{n+1}-s}f, X_{s}\rangle - \sqrt{2Z_s(f)} > \varepsilon_n(f) \right\}.
\end{align*}
Similar to \eqref{step_19} and \eqref{Markov-Ineq}, by the strong Markov property and Markov's inequality, we have
\begin{align}\label{step_19-2}
	&\mathbb{P}_{\delta_x}\left( \langle T_{t_{n+1}-t_n}f, X_{t_n}\rangle -\langle f, X_{t_{n+1}}\rangle > \varepsilon_n(f) \big| \mathcal{F}_{t_n}\right) \nonumber\\
	&\geq \mathbb{P}_{\delta_x}\left( \langle T_{t_{n+1}-t_n}f, X_{t_n}\rangle -\langle f, X_{t_{n+1}}\rangle > \varepsilon_n(f), T_n(f)< t_{n+1} \big| \mathcal{F}_{t_n}\right) \nonumber \\
	&\geq  \mathbb{P}_{\delta_x}\left(\langle T_{t_{n+1}-T_n(f)}f, X_{
	T_n(f)}\rangle -\langle f, X_{t_{n+1}}\rangle > -\sqrt{2Z_{T_n(f)}}, T_n(f) <t_{n+1}
	\big|\mathcal{F}_{t_n} \right)\nonumber \\
	&\geq \frac{1}{2} \mathbb{P}_{\delta_x}\left( T_n(f)< t_{n+1} \big|\mathcal{F}_{t_n} \right)= \frac{1}{2} \mathbb{P}_{\delta_x}\left( B_n(f)> \varepsilon_n(f) \big|\mathcal{F}_{t_n} \right),
\end{align}
where the second inequality follows  from an argument similar to that leading to  \eqref{step_19},
and the last inequality follows  from an argument similar to that leading to \eqref{Markov-Ineq}.
Combining \eqref{Step9} and \eqref{step_19-2}, we get that $\mathbb{P}_{\delta_x}$-almost surely on $\mathcal{E}^c$,
$$
\sum_{n=1}^\infty \mathbb{P}_{\delta_x}\left( B_n(f)> \varepsilon_n(f) \big|\mathcal{F}_{t_n} \right) < +\infty.
$$
Applying again the second Borel-Cantelli lemma, we get that  $\mathbb{P}_{\delta_x}\left(\cdot | \mathcal{E}^c\right)$-almost surely,
\begin{align}\label{Step10}
	& \limsup_{n\to\infty} \sup_{t\in [t_n, t_{n+1})} \frac{e^{\lambda_1 t_n/2}\left(\langle T_{t_{n+1}-t_n}f, X_{t_n}\rangle  - \langle T_{t_{n+1}-t}f, X_{t}\rangle  \right)}{\sqrt{2\log (t_n)}}\nonumber\\
	&\leq
	\limsup_{n\to\infty}	\sup_{t\in [t_n, t_{n+1})}	\frac{\sqrt{e^{\lambda_1 t_n}Z_t(f)}}{\sqrt{2 \log (t_n)}}+ 10\sqrt{\sigma_{sm}^2(f)W_\infty}.
\end{align}
It follows from Lemma \ref{lemma3} (2) and Lemma
\ref{Useful-known-result} (3)
that $f^2\lesssim_f b_{4t_0}\lesssim_{f, t_0} T_{3t_0}(a_{t_0})\lesssim b_{t_0}^{1/2}$.
Therefore, using
 the inequality $\mbox{Var}(Y^2)\leq \mathbb{E}(Y^2)$, the branching property and
 Lemma \ref{Useful-known-result} (1), we obtain
\begin{align}\label{Step11}
	 e^{\lambda_1 t} Z_t(f)& \leq  e^{\lambda_1 t} \langle \mathbb{E}_{\delta_{\cdot}} \left(\langle f, X_{t_{n+1}-t}\rangle^2\right) , X_{t}\rangle \lesssim  e^{\lambda_1 t}\langle  T_{t_{n+1}-t}(f^2) , X_t\rangle \nonumber\\
	 & = e^{\lambda_1 t}\mathbb{E}_{\delta_x}\left(\langle  f^2 , X_{t_{n+1}}\rangle\Big|\mathcal{F}_t\right)
	 \lesssim_{f, t_0} e^{\lambda_1 t}\mathbb{E}_{\delta_x}\left(\langle  b_{t_0}^{1/2} , X_{t_{n+1}}\rangle\Big|\mathcal{F}_t\right).
\end{align}
 Noticing that $b_{t_0}^{1/2}- (b_{t_0}^{1/2})_{sm}$ is of form \eqref{e:rsspecial}, combining Theorems \ref{thm2} and \ref{thm4}, we get
$$
\lim_{n\to\infty} \sup_{t\in [t_n, t_{n+1})} e^{\lambda_1 t}\langle  T_{t_{n+1}-t}
(b_{t_0}^{1/2}- (b_{t_0}^{1/2})_{sm}), X_t\rangle = \langle
b_{t_0}^{1/2}, \widehat{\phi}_1\rangle_\mu W_\infty.
$$
Since $(b_{t_0}^{1/2})_{sm}
\in L^2(E,\mu)\cap L^4(E,\mu)$,  combining the $L^2$-maximal inequality and Lemma \ref{lemma2} (1),
\begin{align}
	& \sum_{n> (10 t_0)^{10}}e^{2\lambda_1 t_n} \mathbb{E}_{\delta_x}\Big( \Big| \sup_{t\in [t_n, t_{n+1})} \mathbb{E}_{\delta_x}\left(\langle
	(b_{t_0}^{1/2})_{sm},
	X_{t_{n+1}}\rangle\Big|\mathcal{F}_t\right)\Big| ^2 \Big)\nonumber\\
	&\leq 4\sum_{n> (10 t_0)^{10}}e^{2\lambda_1 t_n} \mathbb{E}_{\delta_x}\left(  \langle
	(b_{t_0}^{1/2})_{sm},
	X_{t_{n+1}}\rangle^2\right)\lesssim_{f, t_0} (b_{t_0}^{1/2}(x)+b_{t_0}(x)) \sum_{n> (10 t_0)^{10}}e^{\lambda_1 t_n} <\infty,
\end{align}
which implies that $\lim_{n\to\infty} \sup_{t\in [t_n, t_{n+1})} e^{\lambda_1 t}\langle  T_{t_{n+1}-t}
(b_{t_0}^{1/2})_{sm}, X_t\rangle = 0$ almost surely.
Combining Lemma \ref{lemma12-2}, \eqref{Step10} and the above arguments, we conclude that
\begin{align}
	-\liminf_{n\to\infty} \inf_{t\in [t_n, t_{n+1})} \frac{
		e^{\lambda_1 t_n/2} \langle T_{t_{n+1}-t}f, X_t\rangle}{\sqrt{2\log t_n}} \leq 18 \sqrt{\sigma_{sm}^2(f) W_\infty},\quad \mathbb{P}_{\delta_x}\left(\cdot| \mathcal{E}^c\right)\textup{-a.s.}
\end{align}
Using a similar argument with $f$ replaced by $-f$,  we  complete the proof of the lemma.

\hfill$\Box$

\textbf{Proof of Proposition \ref{prop1}: }
Proposition \ref{prop1} follows from Lemmas \ref{Tec-cor} and   \ref{lemma8}.

\hfill$\Box$

\subsection{Proof of Lemma \ref{Useful-lemma}}\label{ss:5.3}

\textbf{Proof of Lemma \ref{Useful-lemma}:}
Set $h:= \kappa |f| +f \geq 0$ with $\kappa\geq 1$.
Define $T_t^{(k)}h:= \mathbb{E}_{\delta_x}\left(\langle h, X_t\rangle^k\right)$.
By {\bf(H4)}{\bf(b)}, for $z\in [0, 1]$ and $k=1, 2, 3, 4$, $\partial^k_z\psi(\cdot, z)$ are bounded.
 We know that $T_t^{(1)}h=T_th$ and $T^{(2)}_th=\mathbb{E}_{\delta_x}\left(\langle h, X_t\rangle^2\right)$ is given by \eqref{Second-moment}.
 Recall the definition of $A^{(k)}$ in \eqref{def-A}. 
We now derive some formulas for $T_t^{(3)}h$ and  $T_t^{(4)}h$. 
We claim that 
\begin{align}\label{Step33}
	&T_t^{(3)}h =  \int_0^t T_{t-s} \left(A^{(3)}\cdot	\left(T_{s}h \right)^3 \right)\mathrm{d}s  +3\int_0^tT_{t-s}\left( A^{(2)}\cdot \left(T_{s}^{(2)}h	\right)T_{s}h \right)\mathrm{d}s+ T_t(h^3)
\end{align}
and that
\begin{align}\label{Step34}
	T_t^{(4)} h  =& \int_0^tT_{t-s} \left(A^{(4)}\cdot	\left(	T_{s}h \right)^4\right) \mathrm{d}s  + 6\int_0^t T_{t-s}\left( A^{(3)}\cdot\left(	T_{s}h \right)^2 	T_{s}^{(2)}h \right)\mathrm{d}s  \nonumber\\
& +4\int_0^t T_{t-s}\left( A^{(2)}\cdot	T_{s}^{(3)}h	T_{s}h \right) \mathrm{d}s+ 3\int_0^t T_{t-s}\left(A^{(2)} \cdot\left(	T_{s}^{(2)}h \right)^2	\right)\mathrm{d}s+ T_t(h^4).
\end{align}
In fact, if $h\geq 0$ is bounded, the above results were proved in \cite{GHK2022} for the case $\lambda_1=0$, but the argument there works for the case 
$\lambda_1\neq 0$. 
Thus, using a routine limit argument, one can check that \eqref{Step33} and \eqref{Step34} also holds any $h\geq 0$.
Since $|f|\lesssim_f b_{4t_0}^{1/2}$, combining Lemma \ref{lemma3} (2) and
 Lemma \ref{Useful-known-result} (3)  (with $t=t_0$), we see that
\begin{align}\label{useful-ineq}
|f|^2 \lesssim_f b_{4t_0} \lesssim_{t_0} T_{3t_0}(a_{t_0})\lesssim_{t_0} b_{t_0}^{1/2} \quad \Rightarrow\quad |f|^3 \lesssim_f b_{4t_0}^{1/2}b_{t_0}^{1/2} \in L^2 (E,\mu),
\end{align}
which implies  $T_t(|f|^3) \lesssim_{f,t_0}e^{-\lambda_1 t} b_{t_0}^{1/2}$ by  Lemma \ref{lemma3} (2).
Then  $T_t(h^3), T_t(h^4)<\infty$,
and  all integrals on the right side of  \eqref{Step33} and \eqref{Step34}  are finite.

For any $x$ and $t$, both sides \eqref{Step33} and both sides of \eqref{Step34} are polynomials of $\kappa$.
Since \eqref{Step33} is  valid for all $\kappa\geq1$, the corresponding coefficients of the polynomials
on both sides agree.  The same is valid for \eqref{Step34}.
Thus,
\begin{align}\label{Step35}
	&T_t^{(3)}f =  \int_0^t T_{t-s} \left(A^{(3)}\cdot	\left(T_{s}f \right)^3 \right)\mathrm{d}s  +3\int_0^tT_{t-s}\left( A^{(2)}\cdot \left(T_{s}^{(2)}f	\right)T_{s}f \right)\mathrm{d}s+ T_t(f^3)
\end{align}
and
\begin{align}\label{Step36}
	T_t^{(4)} f  =& \int_0^tT_{t-s} \left(A^{(4)}\cdot	\left(	T_{s}f \right)^4\right) \mathrm{d}s  + 6\int_0^t T_{t-s}\left( A^{(3)}\cdot\left(	T_{s}f\right)^2 	T_{s}^{(2)}f \right)\mathrm{d}s  \nonumber\\
	& +4\int_0^t T_{t-s}\left( A^{(2)}\cdot	T_{s}^{(3)}f	T_{s}f \right) \mathrm{d}s+ 3\int_0^t T_{t-s}\left(A^{(2)} \cdot\left(	T_{s}^{(2)}f \right)^2	\right)\mathrm{d}s+ T_t(f^4).
\end{align}

$T_t^{(1)}f =T_t f$ can be bounded from above by using Lemma \ref{lemma3} (2), so we treat $T_t^{(2)} f$ first.
It was proved in \cite[(2.22) and (2.24)]{RSZ2017} (with $t_0$ replaced by $4t_0$ ) that for any $t> 40 t_0$,
\begin{align}
	|T_t^{(2)} f| \lesssim_{f, t_0} (e^{-\lambda_1 t} + t^{2\tau(f)} e^{-2\mathfrak{R}_{\gamma(f)}t})b_{4t_0}^{1/2}+  \int_{8t_0}^{t-8t_0}T_s\left[
	|T_{t-s}f|^2 \right]\mathrm{d}s + T_t(|f|^2).
\end{align}
Applying Lemma \ref{lemma3} (2) (with $t_1=4t_0$) repeatedly
and noticing that $2\mathfrak{R}_{\gamma(f)}>0>\lambda_1$,  for $t>40t_0$,
\begin{align}\label{fine-upp-T-2}
		|T_t^{(2)} f| & \lesssim_{f, t_0} e^{-\lambda_1 t} b_{4t_0}^{1/2}+  \int_{8t_0}^{t-8t_0} s^{2\tau(f)} e^{-2\mathfrak{R}_{\gamma(f)}s}T_{t-s}\left[
 b_{4t_0}\right]\mathrm{d}s \nonumber\\
 &\lesssim_{f, t_0} e^{-\lambda_1 t} b_{4t_0}^{1/2} \left( 1+ \int_{8t_0}^{t-8t_0} s^{2\tau(f)} e^{-2\mathfrak{R}_{\gamma(f)}s} e^{\lambda_1 s} \mathrm{d }s \right) \lesssim_{f, t_0} e^{-\lambda_1 t} b_{4t_0}^{1/2}.
\end{align}

 Now we treat $T_t^{(3)}f$.
For $t>20t_0$, by Lemma \ref{lemma3} (2) (with $t_1=4t_0$) and
\eqref{T-t-upp},
\begin{align}
|T_t f|^2 & = |T_{12 t_0}(T_{t-12t_0}f)|^2 \lesssim_{t_0}  |T_t f|^2  \land T_{12 t_0}(|T_{t-12t_0}f|^2)  \nonumber\\
& \lesssim_{f, t_0} t^{2\tau(f)}e^{-2\mathfrak{R}_{\gamma(f)} t} \left(T_{12 t_0}(b_{4t_0}) \land b_{4t_0} \right)\lesssim t^{2\tau(f)}e^{-2\mathfrak{R}_{\gamma(f)}t} ( b_{4t_0}^{1/2}\land b_{4t_0}).
\end{align}
Therefore, for $t>41 t_0$,
\begin{align}\label{fine-upp-T-1-2}
	|T_t f|^4  & = 	|T_{24t_0} (T_{t-24t_0}f)|^4 \lesssim_{t_0}	|T_{24t_0} (|T_{t-24t_0}f|^2)|^2 \lesssim t^{4\tau(f)}e^{-4\mathfrak{R}_{\gamma(f)}t} |T_{24t_0} (b_{4t_0})|^2 \nonumber\\
	& = t^{4\tau(f)}e^{-4\mathfrak{R}_{\gamma(f)}t} |T_{12t_0} T_{12t_0}(b_{4t_0})|^2  \lesssim_{t_0} t^{4\tau(f)}e^{-4\mathfrak{R}_{\gamma(f)}t} |T_{12t_0} (b_{4t_0}^{1/2})|^2 \nonumber\\
	&\lesssim_{t_0} t^{4\tau(f)}e^{-4\mathfrak{R}_{\gamma(f)}t} \left(T_{12t_0} (b_{4t_0}) \land b_{4t_0}\right) \lesssim_{t_0}t^{4\tau(f)}e^{-4\mathfrak{R}_{\gamma(f)}t}(b_{4t_0}^{1/2}\land b_{4t_0}).
\end{align}
Combining the two inequalities above, we conclude that for $t>41t_0$,
\begin{align}\label{fine-upp-T-1-1}
	|T_tf|^3  = \sqrt{|T_tf|^2 |T_tf|^4}\lesssim_{f,t_0} t^{3\tau(f)}e^{-3\mathfrak{R}_{\gamma(f)}t}(b_{4t_0}^{1/2}\land b_{4t_0}) .
\end{align}
Since $|ab| \leq \frac{2}{3}|a|^{3/2} + \frac{1}{3} |b|^{3}$ and $|T_s f|^{p}\lesssim_{p, t_0}T_s(|f|^p)$ for any $s\leq 41 t_0$ and $p>1$, we have, for $t>41t_0$,
\begin{align}
	& \left| \int_0^{41t_0} T_{t-s} \left(A^{(3)}\cdot	\left(T_{s}f \right)^3 \right)\mathrm{d}s  +3\int_0^{41t_0}T_{t-s}\left( A^{(2)}\cdot \left(T_{s}^{(2)}f	\right)T_{s}f \right)\mathrm{d}s+ T_t(f^3)\right| \nonumber\\
	&\lesssim_{t_0}\int_0^{41t_0} T_{t-s} \left(	\left|T_{s}f \right|^3 \right)\mathrm{d}s  +\int_0^{41t_0}T_{t-s}\left(  \left|T_{s}^{(2)}f	\right| ^{3/2}\right)\mathrm{d}s+ T_t(|f|^3)\nonumber\\
	&\lesssim  \int_0^{41t_0} T_{t-s} \left(T_{s}(|f|^3) \right)\mathrm{d}s  +\int_0^{41t_0}T_{t-s}\left(  \left|T_s(|f|^2)\right| ^{3/2}\right)\mathrm{d}s+ T_t(|f|^3) \lesssim_{t_0} T_t(|f|^3).
\end{align}
We remark here that according to the same argument, we also have for each $t>0$,
\begin{align}\label{Rough-T-3}
	|T_t^{(3)}f| \lesssim_t T_t(|f|^3).
\end{align}
Therefore, combining \eqref{Step35}, \eqref{fine-upp-T-2} and \eqref{fine-upp-T-1-1}, we see that for $t> 41 t_0$,
\begin{align}
	|T_t^{(3)} f| & \lesssim_{f, t_0} T_t(|f|^3)+ \int_{41 t_0}^t s^{3\tau(f)}e^{-3\mathfrak{R}_{\gamma(f)}s} T_{t-s} \left(b_{4t_0} \right)\mathrm{d}s  +\int_{41t_0}^t s^{\tau(f)} e^{-(\lambda_1+\mathfrak{R}_{\gamma(f)}) s}T_{t-s}\left( b_{4t_0} \right)\mathrm{d}s\nonumber\\
	& \lesssim_{f,t_0} T_t(|f|^3) + \int_{41t_0}^t s^{\tau(f)} e^{-(\lambda_1+\mathfrak{R}_{\gamma(f)}) s}T_{t-s}\left( b_{4t_0} \right)\mathrm{d}s,
\end{align}
where in the last inequality  we used the fact that $2\mathfrak{R}_{\gamma(f)}>\lambda_1$.
 Then using  Lemma
\ref{Useful-known-result} (3) (with $t=t_0$),
  for $t>41 t_0$,
\begin{align}
		|T_t^{(3)} f| &  \lesssim_{f,t_0}  e^{-\lambda_1 t} b_{t_0}^{1/2} + \int_{41t_0}^t s^{\tau(f)} e^{-(\lambda_1+\mathfrak{R}_{\gamma(f)}) s}T_{t+2t_0-s}\left( a_{2t_0}	\right)\mathrm{d}s \nonumber\\
		&\lesssim e^{-\lambda_1 t} b_{t_0}^{1/2} + e^{-\lambda_1 t} b_{t_0}^{1/2} \int_{41t_0}^t s^{\tau(f)} e^{-\mathfrak{R}_{\gamma(f)} s}\mathrm{d}s\lesssim_{f, t_0} e^{-\lambda_1 t} b_{t_0}^{1/2}.
\end{align}
With $t_0$ replaced by $4t_0$, we obtain that for $t>164 t_0$,
\begin{align}\label{fine-upp-T-3}
		|T_t^{(3)} f| \lesssim_{f, t_0} e^{-\lambda_1 t} b_{4t_0}^{1/2}.
\end{align}

Finally we bound $T_t^{(4)}f$ from above.
Combining \eqref{Step36} and inequalities $|a|^2 |b|\lesssim |a|^4 + |b|^2$ and $|a| |b|\lesssim |a|^{4/3} +|b|^4$ we obtain that, for $t>164t_0$,
\begin{align}\label{rou-upp-T-4}
	\left|T_t^{(4)} f \right|  \lesssim &\int_0^tT_{t-s} \left(	\left|	T_{s}f \right|^4\right) \mathrm{d}s  + \int_0^t T_{t-s}\left( \left|	T_{s}f\right|^2 	|T_{s}^{(2)}f| \right)\mathrm{d}s  \nonumber\\
	&\quad +\int_0^t T_{t-s}\left( |	T_{s}^{(3)}f| \cdot |	T_{s}f| \right) \mathrm{d}s+ \int_0^t T_{t-s}\left(\left|	T_{s}^{(2)}f \right|^2	\right)\mathrm{d}s+ T_t(f^4)\nonumber\\
	\lesssim  &\int_0^tT_{t-s} \left(	\left|	T_{s}f \right|^4\right) \mathrm{d}s    +\int_0^{164t_0} T_{t-s}\left( |	T_{s}^{(3)}f|^{4/3}  \right) \mathrm{d}s \nonumber\\
	&\quad +  \int_{164t_0}^t T_{t-s}\left( |	T_{s}^{(3)}f| \cdot |	T_{s}f| \right) \mathrm{d}s + \int_0^t T_{t-s}\left(\left|	T_{s}^{(2)}f \right|^2	\right)\mathrm{d}s+ T_t(f^4).
\end{align}
For $t>164 t_0$,  combining \eqref{Rough-T-3}, $|T_s f|^{p}\lesssim_{p, t_0}T_s(|f|^p)$ for any $s\leq 164 t_0$ and $p>1$, and $|T_s^{(2)}f|\lesssim_{t_0} T_s (f^2)$ for $s\leq 164 t_0$, it holds that
\begin{align}
	&\int_0^{164t_0}T_{t-s} \left(	\left|	T_{s}f \right|^4\right) \mathrm{d}s    +\int_0^{164t_0} T_{t-s}\left( |	T_{s}^{(3)}f|^{4/3}  \right) \mathrm{d}s+ \int_0^{164t_0} T_{t-s}\left(\left|	T_{s}^{(2)}f \right|^2	\right)\mathrm{d}s+ T_t(f^4)\nonumber\\
	&\lesssim_{t_0} \int_0^{164t_0}T_{t-s} \left(	T_{s}(f^4)\right) \mathrm{d}s    +\int_0^{164t_0} T_{t-s}\left( |	T_{s}(|f|^3)|^{4/3}  \right) \mathrm{d}s \nonumber\\
	&\qquad + \int_0^{164t_0} T_{t-s}\left(\left|	T_s(|f|^2) \right|^2	\right)\mathrm{d}s+ T_t(f^4) \lesssim_{t_0}T_t(f^4).
\end{align}
Therefore, combining \eqref{fine-upp-T-2},  \eqref{fine-upp-T-1-2},  \eqref{fine-upp-T-3} and \eqref{rou-upp-T-4}, we obtain that
\begin{align}
& 	|T_t^{(4)}f| \nonumber\\
& \lesssim_{t_0} T_t(f^4)+ \int_{164t_0} ^tT_{t-s} \left(	\left|	T_{s}f \right|^4\right) \mathrm{d}s    +\int_{164t_0}^t T_{t-s}\left( |	T_{s}^{(3)}f| \cdot |T_sf| \right) \mathrm{d}s+ \int_{164t_0}^t T_{t-s}\left(\left|	T_{s}^{(2)}f \right|^2	\right)\mathrm{d}s\nonumber\\
&\lesssim_{f,t_0} T_t(f^4)+ \int_{164t_0} ^t s^{4\tau(f)}e^{-4\mathfrak{R}_{\gamma(f)}s} T_{t-s} \left( b_{4t_0}\right) \mathrm{d}s    +\int_{164t_0}^t s^{\tau(f)}e^{-\mathfrak{R}_{\gamma(f)}s}e^{-\lambda_1 s} T_{t-s}\left( b_{4t_0}  \right) \mathrm{d}s\nonumber\\
&\quad + \int_{164t_0}^t e^{-2\lambda_1 s}T_{t-s}\left(b_{4t_0}	\right)\mathrm{d}s \lesssim T_t(f^4) + \int_{164t_0}^t e^{-2\lambda_1 s}T_{t-s}\left(b_{4t_0}	\right)\mathrm{d}s .
\end{align}
Since $b_{4t_0}\lesssim_{t_0} T_{2t_0}
(a_{2t_0})$ and that $|f|^4\lesssim_f b_{t_0}\in L^2(E,\mu)$ by \eqref{useful-ineq}, we deduce that for all $t>164 t_0$,
\begin{align}
	 	|T_t^{(4)}f| \lesssim_{f,t_0} &T_t(f^4) + \int_{164t_0}^t e^{-2\lambda_1 s}T_{t+2t_0-s}\left( a_{2t_0}		\right)\mathrm{d}s \nonumber\\
	\lesssim_{f, t_0} &  e^{-\lambda_1 t} b_{t_0}^{1/2} + e^{-\lambda_1 t} b_{t_0}^{1/2} \int_{164t_0}^t e^{-\lambda_1 s}\mathrm{d}s \lesssim_{f,t_0}e^{-2\lambda_1 t} b_{t_0}^{1/2},
\end{align}
as desired.

\hfill$\Box$

\bigskip
\noindent
{\bf Acknowledgements:}
We thank the referees for the many comments and suggestions that have
greatly helped us to improve this paper.
	
	\bigskip
	\noindent

	\begin{singlespace}
		
	\end{singlespace}

	\vskip 0.2truein
	\vskip 0.2truein

\noindent{\bf Haojie Hou:}  School of Mathematics and Statistics, Beijing Institute of Technology, Beijing 100081, P. R. China.  Email: {\texttt
houhaojie@bit.edu.cn}

\smallskip		
	\noindent{\bf Yan-Xia Ren:} LMAM School of Mathematical Sciences \& Center for
	Statistical Science, Peking
	University,  Beijing, 100871, P.R. China. Email: {\texttt
	yxren@math.pku.edu.cn}
	
	\smallskip
	\noindent {\bf Renming Song:} Department of Mathematics,
	University of Illinois Urbana-Champaign,
	Urbana, IL 61801, U.S.A.
	Email: {\texttt rsong@illinois.edu}

\end{document}